%% file: tail-bound-MAIN-v9_BUMI.tex
\documentclass[12pt,letter]{article}
\usepackage[nottoc]{tocbibind}
\usepackage[latin1]{inputenc}
\usepackage{fancyhdr}
\usepackage{verbatim}
\usepackage{indentfirst}
\usepackage{graphicx}
\usepackage{epstopdf}
\usepackage{color}
\usepackage{newlfont}
\usepackage{amssymb}
\usepackage[intlimits]{amsmath}
\usepackage{latexsym}
\usepackage{amsthm}
\usepackage{amscd}
\usepackage{mathrsfs}
\usepackage[dvips, outer=1in, top=1in, bottom=1in]{geometry}
\usepackage{multirow}
\usepackage{hyperref}
\usepackage{rotating}
\usepackage{mathtools} 
\usepackage{dsfont}

\theoremstyle{plain}                    
\newtheorem{theorem}{Theorem}[section]
\newtheorem*{theorem*}{Theorem}
\newtheorem{lemma}[theorem]{Lemma}    
\newtheorem{prop}[theorem]{Proposition}
\newtheorem{cor}[theorem]{Corollary}
\theoremstyle{definition}
\newtheorem{defin}[theorem]{Definition}
        
\newtheorem{remark}[theorem]{Remark}

\newcommand{\N}{\mathbb{N}}
\newcommand{\Z}{\mathbb{Z}}
\newcommand{\Q}{\mathbb{Q}}
\newcommand{\R}{\mathbb{R}}
\newcommand{\C}{\mathbb{C}}
\newcommand{\A}{\mathcal{A}}

\newcommand{\h}{\mathfrak{H}}
\newcommand{\sltr}{\mathrm{SL}(2,\R)}
\newcommand{\sltz}{\mathrm{SL}(2,\Z)}
\newcommand{\asltr}{\mathrm{ASL}(2,\R)}

\newcommand{\tsltr}{\widetilde{\mathrm{SL}}(2,\R)}
\newcommand{\ha}{\frac{1}{2}}
\newcommand{\tha}{\tfrac{1}{2}}

\newcommand{\wtPhi}{\widetilde{\Phi}}
\newcommand{\calL}{\mathcal{L}}
\newcommand{\calR}{\mathcal{R}}
\newcommand{\calLL}{\calL \overline{\calL}}
\newcommand{\calLR}{\calL \overline{\calR}}
\newcommand{\calRL}{\calR \overline{\calL}}
\newcommand{\calRR}{\calR \overline{\calR}}

\newcommand{\de}{\mathrm{d}}

\newcommand{\Hei}{\mathbb{H}(\mathbb{R})}
\newcommand{\ba}{\begin{array}}
\newcommand{\ea}{\end{array}}
\newcommand{\ve}[2]{\left(\ba{c}\!#1\!\\ \!#2\!\ea\right)}
\newcommand{\sve}[2]{\left(\begin{smallmatrix}\!#1\!\\ \!#2\!\end{smallmatrix}\right)}

\newcommand{\sveH}[3]{\left(\left(\begin{smallmatrix}\!#1\!\\ \!#2\!\end{smallmatrix}\right),#3\right)}
\newcommand{\sgn}{\mathrm{sgn}}

\newcommand{\bm}[1]{\mbox{\boldmath{$#1$}}}
\newcommand{\LtR}{\mathrm L^2(\R)}
\newcommand{\SR}{\mathcal{S}(\R)}
\newcommand{\ma}[4]{
\begin{pmatrix}#1&#2\\#3&#4\end{pmatrix}
}
\newcommand{\sma}[4]{\left(\begin{smallmatrix} #1&#2\\#3&#4\end{smallmatrix}\right)}
\newcommand{\e}[1]{\mathrm{e}\!\left(#1\right)}

\newcommand{\Leb}{\mathrm{Leb}}
\newcommand{\Thetapair}[2]{\Theta_{#1}\overline{\Theta_{#2}}}

\def\veck{{\text{\boldmath$k$}}}

\def\vecxi{{\text{\boldmath$\xi$}}}

\def\Re{\operatorname{Re}}
\def\Im{\operatorname{Im}}

\def\GamG{\Gamma\backslash G}

\def\ASL{\ASL(2,\mathbb{R})}

\def\Onder#1#2#3#4#5{#1 \setbox0=\hbox{$#1$}\setbox1=\hbox{$#2$}
       \dimen0=.5\wd0 \dimen1=\dimen0 \dimen2=\dp0 \dimen3=\dimen2
       \advance\dimen0 by .5\wd1 \advance\dimen0 by -#4
       \advance\dimen1 by -.5\wd1 \advance\dimen1 by -#4
       \advance\dimen2 by -#3 \advance\dimen2 by \ht1
       \advance\dimen2 by 0.3ex \advance\dimen3 by #5
        \kern-\dimen0\raisebox{-\dimen2}[0ex][\dimen3]{\box1}
       \kern\dimen1}

\newcommand{\lcur}{\left\{}
\newcommand{\rcur}{\right\}}
\newcommand{\labs}{\left|}
\newcommand{\rabs}{\right|}
\newcommand{\lp}{\left(}
\newcommand{\rp}{\right)}

\newcommand{\ind}{\mathds{1}}

\newcommand{\GaG}{\Gamma\backslash G}

\newcommand{\Si}{\mathcal{S}}
\newcommand{\bn}{\mathbf{0}}

\newcommand{\tg}{\tilde{g}}

\numberwithin{equation}{section}

\input xy
\xyoption{all}
\usepackage{import}

\newcommand{\footremember}[2]{%
    \footnote{#2}
    \newcounter{#1}
    \setcounter{#1}{\value{footnote}}%
}

\title{Improved Tail Estimates for the Distribution of Quadratic {W}eyl Sums}
\author{
Francesco Cellarosi\footremember{queens}{Department of Mathematics and Statistics. Queen's University. Kingston, ON, Canada.}\footnote{Corresponding author: \texttt{fc19@queensu.ca}}
\and Jory Griffin\footremember{bristol}{School of Mathematics. University of Bristol. Bristol, U.K.}
\and Tariq Osman\footremember{brandeis}{Department of Mathematics. Brandeis University. Waltham, MA, U.S.A.}.
}
\date{}

\begin{document}
\clearpage{\pagestyle{empty}\cleardoublepage}


\maketitle

\begin{center}
\today

\end{center}

\begin{abstract}
We consider quadratic Weyl sums $S_N(x;c,\alpha)=\sum_{n=1}^N \exp\left\{2\pi i\left(  \left(\tfrac{1}{2}n^2+cn\right)\!x\right.\right.$ $\left.\left.+\alpha n\right)\right\}$   for $c=\alpha=0$ (the rational case) or $(c,\alpha)\notin\Q^2$ (the irrational case), where $x\in\R$  is randomly distributed according to a probability measure absolutely continuous with respect to the Lebesgue measure. 
The limiting distribution in the complex plane of $\frac{1}{\sqrt{N}}S_N(x;c,\alpha)$ as $N\to\infty$ was described in \cite{Marklof-1999} (respectively  \cite{Cellarosi-Marklof}) in the rational  (resp. irrational) case. 
According to the limiting distribution, the probability of landing outside a ball of radius $R$ is known to be  asymptotic  to $\frac{4\log 2}{\pi^2}R^{-4}(1+o(1))$ in the rational case and to $\frac{6}{\pi^2}R^{-6}(1+O(R^{-12/31}))$ in the irrational case, as $R\to\infty$. In this work we refine the technique of \cite{Cellarosi-Marklof} to improve the known tail estimates to $\frac{4\log 2}{\pi^2}R^{-4}(1+O_\varepsilon(R^{-2+\varepsilon}))$ and  $\frac{6}{\pi^2}R^{-6}(1+O_\varepsilon(R^{-2+\varepsilon}))$ for every $\varepsilon>0$. In the rational case, we rely on the equidistribution of a rational horocycle lift to a torus bundle over the unit tangent bundle to the classical modular surface.  All the constants implied by the $O_\varepsilon$-notations are made explicit.  
\end{abstract}

\tableofcontents
\newpage

\import{./}{introduction-02.tex}
\import{./}{preliminaries-05.tex}

\import{./}{limit-theorems-07.tex}

\import{./}{growth-in-the-cusp-09.tex}

\import{./}{partition-of-unity-04.tex}

\import{./}{stationary-phase-lemma-9.tex}

\import{./}{tail-approximations-10.tex}
\import{./}{appendix-zeta-01.tex}

\section*{Acknowledgements} The first and third authors acknowledge the support from the NSERC Discovery Grant ``Statistical and Number-Theoretical Aspects of Dynamical Systems''. The authors thank the anonymous referee for the  useful corrections they suggested. \\ 

\bibliographystyle{plain}
\bibliography{improved-tail-bound-bibliography}
\end{document}

%% file: introduction-02.tex
\section{Introduction}\label{section-introduction}

Set $\e{x}:=e^{2\pi i x}$ and consider quadratic Weyl sums of the form 
\begin{align}
S_N(x)=S_N(x;c,\alpha)=\sum_{n=1}^N \e{ \left(\tfrac{1}{2}n^2+cn\right)\!x+\alpha n},\label{theta-sum-intro}
\end{align}
 where $N$ is a positive integer, and $x$, $c$ and $\alpha$ are real. 
We assume that $c$ and $\alpha$ are fixed, and that $x$ in \eqref{theta-sum-intro} is randomly distributed on $\R$ according to a probability measure $\lambda$, absolutely continuous with respect to the Lebesgue measure on $\R$. 
Sums \eqref{theta-sum-intro}  can be thought of as the $N$-th partial sum of a particular sequence of strongly dependent random variables on the unit circle, or as the position after $N$ steps of a deterministic walk in $\C$ with a random seed $x$. It is natural to ask how the parameters $(c,\alpha)$ affect the distribution of  \eqref{theta-sum-intro}, after suitable rescaling, as $N\to\infty$.
 The limiting distribution on the complex plane of $\frac{S_N}{\sqrt{N}}$ as $N\to\infty$ has been studied by Marklof \cite{Marklof-1999} when $c=\alpha=0$ and by Cellarosi and Marklof \cite{Cellarosi-Marklof} when $(c,\alpha)\notin\Q^2$. We shall refer to the two cases as \emph{rational} and \emph{irrational}, respectively, see also Remark \ref{rk-general-rational-case} below.  In both cases the limiting distribution of of $\frac{S_N}{\sqrt{N}}$  is not Gaussian, as  can be seen by the following heavy tailed asymptotics.

\begin{itemize}
\item
Rational case: if $c=\alpha=0$, 
then
\begin{align}
\lim_{N\to\infty}\lambda\left\{x\in\R:\:\frac{|S_N(x;0,0)|}{\sqrt N}>R\right\}
=\displaystyle\frac{4\log2}{\pi^2}\frac{1}{R^4}\left(1+o(1)\right)\label{intro-tail-estimate-rat}
\end{align}
as $R\to\infty$,
see \cite{Marklof-1999}.
\item Irrational case: 
if $(c,\alpha)\notin\Q^2$, then 
\begin{align}
\lim_{N\to\infty}\lambda\left\{x\in\R:\:\frac{|S_N(x;c,\alpha)|}{\sqrt N}>R\right\}
=\displaystyle\frac{6}{\pi^2}\frac{1}{R^6}\!\left(1+O\!\left(R^{-\frac{12}{31}}\right)\right)\label{intro-tail-estimate-irr}
\end{align}
as $R\to\infty$,  see \cite{Cellarosi-Marklof}.
\end{itemize}
More generally, the convergence of $t\mapsto \frac{1}{\sqrt{N}} S_{\lfloor tN\rfloor}$ as $N\to\infty$ to a random process has been shown by Cellarosi \cite{Cellarosi-curlicue} and Cellarosi and Marklof \cite{Cellarosi-Marklof}.
In this paper we are able to improve the tail estimates \eqref{intro-tail-estimate-rat}-\eqref{intro-tail-estimate-irr} of \cite{Marklof-1999} and \cite{Cellarosi-Marklof}  to the following `optimal' ones (see Remark \ref{rk-optimality} below).
\begin{theorem}\label{main-thm-1}
For all  $0<\varepsilon\leq1$
\begin{align}
\lim_{N\to\infty}\lambda\left\{x\in\R:\:\frac{|S_N(x;0,0)|}{\sqrt N}>R\right\}\label{tail-statement-case-0,0}
=\displaystyle\frac{4\log2}{\pi^2R^4}
\!\left(1+O_\varepsilon\!\left(R^{-2+\varepsilon}\right)\right)
\end{align}
as $R\to\infty$.
\end{theorem}
\begin{theorem}\label{main-thm-2}
Let $(c,\alpha)\notin\Q^2$. Then 
for all  $0<\varepsilon\leq 1$
\begin{align}
\lim_{N\to\infty}\lambda\left\{x\in\R:\:\frac{|S_N(x;c,\alpha)|}{\sqrt N}>R\right\}\label{tail-statement-generic case}
=\displaystyle\frac{6}{\pi^2R^6}\!\left(1+O_\varepsilon\!\left(R^{-2+\varepsilon}\right)\right)
\end{align}
as $R\to\infty$.
\end{theorem}

The recent literature includes several examples of heavy-tailed  limiting distribution arising from the study of problems in number theory 
(e.g. the distribution of visible lattice points by Boca, Cobeli and Zaharescu \cite{Boca-Cobeli-Zaharescu-visible-lattice-points}, the gaps distribution between Farey fractions satisfying divisibility constraints by Boca, Heersink and Spiegelhalter
\cite{Boca-Heersink-Spiegelhalter}, the distribution of short Gauss sums by Demirci Akarsu
\cite{DemirciAkarsu}, the celebrated result by Elkies and McMullen on the gaps in the sequence $(\sqrt{n}\mod1)_{n\geq1}$ 
\cite{Elkies-McMullen}, and the distribution of Frobenius numbers by Marklof
\cite{Marklof-Frobenius-numbers}), 
mathematical physics
(e.g. the distribution of free-path lengths in various Lorentz gas models in the small scatterer limit, see the works of Boca and Zaharescu \cite{MR2274553}, Marklof and Str\"{o}mbergsson 
\cite{Marklof-Strombergsson-Annals-10}-\cite{Marklof-Strombergsson-Lorentz-gas-asymptotic-estimates}-\cite{Marklof-Strombergsson-Lorentz-gas-power-law}, and N\'{a}ndori, Sz\'{a}sz and Varj\'{u}
\cite{Nandori-Szasz-Varju-tail-asymptotics}), 
and ergodic theory 
(e.g. the gap distribution between the slopes of saddle connections in flat surfaces as in the works of Athreya, Chaika and Leli\`{e}vre \cite{MR3330337}, of Uyanik and Work \cite{MR3567252}, as well as the deviations of ergodic averages for toral rotations by Dolgopyat and Fayad \cite{MR3177379}-\cite{MR4179835}). 

\begin{figure}[htbp]
\begin{center}
\hspace{0cm}\includegraphics[width=16.4cm]{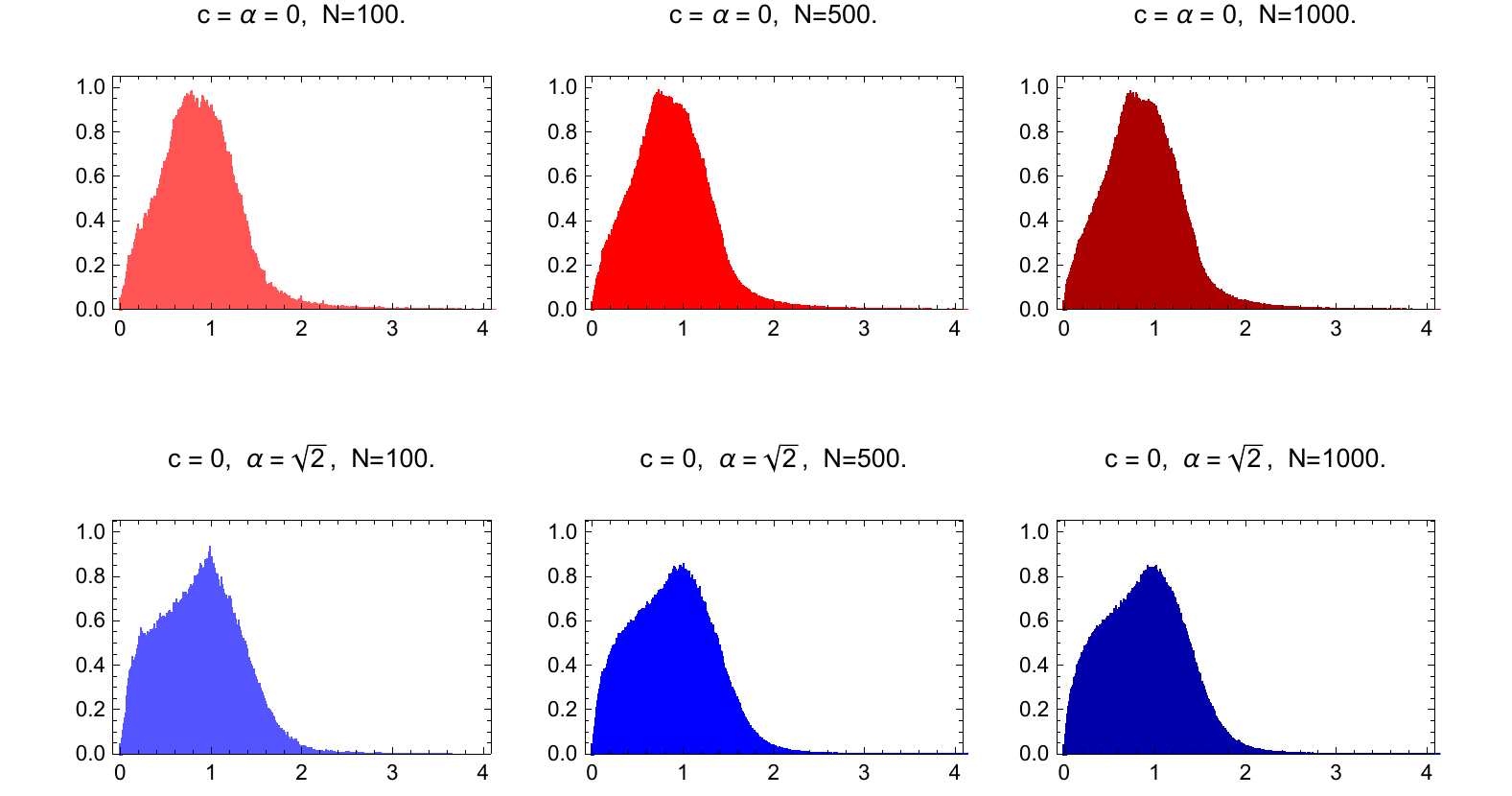}
\caption{\small{Histograms of the empirical distribution of $\frac{1}{\sqrt{N}}|S_N(x;c,\alpha)|$ for $10^6$ uniformly sampled values of $x\in[0,1]$. We use $c=\alpha=0$ (top row),  $(c,
\alpha)=(0,\sqrt{2})$ (bottom row), $N=100$ (left column), $N=500$ (middle column), and $N=1000$ (right column).} \label{fig-six-histograms}}
\end{center}
\end{figure}

\begin{figure}[htbp]
\hspace{-.7cm}\includegraphics[width=17cm]{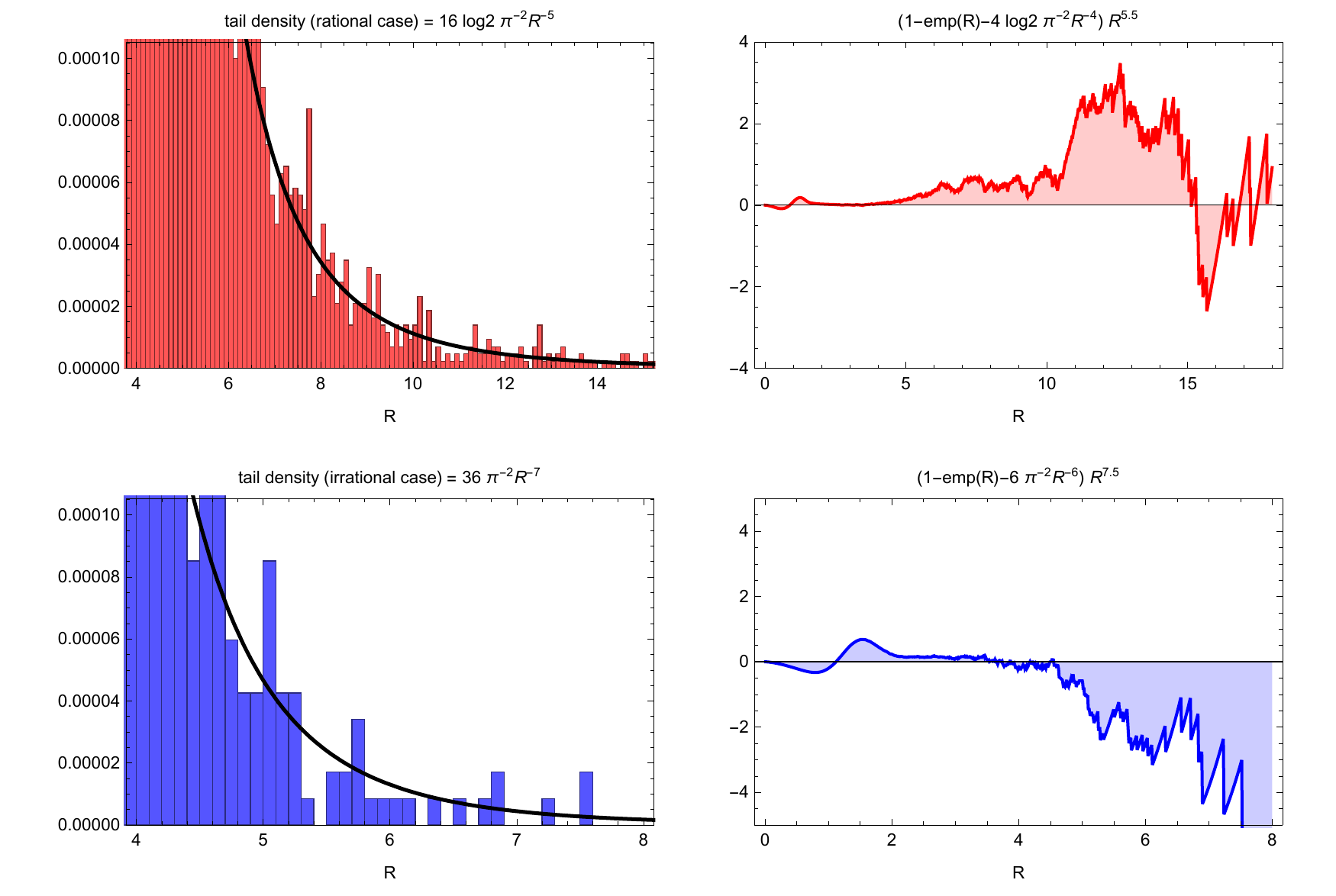}
\caption{\small{We sample  $5\times 10^6$  points according to (an approximation of) the limiting distribution in the rational case (top row) and in the irrational case (bottom row), see Remark \ref{remark-sample-from-limiting-distribution}. The empirical distributions corresponding to our samples are shown in the left panels along with the predicted tail density, obtained by differentiating $1-\frac{4\log2}{\pi^2}R^{-4}$ in the rational case and $1-\frac{6}{\pi^2}R^{-6}$ in the irrational case, with respect to $R$. In the right panels we illustrate Theorems \ref{main-thm-1}-\ref{main-thm-2} for $\varepsilon=\ha$ and plot $(1-\frac{4\log2}{\pi^2}R^{-4}-\mathrm{cdf}(R))R^{5.5}$ in the rational case, and $(1-\frac{6}{\pi^2}R^{-6}-\mathrm{cdf}(R))R^{7.5}$ in the irrational case, where $\mathrm{cdf}(R)$ denotes the empirical distribution function. Our results show the existence of a uniform horizontal band within which these graphs are confined. Given our sample size, we expect to see $\frac{4\log2}{\pi^2}16^{-4} \times5\times 10^6\approx21$ values larger than 16 in the rational case and $\frac{6}{\pi^2}7^{-6}\times 5\times 10^6\approx 26$ values larger than 7 in the irrational case. This is the reason why the ranges of $R$ shown in the two cases end at 16 and 7, respectively. In this way, the accuracy of these graphs is comparable. The bin size in the two histograms on the left is the same.}\label{fig-tails-and-fluctuations}} 
\end{figure}

Theorems \ref{main-thm-1} and \ref{main-thm-2} address the size of fluctuations about the tail of the limiting distribution of rescaled quadratic Weyl sums. Our results can also be interpreted as refined large deviations estimates. In Figure \ref{fig-six-histograms} we illustrate how the distribution of $\frac{|S_N|}{\sqrt{N}}$ approaches  a different limiting distribution as $N\to\infty$ in the rational and irrational cases. The predicted tail asymptotics for the limiting distributions are shown in Figure \ref{fig-tails-and-fluctuations}, along with the fluctuations predicted by Theorems \ref{main-thm-1} and \ref{main-thm-2}.

\begin{remark}
The slightly more general sums 
\begin{align}
\displaystyle\sum_{n=1}^N\e{\left(\tfrac{1}{2}n^2+c_1n+c_0\right)\!x+\alpha n}
\end{align} are considered in \cite{Cellarosi-Marklof}. In this work we are only concerned in the absolute value of such sums, which does not depend on $c_0$, and hence we simplify the notation by setting $c_0=0$ and $c_1=c$.
\end{remark}
\begin{remark}\label{rk-general-rational-case}
The remaining cases not covered by Theorems \ref{main-thm-1} and \ref{main-thm-2} correspond to 
$(c,\alpha)\in\Q^2\smallsetminus\{(0,0)\}$. For \emph{all} rational pairs $(c,\alpha)\in\Q^2$ we have an asymptotic tail formula similar to \eqref{tail-statement-case-0,0}, namely of the form 
\begin{align}
\frac{\mathcal{T}(c,\alpha)}{R^4}\left(1+O_\varepsilon(R^{-2+\varepsilon})\right),\label{tail-statement-other-rationals}
\end{align}
where $\mathcal{T}(c,\alpha)\geq0$ for all $(c,\alpha)\in\Q^2$ and $\mathcal{T}(0,0)=\frac{4\log2}{\pi^2}$.
To find a formula for $\mathcal{T}(c,\alpha)$ in the general rational case, one can use a combination of number-theoretical and group-theoretical ideas which are independent of the techniques used here to obtain  the error term $O_\varepsilon(R^{-2+\varepsilon})$ in  \eqref{tail-statement-case-0,0} and \eqref{tail-statement-other-rationals}. 
We therefore defer the determination of $\mathcal{T}(c,\alpha)$ to a separate work \cite{Cellarosi-Osman-rational-tails}  and here we  only consider the rational case  $c=\alpha=0$ as in Theorem  \ref{main-thm-1}.  It is nevertheless worthwhile to mention that there exist infinitely many rational pairs $(c,\alpha)$ for which $\mathcal{T}(c,\alpha)=0$. For such parameters, such as $(c,\alpha)=(\frac{1}{2(2k-1)},\frac{1}{2(2k-1)})$ with $k\in\N$, the liming distribution of $\frac{S_N}{\sqrt{N}}$ is therefore compactly supported. Other examples of compactly supported limiting distributions can be found in the work of Demirci Akarsu and Marklof on incomplete Gauss sums \cite{DAM2013}, and of Kowalski and Sawin on Kloosterman and Birch sums \cite{Kowalksi-Sawin}.
\end{remark}

In this paper we actually  obtain precise tail estimates for the limiting distribution of products of sums of the form \eqref{theta-sum-intro}, i.e. $\displaystyle\lim_{N\to\infty}\lambda\{x\in\R:\:|\tfrac{1}{N}S_{\lfloor aN\rfloor}(x;c,\alpha) \overline{S_{\lfloor bN\rfloor}(x;c,\alpha)}|>R^2\}$ where, up to relabelling,  $0<a\leq b$. Furthermore, without loss of generality, we may assume that $1= a \leq b$ since 
\begin{align}
&\lim_{N\to\infty}\lambda\left\{x\in\R:\:\left|\tfrac{1}{N}S_{\lfloor aN\rfloor}(x;c,\alpha) \overline{S_{\lfloor bN\rfloor}(x;c,\alpha)}\right|>R^2\right\}=\label{rescaling-1}\\
&=\lim_{M\to\infty}\lambda\left\{x\in\R:\:\left|\tfrac{1}{M}S_{M}(x;c,\alpha) \overline{S_{\lfloor \frac{b}{a}M\rfloor}(x;c,\alpha)}\right|>(\sqrt{a}R)^2\right\}.\label{rescaling-2}
\end{align}
The following two theorems are the main result of this paper.

\begin{theorem}\label{main-thm-3-rat} Let $b\geq1$ and $0<\varepsilon\leq 1$. There exists a constant $R^{\mathrm{rat}}(b,\varepsilon)>0$
such that for all $R>R^{\mathrm{rat}}(b,\varepsilon)$ we have
\begin{align}
\lim_{N\to\infty}\lambda\left\{x\in\R:\: \left|\tfrac{1}{N}S_N(x;0,0)\overline{S_{\lfloor b N\rfloor}(x;0,
0)}\right|>R^2\right\}=
\frac{2D_{\mathrm{rat}}(b)}{\pi^2 R^4}\left(1+O_\varepsilon\!\left(R^{-2+\varepsilon}\right)\right),
\label{statement-main-thm-3-rat}
\end{align}
where 
\begin{align}
D_{\mathrm{rat}}(b)=\begin{cases}
2\log 2,&\mbox{if $b=1$,}\\
2b\,\mathrm{coth}^{-1}(b)+\frac{1}{2}\log(b^2-1)+\frac{b^2}{2}\log(1-\frac{1}{b^2}),&\mbox{if $b>1$.}
\end{cases}
\label{constant-D_rat(chi,chi_b)-in-intro}
\end{align}
\end{theorem}

\begin{theorem}\label{main-thm-3-irr} 
Let $(c,\alpha)\notin\Q^2$, $b\geq1$, and $0<\varepsilon\leq1$. There exists a constant $R^{\mathrm{irr}}(b,\varepsilon)>0$ such that for all $R>R^{\mathrm{irr}}(b,\varepsilon)$ we have
\begin{align}
\lim_{N\to\infty}\lambda\left\{x\in\R:\: \left|\tfrac{1}{N}S_N(x;c,\alpha)\overline{S_{\lfloor b N\rfloor}(x;c,
\alpha)}\right|>R^2\right\}=
\frac{2D_{\mathrm{irr}}(b)}{\pi^2 R^6}\left(1+O_\varepsilon\!\left(R^{-2+\varepsilon}\right)\right),
\label{statement-main-thm-3-irr}
\end{align}
where $D_{\mathrm{irr}}(b)$ is given in \eqref{D_irr-def} and $D_{\mathrm{irr}}(1)=3$.
\end{theorem}
In this paper, we aim to make all the constants implied by the $O$-notations explicit and simple, often at the cost of making such constants bigger (see also Section \ref{section-improve-constants}). Explicit versions of Theorems \ref{main-thm-3-rat} and \ref{main-thm-3-irr} are given by Theorems \ref{thm-rat-epsilon-explicit} and \ref{thm-irr-epsilon-explicit}, respectively.
By taking $b=1$ we obtain
 Theorems \ref{main-thm-1} and \ref{main-thm-2}  as particular cases of Theorems \ref{main-thm-3-rat} and \ref{main-thm-3-irr}.  Their explicit versions are given in Corollaries \ref{cor-explicit-version-thm-1} and \ref{cor-explicit-version-thm-2}, respectively.
To obtain Theorems \ref{main-thm-3-rat} and \ref{main-thm-3-irr}, 
we first prove a more general tail estimate for products of Weyl sums of the form 
\begin{align}
S_N(x;c,\alpha;f)=\sum_{n\in\Z}f\!\left(\frac{n}{N}\right)\e{ \left(\tfrac{1}{2}n^2+cn\right)\!x+n\alpha},\label{def-S_N-f}
\end{align}
where $f:\R\to\R$ is bounded and of sufficient decay at $\pm\infty$ so that the series \eqref{def-S_N-f} is absolutely convergent. Note that  \eqref{def-S_N-f} reduces to \eqref{theta-sum-intro} when $f=\mathbf{1}_{(0,1]}$.
Assuming a certain degree of decay and regularity for two functions $f_1$ and $f_2$,  we can accurately study the tails of the limiting distribution of the product $\tfrac{1}{N} S_N(x;c,\alpha;f_1)\overline{S_N(x;c,\alpha;f_2)}$ as $N\to\infty$, where $x$ is distributed on $\R$ according to the probability measure $\lambda$. Using the uniform decay of a family of  Fourier-like transforms, we define the regularity class $\mathcal{S}_\eta(\R)$, see Section \ref{section-Jacobi-theta-functions}. We then prove the following
\begin{theorem}\label{main-thm-4}
Let $\eta>1$ and let $f_1,f_2\in\mathcal{S}_\eta(\R)$. Then
\begin{align}
&&\lim_{N\to\infty}\lambda\left\{x\in\R:\: \left|\tfrac{1}{N}S_N(x;c,\alpha;f_1)\overline{S_N(x;c,\alpha;f_2)}\right|>R^2\right\}=\nonumber\\
&&=\begin{cases}
\displaystyle\frac{2D_{\mathrm{rat}}(f_1,f_2)}{\pi^2R^4}\left(1+O\!\left(R^{-2\eta}\right)\right),&\mbox{if $c=\alpha=0$,}\\
\\
\displaystyle\frac{2D_{\mathrm{irr}}(f_1,f_2)}{
\pi^2 R^6}\left(1+O\!\left(R^{-2\eta}\right)\right),&\mbox{if $(c,\alpha)\notin \Q^2$,}
\end{cases}
\label{statement-main-thm-4}
\end{align}
where $D_{\mathrm{rat}}(f_1,f_2)$ and $D_{\mathrm{irr}}(f_1,f_2)$ are given  in \eqref{def-D_rat(f_1,f_2)} and \eqref{def-D_irr(f_1,f_2)}, respectively.  The constants implied by the $O$-notations in \eqref{statement-main-thm-4} depend explicitly on $\eta, f_1, f_2$.
\end{theorem}

Both instances of Theorem \ref{main-thm-4} are a combination of a dynamical statement and an analytical tail estimate. The dynamical ingredient is the equidistribution of horocycle lifts under the `stretching' action of the geodesic flow in a suitably defined 5-dimensional homogeneous space $\GamG$. 
In the rational case, the equidistribution takes place on a codimension-2 submanifold (see Theorem \ref{rat_lim_portmenteau}, extending Theorem 5.2 in \cite{Marklof-1999}), while in the irrational case, already studied in \cite{Cellarosi-Marklof},  horocycle lifts equidistribute in the whole space according to its normalised Haar measure (see Theorem \ref{irr_limit_portmenteau}). Therefore, the dynamics singles out a limiting probability measure on $\GamG$, which is different in the two cases we consider, and we are left to study the push-forward of such probability measures via a product of two Jacobi theta functions $\Theta_{f_1}\overline{\Theta_{f_2}}:\GamG\to\C$. The corresponding probability measures on $\C$ are the ones whose tail behaviour we address in Lemmata \ref{L2.2} and \ref{L2.2_irrational}.


\begin{remark}\label{rk-optimality}
 The assumption $\eta>1$ in Theorem \ref{main-thm-4} is crucial, and hence Theorems \ref{main-thm-3-rat}-\ref{main-thm-3-irr}, in which $\eta=1$, do not follow directly. Consequently, the error terms in Theorems \ref{main-thm-3-rat}-\ref{main-thm-3-irr} are `optimal', in the sense that improving the statements to include $\varepsilon=0$ is not possible. Following the strategy of \cite{Cellarosi-Marklof}, we use a careful smoothing procedure and a dynamically defined partition of unity to, in effect, increase the regularity $\eta$. This then allows us to derive \eqref{statement-main-thm-3-rat} and \eqref{statement-main-thm-3-irr} from \eqref{statement-main-thm-4}. By improving a key estimate of \cite{Cellarosi-Marklof}, in which the case $\eta=2$ is considered, and instead pushing $\eta\to1^+$, we achieve the improved error terms in \eqref{statement-main-thm-3-rat} and \eqref{statement-main-thm-3-irr}. 
 
The above-mentioned improvement follows from our Theorem \ref{theorem-uniform-bound-f(t,x)-new}, which we also 
use to illustrate the divergence of the expected displacement of a quantum harmonic oscillator as the initial condition approaches a piecewise $C^1$ function with finitely many jump discontinuities.
\end{remark}

%

%% file: preliminaries-05.tex
\section{Preliminaries}
In this section we define the Jacobi theta function $\Theta_f$, which plays a central role in our analysis. We will see how the quantity  $\tfrac{1}{N} S_N(x;c,\alpha;f_1)\overline{S_N(x;c,\alpha;f_2)}$ can be seen as the product of two  theta functions $\Theta_{f_1}\overline{\Theta_{f_2}}$ evaluated at certain points in the Lie group $G=\sltr\ltimes\R^2$. We also introduce \emph{two} measures on $G$. 
The first one is the Haar measure $\mu$, which is classical and is used in the analysis of the irrational case $(c,\alpha)\notin\Q^2$. The second one, $\mu^{\bn}$,  is the product of the Haar measure on $\sltr$ with an atomic measure on $\R^2$ and is used when addressing the rational case $(c,\alpha)=(0,0)$. The analysis of all other rational cases (see Remark \ref{rk-general-rational-case}) requires similar measures $\mu^{(c,\alpha)}$. We also introduce a lattice $\Gamma<G$ under which $\Theta_{f_1}\overline{\Theta_{f_2}}$ is invariant.

\subsection{The universal Jacobi group 
and $\asltr$} 
\label{section-universal-Jacobi-group-and-ASL}
Let $\h:=\{w\in\C:\:\Im(w)>0\}$ denote the upper half plane. For each $M=\sma{a}{b}{c}{d}\in\sltr$ and $z\in\h$, we let $Mz=\frac{az+b}{cz+d}$ and $\epsilon_M(z)=(cz+ d)/|cz+ d|$.
Let $\tsltr$ denote the universal cover of $\sltr$, i.e.
\begin{align}
\tsltr:=\{[M,\beta_M]:\:M\in\sltr,\:\beta_M \mbox{ a continuous function on $\h$ s.t. }e^{i\beta_M(z)}=\epsilon_M(z)\}.
\end{align}
The group law on $\tsltr$ is  given by
\begin{align}
[M_1,\beta_{M_1}'][M_{2},\beta_{M_2}'']&=[M_1M_2,\beta_{M_1M_2}'''],\qquad\beta_{M_1M_2}'''(z)=\beta_{M_1}'(M_2 z)+\beta_{M_2}''(z),\label{mult-tsltr1}\\
[M,\beta_M]^{-1}&=[M^{-1},\beta_{M^{-1}}'],\qquad \beta_{M^{-1}}'(z)=-\beta_{M}(M^{-1}z).
\end{align}
We identify 
$\tsltr$ with $\h\times\R$ via $[M,\beta_M]\mapsto(z,\phi)=(M i,\beta_M(i))$ and we have an action of $\tsltr$ on $\h\times\R$ given by
\begin{align}[M,\beta_M](z,\phi)=(M z,\phi+\beta_M(z))\label{act-tsltr-on-hxR}.\end{align}
%
Let $\omega:\R^2\times\R^2\to\R$ be the standard symplectic form $\omega(\vecxi,\bm{\xi'})=xy'-yx'$, where $\vecxi=\sve{x}{y},\bm{\xi'}=\sve{x'}{y'}$. The Heisenberg group $\Hei$ is defined as $\R^2\times\R$ with  group law
\begin{align}
(\vecxi,t)(\bm{\xi'},t')=\left(\vecxi+\bm{\xi'},t+t'+\tfrac{1}{2}\omega(\vecxi,\bm{\xi'})\right).\label{mult-Heisenberg}
\end{align}
We consider the \emph{universal Jacobi group} 
$\widetilde{G}=\tsltr\ltimes\Hei$,
with group law
\begin{align}
([M,\beta_M];\vecxi,\zeta)([M',\beta'_{M'}];\bm{\xi'},\zeta')=\left([MM',\beta''_{MM'}];\vecxi+M\bm{\xi'},\zeta+\zeta'+\tha\omega(\vecxi,g\bm{\xi'})\right)\label{mult--univ-Jacobi},
\end{align}
where $\beta''_{MM'}(z)=\beta_M(M'z)+\beta'_{M'}(z)$ as in \eqref{mult-tsltr1}.
Identifying $\widetilde{G}$ with $\h\times\R\times\R^2\times\R$, the Haar measure $\tilde\mu$ on $\widetilde{G}$ is given in coordinates $\tg=(x+i y,\phi;\scriptsize{\sve{\xi_1}{\xi_2}},\zeta)$  by 
$\de \tilde{\mu}(\tilde g)=\frac{1}{y^2}d x\,\de y\,\de\phi\,\de\xi_1\,\de\xi_2\,\de\zeta$. 
Later, we will focus our attention on the special affine group 
$G:=\asltr=\sltr\ltimes\R^2$,
which is the range of the homomorphism 
\begin{align}
\pi:\widetilde{G}\to G,\hspace{.4cm}\pi([M,\beta_M];\vecxi,\zeta)= (M,\vecxi)\in G\label{projection-from-Gtilde-to-G}
\end{align} and hence 
\begin{align}
G\cong\widetilde{G}/\{([I,2\pi k];\bn,\zeta):\: k\in\Z,\: \zeta\in\R\}.\label{from-G-tilde-to-G}
\end{align}
The group law on $G$ is
\begin{align}
(M;\vecxi)(M';\bm{\xi'})=\left(MM';\vecxi+M\bm{\xi'}\right)\label{mult-ASL2R}
\end{align}
and by identifying $G$ with $\h\times[0,2\pi)\times\R^2$ the Haar measure $\mu$ on $G$ is given in coordinates $g=(x+i y,\phi;\scriptsize{\sve{\xi_1}{\xi_2}})$, up to a multiplicative factor, by 
\begin{align}
\de\mu(g)=\frac{3}{\pi^2}\frac{d x\,\de y\,\de\phi\,\de\xi_1\,\de\xi_2}{y^2}.\label{Haar-measure-mu-on-G}
\end{align}
We will also need the measure $\mu^{\bn}$ on $G$, given by 
\begin{align}
\de\mu^{\bn}(g)=\frac{1}{\pi^2} \frac{d x\,\de y\,\de\phi}{y^2}\,\de\nu,\hspace{.3cm}\nu=
\sum_{\mathbf{k}\in\Z^2}\left(\delta_{\sve{0}{0}+\mathbf{k}}+\delta_{\sve{1/2}{0}+\mathbf{k}}+\delta_{\sve{0}{1/2}+\mathbf{k}}\right),\label{measure-mu0-on-G}
\end{align}
which projects to the Haar measure on $\sltr=\h\times[0,2\pi)$ and to an atomic measure 
on $\R^2$. The reason for the normalising factors $\frac{3}{\pi^2}$ and $\frac{1}{\pi^2}$ in \eqref{Haar-measure-mu-on-G} and \eqref{measure-mu0-on-G}  will be clear in Section \ref{section-invariance-of-Theta_f1Theta_2}.

\subsection{Integral transforms on $\LtR$}
For $f\in\LtR$ and $\phi\in\R$, let us define the  operator $\mathscr{I}_\phi:\LtR\to\LtR$ as
\begin{align}
[\mathscr{I}_\phi f](w)=\begin{cases} 
f(w),&\mbox{if $\phi\equiv0$ mod $2\pi$,}\\
f(-w),&\mbox{if $\phi\equiv\pi$ mod $2\pi$,}\\ \\\displaystyle
|\sin\phi|^{-1/2}\!\int_\R \e{\frac{\ha(w^2+w'^2)\cos\phi-w w'}{\sin\phi}}f(w')\de w',&\mbox{if $\phi\equiv\hspace{-.38cm}/\hspace{.2cm}0$ mod $\pi$.}
\end{cases}\:\:\label{def-R-tilde-k-f}
\end{align}
We will see in Section \ref{section-stationary-phase-and-estimate-of-kappa_eta} how to interpret $w\mapsto[\mathscr{I}_\phi f](w)$ as the solution at time $\phi$ to a quantum harmonic oscillator with initial condition $f$. Let us also introduce the function $\phi\mapsto\sigma_\phi$, given by
\begin{align}
\sigma_\phi:=\begin{cases}2\nu,&\mbox{if $\phi=\nu\pi$, $\nu\in\Z$;}\\
2\nu+1,&\mbox{if $\nu\pi<\phi<(\nu+1)\pi$, $\nu\in\Z$.}\end{cases}\label{def-sigma-phi}
\end{align}
Finally, define the operator $\mathscr F_\phi:\LtR\to\LtR$ as 
\begin{align}
[\mathscr F_\phi f] (w)=\e{-\sigma_\phi/8}[\mathscr{I}_\phi f](w)\label{def-script-F-phi}
\end{align}
The reason for the factor $\e{-\sigma_\phi/8}$ in \eqref{def-script-F-phi} is that %
$\lim_{\phi\to 0^\pm}[\mathscr{I}_\phi f](w)=\e{\pm\tfrac{1}{8}}f(w)$ 
 for $f\in\SR$. 
The 1-parameter family  $(\mathscr F_\phi)_{\phi\in\R}$ is a continuous group of unitary operators on $\LtR$. The particular case $\mathscr{F}_\frac{\pi}{2}$ coincides with the classical unitary Fourier transform.   Note that $\phi\mapsto \mathscr F_\phi$ is $4\pi$-periodic and that $\overline{\mathscr{F}_\phi f}=\mathscr{F}_{-\phi}\overline f$. For a more compact notation, unless explicitly stated, we will write $f_\phi$ instead of $\mathscr{F}_\phi f$.

\subsection{The Schr\"{o}dinger-Weil representation of $\widetilde{G}$ on $\LtR$}
By  Iwasawa decomposition, every $M\in\sltr$ can be written uniquely as
\begin{align}M=n_xa_yk_\phi ,\label{Iwasawa-dec}
\end{align}
where 
\begin{align}n_x=\ma{1}{x}{0}{1},\:\:a_y=\ma{y^{1/2}}{0}{0}{y^{-1/2}},\:\:k_\phi=\ma{\cos\phi}{-\sin\phi}{\sin\phi}{\cos\phi},\end{align}
and $z=x+i y\in\h$, $\phi\in[0,2\pi)$.
This allows us to parametrize $\mathrm{SL}(2,\R)$ by $\h\times[0,2\pi)$. We will use the shorthand $(z,\phi)=n_xa_yk_\phi$. Recall the identification of $\tsltr$ with $\h\times\R$. We can extend \eqref{Iwasawa-dec}  to a decomposition of $\tsltr$: for every $\tilde M=[M,\beta_M]\in\tsltr$ we have
\begin{align}
\tilde M=[M,\beta_M]=\tilde n_x\tilde a_y\tilde k_\phi=[n_x,0][a_y,0][k_\phi,\beta_{k_\phi}].\label{Iwasawa-tlstr}
\end{align}

The classical (projective) \emph{Shale-Weil representation} of $\sltr$  can be lifted to a genuine representation of $\tsltr$  using \eqref{Iwasawa-tlstr}, see \cite{Lion-Vergne}: for $f\in\LtR$ let
\begin{align}
[R( \tilde n_x)f](w)&=
\e{\tfrac{1}{2}w^2 x}f(w),\\
[R ( \tilde a_y)f](w)&=
y^{1/4}f(y^{1/2}w),\label{def-R-a-tilde-y}\\
[R(\tilde k_\phi)f](w)&=f_\phi(w).
\end{align}
Recalling \eqref{mult-Heisenberg}, elements of $\Hei$ can be decomposed as 
\begin{align}
\sveH{x}{y}{\zeta}=\sveH{x}{0}{0}\sveH{0}{y}{0}\sveH{0}{0}{\zeta-\tfrac{xy}{2}}
\end{align}
and this allows us to define the \emph{Schr\"{o}dinger representation} of $\Hei$ on $\mathrm L^2(\R)$ as
\begin{align}
[W\!\sveH{x}{0}{0}f](w)&=\e{xw}f(w),\\
[W\!\sveH{0}{y}{0}f](w)&=f(w-y),\label{def-W(0,y)}\\
[W\!\sveH{0}{0}{\zeta}f](w)&= \e{\zeta}f(w).
\end{align}
Combining the two above representations, we can define the \emph{Schr\"{o}dinger-Weil representation} of the universal Jacobi group. 
Recalling the identification of $\widetilde{G}$ with $\h\times\R\times\R^2\times\R$, we  set
\begin{align}
\tilde R(z,\phi;\vecxi,\zeta)=W(\vecxi,\zeta) R(z,\phi),\label{def-tilde-R}
\end{align}
for $(z,\phi;\vecxi,\zeta)\in\widetilde{G}$.

\subsection{Jacobi theta functions}\label{section-Jacobi-theta-functions}
We consider the space of functions $f:\R\to\R$ for which $f_\phi$ defined in \eqref{def-script-F-phi} has a certain decay at infinity, uniformly in $\phi$: denote \begin{align}
\kappa_\eta(f)=\sup_{w,\phi}\left|f_\phi(w)\right|(1+|w|^2)^{\frac{\eta}{2}}.\label{def-kappa_eta(f)}
\end{align} and define 
\begin{align}
\mathcal S_{\eta}(\mathbb{R}):=\left\{f:\R\to\R:\:\kappa_\eta(f)<\infty\right\},\label{def-S_eta}
\end{align}
see \cite{Marklof2007b}. These function spaces play a key role in our analysis when $\eta>1$. 
We stress that the indicator functions $\mathbf{1}_{(0,s]}$ and $\chi_s=\mathbf{1}_{(0,s)}$ belong to $\mathcal S_1$ but not to $\mathcal S_\eta$ for any $\eta>1$. In Section \ref{section-approximation-of-indicators}, we introduce approximations $\chi_s^{(J)}$ of $\chi_s$ which  belong to $\mathcal S_\eta$ with $\eta>1$ and it is crucial to estimate  $\kappa_\eta(\chi_a^{(J)})$. 
Most of the analysis in \cite{Cellarosi-Marklof} deals with the case $\eta=2$; here we see how to obtain the optimal tail bounds \eqref{tail-statement-case-0,0} and \eqref{tail-statement-generic case} by pushing $\eta\to1^+$.

We can now use the representation \eqref{def-tilde-R} to define the most important function in our investigation.
For $(z,\phi;\vecxi,\zeta)\in\widetilde{G}$ and $f\in\mathcal S_\eta(\R)$, with $\eta>1$,  define the \emph{Jacobi theta function} as
\begin{align}
\Theta_f(z,\phi;\vecxi,\zeta)&=\sum_{n\in\Z}\,[\tilde R(z,\phi;\vecxi,\zeta)f](n)=\label{Jacobi-theta-sum-1}\\
&=y^{1/4}e(\zeta-\tha \xi_1\xi_2)\sum_{n\in\Z}\,f_\phi\!\left((n-\xi_2)y^{1/2}\right)\e{\tha(n-\xi_2)^2x+n\xi_1}\label{Jacobi-theta-sum-2},
\end{align}
where $z=x+ i y$, and $\vecxi=\sve{\xi_1}{\xi_2}$. 
The reason for introducing theta functions is that we can rewrite \eqref{def-S_N-f}, after normalization, as 
\begin{align}
N^{-\ha} S_N(x;c,\alpha;f)=\Theta_{f}(x+i N^{-2},0;\sve{\alpha+c x}{0},0).\label{writing-normalized-sums-as-Theta-functions}
\end{align}

\subsection{Invariance of $\Theta_{f_1}\overline{\Theta_{f_2}}$}\label{section-invariance-of-Theta_f1Theta_2}
It is shown in \cite{Marklof2003b} that if $f\in\mathcal{S}_\eta$ with $\eta>1$, then that the function $\Theta_f:\widetilde{G}\to\C$ is invariant under the left action by the lattice $\widetilde\Gamma=\langle\widetilde\gamma_1,\widetilde\gamma_2,\widetilde\gamma_3,\widetilde\gamma_4,\widetilde\gamma_5\rangle$, where
\begin{align}
\widetilde\gamma_1&=\left(\left[\ma{0}{-1}{1}{0},\arg\right];\ve{0}{0},\frac{1}{8}\right)=\left(i,\frac{\pi}{2};\ve{0}{0},\frac{1}{8}\right),\label{invariance-by-gamma1}\\
\widetilde\gamma_2&=\left(\left[\ma{1}{1}{0}{1},0\right];\ve{1/2}{0},0\right)=\left(1+i,0;\ve{1/2}{0},0\right),\\
\widetilde\gamma_3&=\left(\left[\ma{1}{0}{0}{1},0\right];\ve{1}{0},0\right)=\left(i,0;\ve{1}{0},0\right),\\
\widetilde\gamma_4&=\left(\left[\ma{1}{0}{0}{1},0\right];\ve{0}{1},0\right)=\left(i,0;\ve{0}{1},0\right),\label{def-gamma-4-tilde}\\
\widetilde\gamma_5&=\left(\left[\ma{1}{0}{0}{1},0\right];\ve{0}{1},0\right)=\left(i,0;\ve{0}{0},1\right).
\end{align}
Now, let us fix two functions $f_1,f_2\in\mathcal S_\eta$ with $\eta>1$ and consider the product $
\Theta_{f_1}\overline{\Theta_{f_2}}$.
Since $\Theta_{f_1}\overline{\Theta_{f_2}}(z,\phi;\vecxi,\zeta)$ 
does not depend on $\zeta$ and depends only on the value of $\phi\bmod 2\pi$, recalling \eqref{from-G-tilde-to-G}, we consider $\Theta_{f_1}\overline{\Theta_{f_2}}$ as a function on $G$.
This function 
is therefore invariant under the left action by the group $\Gamma=\pi(\widetilde\Gamma)=\langle \gamma_1,\gamma_2,\gamma_3,\gamma_4\rangle$, where $\gamma_i=\pi(\widetilde\gamma_i)$ for $1\leq i\leq 4$.
Therefore 
$\Theta_{f_1}\overline{\Theta_{f_2}}$  is  well-defined 
on the quotient $\GamG$, which is a non-compact 5-dimensional manifold with finite volume.  
A fundamental domain for the action of $\Gamma$ on $G$ is 
\begin{align}
\mathcal F_{\Gamma}=\left\{(z,\phi;\vecxi)\in\mathcal F_{\sltz}\times[0,\pi)\times(-\tha,\tha]^2\right\},\label{fund-dom-Gamma-on-G}
\end{align}
where 
\begin{align}
\mathcal{F}_{\sltz}=\left\{ z\in\h:\: |z|>1,\:-\tha \leq \Re(z)<\tha\right\}\cup\left\{z\in\h:\: |z|=1,\: \Re(z)\leq 0\right\}.\label{fund-dom-sltz-on-H}
\end{align}
Note that $\mathcal{F}_\Gamma$ can be seen as a torus bundle over the unit tangent bundle to the classical modular surface. 
\begin{lemma}\label{lem-mu0-is-Gamma-invariant}
The measure $\mu^{\bn}$ is $\Gamma$-invariant.
\end{lemma}
\begin{proof}
Let $A=A_H\times A_R$ be a measurable rectangle, where $A_H$ is a Borel subset of $\h\times[0,2\pi)$ and $A_R$ is a Borel subset of $\R^2$. Suppose that $A$ has finite $\mu^{\bn}$-measure. Let us show that $\mu^{\bn}(\gamma^{-1} A)=\mu^{\bn}(A)$ for every $\gamma\in\Gamma$. The statement of the lemma follows from this fact. Since $A$ is a measurable rectangle, $\mathbf{1}_A(x+iy,\phi;\vecxi)=\mathbf{1}_{A_H}{x+iy,\phi}\mathbf{1}_{A_{R}}(\vecxi)$. Let $\gamma=(M;\mathbf{v})$. We have
\begin{align}
\mu^{\bn}(\gamma^{-1}A)&=\int_{G}\mathbf{1}_{\gamma^{-1}A}(x+iy,\phi;\vecxi)\de\mu^{\bn}=\int_G\mathbf{1}_A(\gamma(x+iy,\phi;\vecxi))\de\mu^{\bn}\\
&=\int_{G}\mathbf{1}_{A_H}(M(x+iy,\phi))\mathbf{1}_{A_R}(\gamma\cdot\vecxi)\de\mu^{\bn}\\
&=\int_{\R^2}\mathbf{1}_{A_R}(\gamma\cdot\vecxi)\,\de\nu\int_0^\pi\int_0^\infty\int_{-\infty}^\infty \mathbf{1}_{A_H}(M(x+iy,\phi))\frac{\de x\de y\de\phi}{y^2}\label{mu0gammainverseA}
\end{align}
Note that the action of $M$ on $\sltr$ preserves the measure $\frac{\de x\de y\de \phi}{y^2}$. Moreover, using the generators of $\Gamma$, we see that the action of each $\gamma\in\Gamma$ is a bijection on the support of $\nu$, i.e.
\begin{align}
\bigcup_{\veck\in\Z^2}\left\{\sve{0}{0}+\veck,\,\sve{1/2}{0}+\veck,\,\sve{0}{1/2}+\veck\right\}.
\end{align}
Since $\nu$ is atomic with equal weights for its atoms, we obtain $\gamma^{-1}_*\nu=\nu$. Therefore \eqref{mu0gammainverseA} implies $\mu^{\bn}(\gamma^{-1}A)=\mu^{\bn}(A)$.
\end{proof}
The projection of the measure $\mu$ (respectively $\mu^{\bn}$) from $G$ to $\GamG$ will be denoted by $\mu_{\GamG}$ (respectively $\mu^{\bn}_{\GamG}$). These projection measures are well-defined thanks to the $\Gamma$-invariance of the Haar measure $\mu$ and to Lemma \ref{lem-mu0-is-Gamma-invariant}. Our choice of normalisation in \eqref{Haar-measure-mu-on-G} and \eqref{measure-mu0-on-G} yields $\mu_{\GamG}(\mathcal{F}_\Gamma)=\mu^{\bn}_{\GamG}(\mathcal{F}_\Gamma)=1$ and therefore $\mu_{\GamG}$ and $\mu^{\bn}_{\GamG}$ are both probability measures on $\GamG$.

\begin{remark}\label{remark-sample-the-two-measures-mu^0_GamG-and-mu_GamG}
It is not hard to numerically  sample points $(x,y,\phi;\sve{\xi_1}{\xi_2})\in\mathcal{F}_\Gamma$ according to $\mu^{\bn}_{\GamG}$ or $\mu_{\GamG}$, assuming that we have a generator of random numbers uniformly distributed on $[0,1]$, such as $\mathtt{RandomReal[]}$ in the Mathematica$^\copyright$ software. The only nontrivial step is to generate pairs of real numbers $(x,y)\in\mathcal{F}_{\sltz}$ distributed according to the joint density $\frac{3}{\pi y^2}\mathbf{1}_{\mathcal{F}_{\sltz}}(x,y)$.
The marginal density of $x\in[-\ha,\ha)$ is $\int_{\sqrt{1-x^2}}^\infty \frac{3}{\pi y^2}\de y=\frac{3}{\pi\sqrt{1-x^2}}$ and hence its distribution function is $[-\ha,\ha)\ni t\mapsto\ha+\frac{3\arcsin(t))}{\pi}\in[0,1]$. If $x_0$ is uniformly distributed on $[0,1]$, then evaluating the inverse distribution function $[0,1]\ni s\mapsto\sin\!\left(\frac{\pi s}{3}-\frac{\pi}{6}\right)\in[-\ha,\ha)$ at $s=x_0$ yields a random number distributed according to the marginal density of $x$. Since the conditional density of $y$ given $x$ is $\frac{3}{\pi y^2}\mathbf{1}_{\mathcal{F}_{\sltz}}(x,y)/\frac{3}{\pi\sqrt{1-x^2}}=\frac{\sqrt{1-x^2}}{y^2}\mathbf{1}_{\mathcal{F}_{\sltz}}(x,y)$, we find the conditional distribution function of $y$ given $x$ as $[\sqrt{1-x^2},\infty)\ni t\mapsto 1-\frac{\sqrt{1-x^2}}{t}\in[0,1]$. If $y_0$ is uniformly distributed on $[0,1]$, then evaluating the inverse conditional distribution function $[0,1]\ni s\mapsto \frac{\sqrt{1-x^2}}{1-s}\in [\sqrt{1-x^2},\infty)$ at $s=y_0$ yields a random number distributed according to the conditional density of $y$ given $x$. In summary, if $(x_0,y_0)=(\mathtt{RandomReal[]},\mathtt{RandomReal[]})$ is a pair of independent uniformly distributed numbers on $[0,1]$ and we define $x=\sin(\frac{\pi x_0}{3}-\frac{\pi}{6})$ and $y=\frac{\sqrt{1-x^2}}{1-y_0}$, then $(x,y)$ is distributed according to the joint density $\frac{3}{\pi y^2}\mathbf{1}_{\mathcal{F}_{\sltz}}(x,y)$, see Figure \ref{fig-samplexy}.
\begin{figure}[htbp]
\begin{center}
\includegraphics[width=4cm]{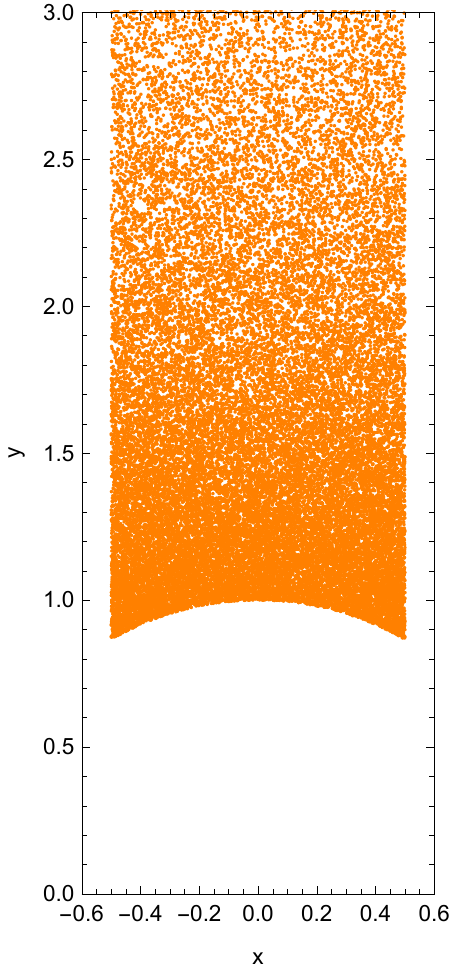}
\caption{\small{A sample of $5\times10^4$ points $(x,y)\in\mathcal{F}_{\sltz}$ according to the probability density  $\frac{3}{\pi y^2}\mathbf{1}_{\mathcal{F}_{\sltz}}(x,y)$ as described in Remark \ref{remark-sample-the-two-measures-mu^0_GamG-and-mu_GamG}. About $32\%$ of the  points  sampled lie above the line $y=3$ and are not shown.}\label{fig-samplexy}}
\end{center}
\end{figure}
Generating the remaining three variables $(\phi,\xi_1,\xi_2)$, independently from $(x,y)$, is easy. In the rational case, if $(\phi_0,\xi_0)=(\mathtt{RandomReal[]},\mathtt{RandomReal[]})$, then we take $\phi=\pi \phi_0$ and we define $(\xi_1,\xi_2)=(0,0)$ if $\xi_0\in[0,\frac{1}{3})$, $(\xi_1,\xi_2)=(\ha,0)$ if $\xi_0\in[\frac{1}{3},\frac{2}{3})$, and $(\xi_1,\xi_2)=(0,\ha)$ if $\xi_0\in[\frac{2}{3},1]$. In the irrational case, we take $(\phi,\xi_1,\xi_2)=(\pi\mathtt{RandomReal[]}, \mathtt{RandomReal[]}, \mathtt{RandomReal[]})$.  
\end{remark}

\subsection{Geodesic and horocycle flows 
}\label{section:geodesic+horocycle+flows}
The group $G$ naturally acts on itself by multiplication. 
Let us consider the action by right-multiplication  by elements of the 1-parameter groups
$\{\Phi^t:\:t\in\R\}$ (the \emph{geodesic flow}) and $\{\Psi^x:\:x\in\R\}$ (the \emph{horocycle flow}), where
\begin{align}
\Phi^t&=\left(a_{e^{-t}};\ve{0}{0}\right)=\left(\left(\ba{cc}e^{-t/2}&0\\0&e^{t/2}\ea\right);\ve{0}{0}\right),\label{def-geodesic-flow}\\
\Psi^{x}&=\left(n_x;\ve{0}{0}\right)=\left(\left(\ba{cc}1&x\\0&1\ea\right);\ve{0}{0}\right).\label{def-horocycle-flow}
\end{align}
We will also consider the right-multiplication by \eqref{def-geodesic-flow} and \eqref{def-horocycle-flow} as flows on  $\GamG$.
In Section \ref{section-approximation-of-indicators} we will also use the flows $\widetilde\Phi^t=\left(\tilde{a}_{e^{-t}};\sve{0}{0},0\right)=\left(\left[\sma{e^{-t/2}}{0}{0}{e^{t/2}},0\right];\sve{0}{0},0\right)$ and $\widetilde\Psi^t=\left(\tilde{n}_x;\sve{0}{0},0\right)=\left(\left[\sma{1}{x}{0}{1},0\right];\sve{0}{0},0\right)$ on $\widetilde G$ and on $\widetilde\Gamma\backslash \widetilde G$.

%% file: limit-theorems-07.tex
\section{Limit Theorems}

In this section we discuss some fundamental dynamical results that will allow us 
to take  $N\to\infty$ in our main theorems. 
Such limit theorems rely on the equidistribution of horocycle lifts in the homogeneous space $\GamG$ under the action of the geodesic flow. 
\subsection{Equidistribution in the smooth rational case}
\begin{theorem}[Limit Theorem, rational case]\label{rational_lim}
Let $F : \Gamma \setminus G \longrightarrow \R$ be a bounded continuous function, and $\lambda$ a Borel probability measure on $\R$ absolutely continuous with respect to the Lebesgue measure. Then
\begin{align*}
\lim_{t \rightarrow \infty} \int_{\R} F(\Gamma(I; \pmb{0})\Psi^u \Phi^t)\: \mathrm{d}\lambda(u) = \int_{\Gamma \setminus G} F \: \mathrm{d}\mu_{\GamG}^{\bn},
\end{align*}
where $\mu_{\GamG}^{\bn}$ is the probability measure defined in  \eqref{measure-mu0-on-G}. 
\end{theorem}

\begin{proof}
Consider the action of $\Gamma$ on $\mathbb{T}^2=\R^2/\Z^2$ via affine transformations, i.e. $(M,\mathbf{v})(\vecxi+\Z^2)=(\mathbf{v}+M\vecxi+\Z^2)$,
and let  $\Sigma = \mathrm{Stab}_{\Gamma}(\bn+\Z^2)$. Using the generators of $\Gamma$ it is not hard to see that $\Sigma = \Gamma_{\theta}\ltimes \Z^2$, where $\Gamma_\theta=\langle\sma{0}{-1}{1}{0},\sma{1}{2}{0}{1}\rangle$ 
is the so-called \emph{theta group}.
Note that the orbit of $\bn+\Z^2$ under $\Gamma$ consists of three points in $\mathbb{T}^2=(-\ha,\ha]^2$, namely $S=\mathrm{Orb}_\Gamma(\bn+\Z^2)=\{\sve{0}{0}, \sve{1/2}{0}, \sve{0}{1/2}\}+\Z^2$. 
By the orbit-stabiliser theorem we have $[\Gamma :\Sigma] = 3$. 
A fundamental domain $\mathcal{F}_\Sigma$ for the action of $\Sigma$ on $G$ is therefore given by three copies of \eqref{fund-dom-Gamma-on-G}. More explicitly,
$\mathcal{F}_\Sigma=\mathcal{F}_\Gamma\cup \rho_1\mathcal{F}_\Gamma+\rho_2\mathcal{F}_\Gamma$,
where 
$\rho_1=\gamma_3^{-1}\gamma_2=\left(\sma{1}{1}{0}{1};\sve{-1/2}{0}\right)$ and 
$\rho_2=\gamma_2\gamma_1=\left(\sma{1}{-1}{1}{0};\sve{1/2}{0}\right)$.
Let $A\subset \sltz\backslash\sltr=\mathcal F_{\sltz}\times[0,\pi)$. We have
\begin{itemize} 
\item[(i)] $(z,\phi;\vecxi)\in A\times\{\sve{1/2}{0}\}$ if and only if $\rho_1(z,\phi;\vecxi)\in \sma{1}{1}{0}{1}A\times\{\sve{0}{0}\}$ and
\item[(ii)] $(z,\phi;\vecxi)\in A\times\{\sve{0}{1/2}\}$ if and only if $\rho_2(z,\phi;\vecxi)\in \sma{1}{-1}{1}{0}A\times\{\sve{0}{0}\}$.
\end{itemize}
If $A$ 
 is a set whose boundary has zero $\mu_{\sltz\backslash\sltr}$-measure, by Sarnak's theorem on equidistribution of long horocycles \cite{Sarnak-equidistribution-horocycles}, we obtain
\begin{align}
&\lim_{t\to\infty}\int_{\R}\mathbf{1}_{A\times\left\{\sve{0}{0}\right\}}(\Gamma(I;\bn)\Psi^x\Phi^t)\de\lambda(x)=\lim_{y\to0}\int_{\R}\mathbf{1}_{A}(\Gamma_\theta n_x a_y)\de\lambda(x)=\\&=\mu_{\Gamma_\theta\backslash\sltr}(A)=\tfrac{1}{3}\mu_{\sltz\backslash\sltr}(A)\label{pf-equidistribution-zero-1}
\end{align}
Similarly, using the fact that M\"{o}bius transformations are isometries on quotients of $\sltr$, 
\begin{align}
&\lim_{t\to\infty}\int_{\R}\mathbf{1}_{A\times\left\{\sve{1/2}{0}\right\}}(\Gamma(I;\bn)\Psi^x\Phi^t)\de\lambda(x)=\lim_{y\to0}\int_{\R}\mathbf{1}_{\sma{1}{1}{0}{1}A}(\Gamma_\theta n_x a_y)\de\lambda(x)=\\&=\mu_{\Gamma_\theta\backslash\sltr}(\sma{1}{1}{0}{1}A)=\tfrac{1}{3}\mu_{\sltz\backslash\sltr}(A),\label{pf-equidistribution-zero-2}\\
&\lim_{t\to\infty}\int_{\R}\mathbf{1}_{A\times\left\{\sve{0}{1/2}\right\}}(\Gamma(I;\bn)\Psi^x\Phi^t)\de\lambda(x)=\lim_{y\to0}\int_{\R}\mathbf{1}_{\sma{1}{-1}{1}{0}A}(\Gamma_\theta n_x a_y)\de\lambda(x)=\\&=\mu_{\Gamma_\theta\backslash\sltr}(\sma{1}{-1}{1}{0}A)=\tfrac{1}{3}\mu_{\sltz\backslash\sltr}(A).\label{pf-equidistribution-zero-3}
\end{align}
The limits \eqref{pf-equidistribution-zero-1}-\eqref{pf-equidistribution-zero-2}-\eqref{pf-equidistribution-zero-3}, along with a standard approximation argument, imply that the limiting distribution, as $t\to\infty$, of $\Gamma(I;\bn)\Psi^x\Phi^t$ (where $x\in\R$ is randomly distributed according to $\lambda$) is given by the product probability  measure $\mu_{\sltz\backslash\sltr}\times\displaystyle\frac{1}{3}\left(\delta_{\sve{0}{0}}+\delta_{\sve{1/2}{0}}+\delta_{\sve{0}{1/2}}\right)=\mu^{\bn}_{\GamG}$. 
\end{proof}

Now we consider $\mathcal{A}\subset\C$ and aim to replace $F$ in the previous theorem by the indicator $\mathbf{1}_{(\Theta_{f_1}\overline{\Theta_{f_2}})^{-1}(\mathcal{A})}$, which is bounded but not continuous. This is done  in the following 
\begin{theorem}\label{rat_lim_portmenteau}
Assume $\lambda$ is any Borel probability measure on $\R$ absolutely continuous with respect to Lebesgue measure. Then, for any $f_1, f_2 \in \mathcal{S}_{\eta}(\R)$ with $\eta > 1$, we have
\begin{align}
\lim_{t\to\infty}\lambda\left\{x \in \R : \Theta_{f_1}\overline{\Theta_{f_2}}(\Gamma(I; \pmb{0})\Psi^u \Phi^t) \in \mathcal{A} \right\}=\mu_{\GamG}^{\bn}((\Theta_{f_1}\overline{\Theta_{f_2}})^{-1}(\mathcal{A}))\label{statement-rat_lim_portmenteau}
\end{align}
whenever $\mathcal{A}\subset\C$ is a measurable set whose boundary is  a null set  
with respect to the push-forward measure 
$(\Theta_{f_1}\overline{\Theta_{f_2}})_*\mu_{\GamG}^{\bn}$. 
\end{theorem}
\begin{proof}
The regularity assumption on $f_1,f_2$ implies that $\Theta_{f_1}\overline{\Theta_{f_2}}:\GamG\to\C$ defined via \eqref{Jacobi-theta-sum-2} is a continuous function. Theorem \eqref{rational_lim} implies that for every bounded continuous $B:\C\to\R$, we have 
\begin{align}
\lim_{t\to\infty}\int_{\R}B(\Theta_{f_1}\overline{\Theta_{f_2}}(\Gamma(I;\bn)\Psi^u\Phi^t)\de\lambda(u)=\int_{\GamG} B\circ (\Theta_{f_1}\overline{\Theta_{f_2}})\, \de\mu_{\GamG}^{\bn}.\label{pf-rat_lim_portmenteau}
\end{align}
The Portmanteau theorem from probability theory (see, e.g., Theorem 11.1.1 in \cite{dudley_2002}) 
allows us to replace $B$ in \eqref{pf-rat_lim_portmenteau} by $\mathbf{1}_{\mathcal A}$, yielding \eqref{statement-rat_lim_portmenteau}.
\end{proof}
\subsection{Equidistribution in the smooth irrational case}
While Theorem \ref{rational_lim} is used in the special rational case $\alpha=c=0$, the following theorem addresses the irrational case $(\alpha,c)\notin\Q^2$.

\begin{theorem}[Limit Theorem, irrational case]\label{irr_limit_portmenteau}
Assume $\lambda$ is any Borel probability measure on $\R$ absolutely continuous with respect to Lebesgue measure. Then for any $(\alpha, c) \notin\Q^2$  and any $f_1, f_2 \in \mathcal{S}_{\eta}(\R)$ with $\eta > 1$, we have
\begin{align}
\lim_{t\to\infty}\lambda\left\{u \in \R : \Theta_{f_1}\overline{\Theta_{f_2}}\left(\Gamma\left(I; \sve{\alpha + c u}{0}\right)\Psi^u \Phi^t\right) \in \mathcal{A} \right\} =\mathrm{\mu}_{\GamG}((\Theta_{f_1}\overline{\Theta_{f_2}})^{-1}(\mathcal{A})).
\end{align}
whenever $\mathcal{A}\subset\C$ is a measurable set whose boundary has measure zero with respect to the push-forward measure $(\Theta_{f_1}\overline{\Theta_{f_2}})_*\mu_{\GamG}$. 
\end{theorem}
\begin{proof}
Because of  the regularity assumption on $f_1$ and $f_2$, this theorem follows immediately from  Corollary 4.3 in \cite{Cellarosi-Marklof}. 
\end{proof}

\subsection{Limit theorems for indicators}
We now extend the limit Theorems \ref{rat_lim_portmenteau} and \ref{irr_limit_portmenteau} to the case when $f_1$ and $f_2$ are indicator functions (which do not belong in $\mathcal{S}_\eta$ for $\eta>1$). Here we write $\chi=\mathbf{1}_{(0,1)}$ and $\chi_b=\mathbf{1}_{(0,b)}$.
%
 Although $\tfrac{1}{N}S_N(x;c,\alpha)\overline{S_{\lfloor b N\rfloor}(x;c,\alpha)}=\Theta_{f_1}\overline{\Theta_{f_2}}(x+i N^{-2},0;\sve{\alpha+c x}{0},0)$ when $f_{1}=\mathbf{1}_{(0,1]}$ and $f_2=\mathbf{1}_{(0,b]}$, here we use indicators of open intervals. The arguments used to prove Theorems \ref{rat_limit_portmenteau-indicators}-\ref{irr_lim_portmenteau-indicators} below will work with either kind of indicators. We will see in Section \ref{section-approximation-of-indicators} how indicators of open intervals are a natural choice. However, replacing one kind of indicator by another will have no effect in the limit, see Corollaries \ref{key-tail-limit-theorem-rat} and \ref{key-tail-limit-theorem-irr}.
%
We also note that in this case the function $\Theta_{\chi}\overline{\Theta_{\chi_b}}$ is not continuous and in fact is only defined $\mu_{\GamG}$-almost everywhere (as well as $\mu_{\GamG}^{\bn}$-almost everywhere as we will see in Section \ref{defining-Theta_chi-almost-everywhere}) on $\GamG$. 
It follows from Corollary 2.4 in \cite{Cellarosi-Marklof} that $\Theta_{f}\in L^4(\GamG,\mu_{\GamG})$ for every $f\in L^2(\R)$ and hence $\Theta_{\chi}\overline{\Theta_{\chi_b}}$ is a well defined element of $L^{2}(\GamG,\mu_{\GamG})$.

\begin{theorem}[Limit theorem, rational case, for indicators]\label{rat_limit_portmenteau-indicators}
Assume $\lambda$ is any Borel probability measure on $\R$ absolutely continuous with respect to Lebesgue measure. Then
\begin{align}
\lim_{t\to\infty}\lambda\left\{u \in \R : \Theta_{\chi}\overline{\Theta_{\chi_b}}(\Gamma(I; \bn)\Psi^u \Phi^t) \in \mathcal{A} \right\} = \mu_{\GamG}^{\bn}((\Theta_{\chi}\overline{\Theta_{\chi_b}})^{-1}(\mathcal{A})).\label{rat_limit_portmenteau-indicators-statement}
\end{align}
whenever $\mathcal{A}\subset\C$ is a measurable set whose boundary has measure zero with respect to the push-forward measure 
$(\Theta_{\chi}\overline{\Theta_{\chi_b}})_*\mu_{\GamG}^{\bn}$. 
\end{theorem}

\begin{theorem}[Limit theorem, irrational case, for indicators]\label{irr_lim_portmenteau-indicators}
Assume $\lambda$ is any Borel probability measure on $\R$ absolutely continuous with respect to Lebesgue measure. Then for any $(\alpha, c) \in \R^2 \setminus \Q^2$
\begin{align}
\lim_{t\to\infty}\lambda\left\{u \in \R : \Theta_{\chi}\overline{\Theta_{\chi_b}}\left(\Gamma\left(I; \sve{\alpha + c u}{0}\right)\Psi^u \Phi^t\right) \in \mathcal{A} \right\} =\mu_{\GamG}((\Theta_{\chi}\overline{\Theta_{\chi_b}})^{-1}(\mathcal{A})).\label{irr_lim_portmenteau-indicators-statement}
\end{align}
whenever $\mathcal{A}\subset\C$ is a measurable set whose boundary has measure zero with respect to the push-forward measure 
$(\Theta_{\chi}\overline{\Theta_{\chi_b}})_*\mu_{\GamG}$. 
\end{theorem}
The approximation arguments needed to prove Theorems \ref{rat_limit_portmenteau-indicators}--\ref{irr_lim_portmenteau-indicators}   are given in Lemmata \ref{lemma-approximation-limit-theorem-1}--\ref{lemma-approximation-limit-theorem-4} below, which are similar to Lemmata 4.6--4.9 in \cite{Cellarosi-Marklof}.

\begin{lemma}\label{lemma-approximation-limit-theorem-1}
 Let $f_1, f_2, h$ be compactly supported, real-valued, Riemann-integrable functions, and assume $h\geq0$. Then for all $\alpha, c\in\R$, we have
\begin{align}
\limsup_{t\to\infty}\int_{\R}  \left|\Theta_{f_1}\overline{\Theta_{f_2}}\left((I,\sve{\alpha+c u}{0})\Psi^u\Phi^t\right)\right| h(u)\de u \leq 2 \left\|f_1\right\|_{2}\left\|f_2\right\|_{2}\left\|h\right\|_{1}.
\end{align}
\end{lemma}
\begin{proof}
Using the Cauchy-Schwarz inequality and Lemma 4.6 in \cite{Cellarosi-Marklof}, we obtain
\begin{align}
&\left|\limsup_{t\to\infty}\int_{\R}\left|\Theta_{f_1}\overline{\Theta_{f_2}}\left((I,\sve{\alpha+c u}{0})\Psi^u\Phi^t\right)\right| h(u)\de u\right|^2\leq \\
&\leq \limsup_{t\to\infty} \int_{\R}\left| \Theta_{f_1}\!\left((I;\sve{\alpha+c u}{0})\Psi^u\Phi^t\right)\right|^2 h(u)\de u \int_{\R}\left| \Theta_{f_2}\!\left((I;\sve{\alpha+c u}{0})\Psi^u\Phi^t\right)\right|^2 h(u)\de u\leq\\
&\leq 4\left\|f_1\right\|_{2}^2 \left\|f_2\right\|_{2}^2 \left\|h\right\|_{1}^2.
\end{align}
\end{proof}

\begin{lemma}\label{lemma-approximation-limit-theorem-2}
Let $f_1,f_2$ be compactly supported, real-valued, Riemann-integrable functions on $\R$, and let $\lambda$ be a Borel probability measure on $\R$ absolutely continuous with respect to Lebesgue measure. Then for all $\alpha,c\in \R$ and every $K>0$ we have
\begin{align}
\limsup_{t\to\infty}\lambda\!\left(\left\{u\in\R:\: |\Theta_{f_1}\overline{\Theta_{f_2}}\left((I;\sve{\alpha+c u}{0})\Psi^u\Phi^t \right)|>K\right\}\right)\leq\frac{2\left\|f_1\right\|_{2}\left\|f_2\right\|_{2}}{K}.\label{lemma-approximation-limit-theorem-2-statement}
\end{align}
\end{lemma}
\begin{proof}
As in the proof of Lemma 4.7 in \cite{Cellarosi-Marklof}, let $\lambda'\in L^1(\R)$ be the density of $\lambda$ and let $m_h$ be the positive measure whose density is $h\in C_c(\R)$. Markov inequality and Lemma \ref{lemma-approximation-limit-theorem-1} yield
\begin{align}
&\lambda\!\left(\left\{u\in\R:\: |\Theta_{f_1}\overline{\Theta_{f_2}}\left((I;\sve{\alpha+c u}{0})\Psi^u\Phi^t \right)|>K\right\}\right)\leq\\
&\leq m_h\!\left(\left\{u\in\R:\: |\Theta_{f_1}\overline{\Theta_{f_2}}\left((I;\sve{\alpha+c u}{0})\Psi^u\Phi^t \right)|>K\right\}\right)+\left\|\lambda'-h\right\|_{1}\leq\\
&\leq \frac{2 \left\|f_1\right\|_{2}\left\|f_2\right\|_{2}\left\|h\right\|_{1}}{K}+\left\|\lambda'-h\right\|_{1}.\label{pf-lemma-approximation-limit-theorem-2}
\end{align}
Finally, using the density of $C_c(\R)$ in $L^1(\R)$ and the continuity of the $L^1$-norm, \eqref{pf-lemma-approximation-limit-theorem-2} implies \eqref{lemma-approximation-limit-theorem-2-statement}
\end{proof}

\begin{lemma}\label{lemma-approximation-limit-theorem-3}
Let $f_1,f_2$ be compactly supported, real-valued, Riemann-integrable functions on $\R$, and let $\lambda$ be a Borel probability measure on $\R$ absolutely continuous with respect to Lebesgue measure. Then for all $\alpha,c\in\R$ and every $\varepsilon>0$ there exists $K_\varepsilon>0$ such that
\begin{align}
\limsup_{t\to\infty}\lambda\!\left(\left\{u\in\R:\:\left|\Theta_{f_1}\overline{\Theta_{f_2}}\left((I;\sve{\alpha+c u}{0})\Psi^u\Phi^t\right)\right|> K_\varepsilon \right\}\right)\leq\varepsilon \left\|f_1\right\|_{2}\left\|f_2\right\|_{2}.
\end{align} 
\end{lemma}
\begin{proof}
The statement follows from Lemma \ref{lemma-approximation-limit-theorem-2} by taking $K_\varepsilon=2\varepsilon^{-1}$.
\end{proof}

\begin{lemma}\label{lemma-approximation-limit-theorem-4}
Let $f_1,f_2$ be compactly supported, real-valued, Riemann-integrable functions on $\R$, and let $\lambda$ be a Borel probability measure on $\R$ absolutely continuous with respect to Lebesgue measure. Let $\eta>1$. Then for all $\alpha,c\in\R$ and every $\varepsilon,\delta>0$ there exist compactly supported functions $\tilde{f}_1,\tilde{f}_2\in\mathcal{S}_\eta(\R)$ such that
\begin{align}
\limsup_{t\to\infty}\lambda\!\left(\left\{u\in\R:\:\left|\Theta_{f_1}\overline{\Theta_{f_2}}\left((I;\sve{\alpha+c u}{0})\Psi^u\Phi^t\right)-\Theta_{\tilde{f}_1}\overline{\Theta_{\tilde{f}_2}}\left( (I;\sve{\alpha+c u}{0})\Psi^u\Phi^t)\right)\right|>\delta\right\}\right)<\varepsilon.
\end{align}
 \end{lemma}
 \begin{proof}
Fix $\varepsilon,\delta>0$. Since $f_1,f_2$ are compactly supported and  $\mathcal{S}_\eta(\R)$ is dense in $L^2(\R)$, we can find compactly supported $\tilde{f}_1,\tilde{f}_2\in\mathcal{S}_\eta(\R)$ such that 
\begin{align}
\frac{4}{\delta}\left(\left\| f_1 \right\|_2 \left\|f_2-\tilde{f}_2\right\|_{2} +\left\|f_1-\tilde{f}_1\right\|_{2} \left\|\tilde{f}_2\right\|_{2}\right) \leq\varepsilon.\label{pf-lemma-approximation-limit-theorem-4-1}
\end{align}
Note that $f\mapsto \Theta_f$ is linear and hence we can write
\begin{align}
\Theta_{f_1}\overline{\Theta_{f_2}}-\Theta_{\tilde{f}_1}\overline{\Theta_{\tilde{f}_2}}= \Theta_{f_1}\overline{\Theta_{f_2 - \tilde f_2}} + \Theta_{f_1 - \tilde f_1}\overline{\Theta_{\tilde f_2}}.
\label{pf-lemma-approximation-limit-theorem-4-2}
\end{align}
Using \eqref{pf-lemma-approximation-limit-theorem-4-2} we obtain the union bound
\begin{align}
\lambda&\!\left(\left\{u\in\R:\:\left|\Theta_{f_1}\overline{\Theta_{f_2}}\left((I;\sve{\alpha+c u}{0})\Psi^u\Phi^t\right)-\Theta_{\tilde{f}_1}\overline{\Theta_{\tilde{f}_2}}\left( (I;\sve{\alpha+c u}{0})\Psi^u\Phi^t)\right)\right|>\delta\right\}\right)\leq\\
&\leq \lambda\!\left(\left\{u\in\R:\:\left|\Theta_{f_1}\overline{\Theta_{f_2-\tilde{f}_2}}\left((I;\sve{\alpha+c u}{0})\Psi^u\Phi^t\right)\right|>\tfrac{\delta}{2}\right\}\right)+\label{pf-lemma-approximation-limit-theorem-4-3}\\
&\hspace{.4cm}+\lambda\!\left(\left\{u\in\R:\:\left|\Theta_{f_1-\tilde{f}_1}\overline{\Theta_{\tilde{f}_2}}\left((I;\sve{\alpha+c u}{0})\Psi^u\Phi^t\right)\right|>\tfrac{\delta}{2}\right\}\right).\label{pf-lemma-approximation-limit-theorem-4-4}
\end{align}
Now we apply Lemma \ref{lemma-approximation-limit-theorem-2} three times with $K=\tfrac{\delta}{2}$ to  the limsup of the expressions in \eqref{pf-lemma-approximation-limit-theorem-4-3}-\eqref{pf-lemma-approximation-limit-theorem-4-4}, thus obtaining \eqref{pf-lemma-approximation-limit-theorem-4-1} as an upper bound.
\end{proof}

We can now proceed to the proof of Theorems \ref{rat_limit_portmenteau-indicators} and \ref{irr_lim_portmenteau-indicators}. Since the Lemmata  \ref{lemma-approximation-limit-theorem-1}--\ref{lemma-approximation-limit-theorem-4}  above hold regardless of the values of $\alpha$ and $c$, we combine the proofs into one.

\begin{proof}[Proof of Theorems \ref{rat_limit_portmenteau-indicators} and \ref{irr_lim_portmenteau-indicators}]
We follow the strategy of the proof of Theorem 4.4 in \cite{Cellarosi-Marklof}. Note that \eqref{rat_limit_portmenteau-indicators-statement}--\eqref{irr_lim_portmenteau-indicators-statement} are equivalent to the following statement: for every bounded continuous function $F:\C\to\R$, we have
\begin{align}
\lim_{t\to\infty}\int_{\R} F\!\left(\Theta_{\chi}\overline{\Theta_{\chi_b}}\left(\Gamma(I; \sve{\alpha+c u}{0})\Psi^u \Phi^t)\right)\right)\de\lambda(u)=
\begin{cases}
\displaystyle\int_{\GamG} F\!\circ\! (\Theta_{\chi}\overline{\Theta_{\chi_b}})\, \de\mu_{\GamG}^{\bn}&\mbox{if $\alpha=c=0$},\\\\
\displaystyle\int_{\GamG} F\!\circ\! (\Theta_{\chi}\overline{\Theta_{\chi_b}})\, \de\mu_{\GamG}&\mbox{if $(\alpha,c)\notin\Q^2$}.
\end{cases}
\label{both_limit_portmenteau-indicators-proof-restatement}
\end{align}
Lemma \ref{lemma-approximation-limit-theorem-3} and the Helly-Prokhorov theorem imply that for every sequence $(t_j)_{j\geq1}$ such that $t_j\to\infty$ as $j\to\infty$ there exists subsequence $(t_{j_\ell})_{\ell\geq1}$ and a probability measure $\rho$ on $\C$ such that for every bounded continuous function $F:\C\to\R$ we have
\begin{align}
\lim_{\ell\to\infty}\int_{\R}F\!\left(\Theta_{\chi}\overline{\Theta_{\chi_b}}\left(\Gamma(I; \sve{\alpha+c u}{0})\Psi^u \Phi^{t_{j_\ell}})\right)\right)\de\lambda(u)=\int_{\C}F\de\rho.
\end{align}
We want to make sure that the limit along the original sequence
\begin{align}
I(\chi,\chi_b)=\lim_{j\to\infty}\int_{\R}F\!\left(\Theta_{\chi}\overline{\Theta_{\chi_b}}\left(\Gamma(I; \sve{\alpha+c u}{0})\Psi^u \Phi^{t_{j}})\right)\right)\de\lambda(u)\label{both_limit_portmenteau-indicators-proof-1}
\end{align}
exists and does not depend on the sequence $(t_j)_{j\geq1}$. To this end, suppose that $F\in C_c^\infty(\C)$. By Lemma \ref{lemma-approximation-limit-theorem-4}, for every $\varepsilon,\delta>0$, there exist compactly supported functions $f_1,f_2\in\mathcal{S}_\eta$, $\eta>1$, such that
\begin{align}
&\int_{\R}\left|F\!\left(\Theta_{\chi}\overline{\Theta_{\chi_b}}\left(\Gamma(I; \sve{\alpha+c u}{0})\Psi^u \Phi^{t_{j}})\right)\right)-F\!\left(\Theta_{f_1}\overline{\Theta_{f_2}}\left(\Gamma(I; \sve{\alpha+c u}{0})\Psi^u \Phi^{t_{j}})\right)\right) \right|\de\lambda(u)\leq\\
&\leq L\int_{\R}\left|\Theta_{\chi}\overline{\Theta_{\chi_b}}\left(\Gamma(I; \sve{\alpha+c u}{0})\Psi^u \Phi^{t_{j}})\right)-\Theta_{f_1}\overline{\Theta_{f_2}}\left(\Gamma(I; \sve{\alpha+c u}{0})\Psi^u \Phi^{t_{j}})\right) \right|\de\lambda(u)\leq\\
&\leq L(\varepsilon+\delta),\label{both_limit_portmenteau-indicators-proof-2}
\end{align}
where $L$ is the Lipschitz constant of $F$. Observe that the limit
\begin{align}
I(f_1,f_2)=\lim_{j\to\infty}\int_{\R}F\!\left(\Theta_{f_1}\overline{\Theta_{f_2}}\left(\Gamma(I; \sve{\alpha+c u}{0})\Psi^u \Phi^{t_{j}})\right)\right)\de\lambda(u)
\end{align}
exists by either Theorem \ref{rat_lim_portmenteau} if $\alpha=c=0$ or Theorem \ref{irr_limit_portmenteau} if $(\alpha,c)\notin\Q^2$. This fact, together with the bound \eqref{both_limit_portmenteau-indicators-proof-2}, imply that the sequence 
\begin{align}
\left(\int_{\R}F\!\left(\Theta_{\chi}\overline{\Theta_{\chi_b}}\left(\Gamma(I; \sve{\alpha+c u}{0})\Psi^u \Phi^{t_{j}})\right)\right)\de\lambda(u)\right)_{j\geq1}
\end{align} is Cauchy and hence the limit \eqref{both_limit_portmenteau-indicators-proof-1} exists. Note that, as $f_1\to\chi$ and $f_2\to\chi_b$ in $L^2(\R)$ we have from \eqref{both_limit_portmenteau-indicators-proof-2} that $I(f_1,f_2)\to I(\chi,\chi_b)$. Therefore 
\begin{align}
I(\chi,\chi_b)=\int_{\C}F\de\rho=
\begin{cases}
\displaystyle\int_{\GamG} F\!\circ\! (\Theta_{\chi}\overline{\Theta_{\chi_b}})\, \de\mu_{\GamG}^{\bn}&\mbox{if $\alpha=c=0$},\\\\
\displaystyle\int_{\GamG} F\!\circ\! (\Theta_{\chi}\overline{\Theta_{\chi_b}})\, \de\mu_{\GamG}&\mbox{if $(\alpha,c)\notin\Q^2$}.
\end{cases}
\end{align}
We have shown that, given $\alpha$ and $c$, the limiting measure $\rho$ is the same for every sequence $(t_j)_{j\geq1}$ tending to infinity, and this implies that  \eqref{both_limit_portmenteau-indicators-proof-restatement} holds. 
\end{proof}

\begin{remark}\label{remark-sample-from-limiting-distribution}
In order to numerically sample points in $\C$ whose law is the push forward of $\mu_{\GamG}^{\bn}$ or $\mu_{\GamG}$ via  $\Theta_{\chi}\overline{\Theta_{\chi_b}}$, it is enough to generate sample points  $(x,y,\phi;\sve{\xi_1}{\xi_2})\in\mathcal{F}_\Gamma$ as explained in Remark \ref{remark-sample-the-two-measures-mu^0_GamG-and-mu_GamG} and then evaluate the function $\Theta_{\chi}\overline{\Theta_{\chi_b}}$ at those points. The evaluation is numerically challenging since the convergence of the series \eqref{Jacobi-theta-sum-2} is very delicate when $(f_1,f_2)=(\chi,\chi_b)$.   In order to obtain Figure \ref{fig-tails-and-fluctuations} (in which $b=1$ and hence $|\Theta_\chi|^2\in\R_{\geq0}$) we use Mathematica$^\copyright$ to find a closed form for $\chi_\phi(w)$ (involving special functions $\mathtt{Erf}$ and $\mathtt{Erfi}$ that can be evaluated to arbitrary precision by the software) and then approximate $|\Theta_\chi|$ by truncating the sum over ${n\in\Z}$ in \eqref{Jacobi-theta-sum-2} to ${-100\leq n\leq 100}$.
\end{remark}


\begin{cor}\label{key-tail-limit-theorem-rat}
Let $\alpha = c = 0$, and $b \geq 1$. Let $\lambda$ be a probability measure, absolutely continuous with respect to Lebesgue measure. Then for every $b\geq1$ and  $R>0$ we have that
\begin{align}
\lim_{N\to \infty}\lambda \!\left\{ x \in \R : \tfrac{1}{N} \left| S_{N}\overline{S_{\lfloor bN \rfloor}} (x; 0, 0)\right| > R^2\right\} = \mu^{\bn}_{\GaG}\!\left\{ |\Thetapair{\chi}{\chi_b}| > R^2\right\}.
\end{align}
\end{cor}

\begin{cor}\label{key-tail-limit-theorem-irr}
Let $( c, \alpha) \in \R^2 \setminus \Q^2$, and $b \geq 1$. Let $\lambda$ be a probability measure, absolutely continuous with respect to Lebesgue measure. Then for every $b\geq1$ and  $R>0$ we have that
\begin{align}
\lim_{N\to \infty}\lambda \!\left\{ x \in \R : \tfrac{1}{ N} \left| S_{N}\overline{S_{\lfloor bN \rfloor}} (x; c, \alpha)\right| > R^2\right\} = \mu_{\GaG}\!\left\{ |\Thetapair{\chi}{\chi_b}| > R^2\right\}.
\end{align}
\end{cor}

\begin{proof}[Proof of Corollaries \ref{key-tail-limit-theorem-rat} and \ref{key-tail-limit-theorem-irr}]
Let $s\in\{1,b\}$. Observe that  $S_{\lfloor s N\rfloor}=S_N^{\ind_{(0,s]}}$ and $S_N^{\chi_s}$ are both finite sums and differ in absolute value by at most $1$. Therefore, for any $(c, \alpha) \in \R^2,$
\begin{align}
\tfrac{1}{\sqrt N} \left| S_{\lfloor sN \rfloor} (x; c, \alpha)\right| =\tfrac{1}{\sqrt N}  \left| S_N^{\chi_s}(x; c, \alpha) \right| + O\! \left( \tfrac{1}{\sqrt N}\right),
\end{align}
and so
\begin{align}
&\tfrac{1}{ N} \labs S_{N}\overline{S_{\lfloor sN \rfloor}} (x; c, \alpha)\rabs =\\
&=\tfrac{1}{ N}  \labs S_N^{\chi}\overline{S_N^{\chi_s}}(x; c, \alpha)\rabs+ O\lp\tfrac{1}{N}\rp+  O \lp \tfrac{1}{\sqrt N} \labs \tfrac{S_N^{\chi}(x; c, \alpha)}{\sqrt N}\rabs\rp+O \lp \tfrac{1}{\sqrt N} \labs\tfrac{S_N^{\chi_s}(x; c, \alpha)}{\sqrt N}\rabs\rp\label{three-O-terms-that-go-to-zero-in-law}.
\end{align}

Theorems \ref{rat_limit_portmenteau-indicators} and \ref{irr_lim_portmenteau-indicators} and \eqref{rescaling-1}-\eqref{rescaling-2} imply that $\labs \tfrac{S_{N}(x; c, \alpha)}{\sqrt N}\rabs$ and $\labs\tfrac{S_N^{\chi_s}(x; c, \alpha)}{\sqrt N}\rabs$ converge in law when $\alpha = c =0$ and $(c, \alpha)\in \R^2 \setminus \Q^2$. By the Cr\'{a}mer-Slutsky theorem (see  \cite{Gut-Probability}, Theorem 11.4 therein) we have that the sum of the three $O$-terms in \eqref{three-O-terms-that-go-to-zero-in-law} goes to zero in law, as $N \to \infty$. Therefore, as $|\cdot|$ is continuous, both $\tfrac{1}{ N} \labs S_{N}\overline{S_{\lfloor sN \rfloor}} (x; c, \alpha)\rabs$ and $\frac{1}{ N}  \labs S_N^{\chi}\overline{S_N^{\chi_s}}(x; c, \alpha)\rabs$ converge in law to the same limit as $N\to\infty$. More precisely,
\begin{align}
\lim_{N\to \infty}\lambda \!\left\{ x \in \R : \tfrac{1}{ N} \left| S_{N}\overline{S_{\lfloor bN \rfloor}} (x; c, \alpha)\right| \in \mathcal{A}\right\} = \begin{dcases*} \mu^{\bn}_{\GaG}\!\left\{ |\Thetapair{\chi}{\chi_b}| \in \mathcal A \right\} & if $\alpha = c = 0$, \\\\ \mu_{\GaG}\!\left\{ |\Thetapair{\chi}{\chi_b}| \in \mathcal{A}\right\} & if $(c, \alpha) \not\in \Q^2$\end{dcases*}
\end{align}
for any measurable set $\mathcal A \subset \R$ whose boundary is of measure zero with respect to the appropriate limiting measure, namely $|\Thetapair{\chi}{\chi_b}|_*\mu^{\bn}_{\GaG}$ or $|\Thetapair{\chi}{\chi_b}|_*\mu_{\GaG}$ as the case may be. Taking $\mathcal A = (R,\infty)$ gives the result.
\end{proof}

%% file: growth-in-the-cusp-09.tex
\section{Growth in the cusp and tail bounds for $\Theta_{f_2} \overline{\Theta_{f_2}}$}
In this section we focus on the size of $\Theta_{f_1}\overline{\Theta_{f_2}}$ when $f_1$ and $f_2$ are of sufficient regularity. The following lemma, which extends Lemma 2.1 in \cite{Cellarosi-Marklof}, is crucial and will be used later. Recall the definition of $\kappa_\eta(f)$ \eqref{def-kappa_eta(f)}.

\begin{lemma}\label{L2.1}
Given $\xi_2 \in \R$, write $\xi_2 = k + \theta$, with $k \in \Z$, and $\theta \in (-\ha,\ha]$. Let $\eta>1$, then, there exist constants $d > 0$ and $C_\eta = d^\eta \zeta(\eta)^2$ such that for every $f_1,f_2 \in \mathcal{S}_{\eta}(\R)$ we have
\begin{align}
\left|\Theta_{f_1}\overline{\Theta_{f_2}} (x + iy, \phi, \pmb{\xi}) - y^{\frac{1}{2}}\,(f_1)_\phi(-\theta y^{\frac{1}{2}})\,\overline{(f_2)_\phi(-\theta y^{\frac{1}{2}})}\right| \leq C_{\eta}\kappa_{\eta}(f_1)\kappa_{\eta}(f_2) y^{-\frac{(\eta - 1)}{2}}\label{L2.1-statement}
\end{align}
for every $y \geq \frac{1}{2}$, and every $x, \phi, \pmb{\xi}$. We may take $d=2^6$.
\end{lemma}

\begin{proof}
The regularity assumption on $f_1$ and $f_2$ 
ensures that the series defining $\Theta_{f_1}$ and $\Theta_{f_2}$ are absolutely convergent. 
Therefore we may write 
\begin{align}
&\Theta_{f_1} \overline{\Theta_{f_2}}(x + iy, \phi, \pmb{\xi})=
 y^{\frac{1}{2}} \sum_{n,m\in\Z} (f_1)_\phi((n-\xi_2)y^{\frac{1}{2}})\,\overline{(f_2)_\phi((m-\xi_2)y^{\frac{1}{2}})}\label{proof-L2.1-0}\\
&\hspace{5.5cm} \e{\tha(n-\xi_2)^2 x + n\xi_1-\tha(m-\xi_2)^2 x - m\xi_1}.\nonumber 
\end{align} 
Since the term $y^{\frac{1}{2}}\,(f_1)_\phi(-\theta y^{\frac{1}{2}})\,\overline{(f_2)_\phi(-\theta y^{\frac{1}{2}})}$ in \eqref{L2.1-statement} corresponds to the term $n=m=k$,  it is enough to prove that 
the above sum restricted to $(n,m)\in\Z^2\smallsetminus \{(k,k)\}$ is bounded in absolute value by $C_{\eta}\kappa_{\eta}(f_1)\kappa_\eta(f_2) y^{\frac{-(\eta - 1)}{2}}$.
Indeed, by  \eqref{def-kappa_eta(f)}, 
\begin{align}
&\:y^{\ha}\sum_{(n,m) \in \Z^2 \setminus{\{(k,k)\}}} \left| (f_1)_\phi((n-\xi_2)y^{\frac{1}{2}})\, \overline{(f_2)_\phi((m-\xi_2)y^{\frac{1}{2}})}\right|\leq\\
 &\leq \:y^{\frac{1}{2}} \sum_{(n,m) \in \Z^2 \setminus{\{(k,k)\}}}\frac{\kappa_\eta(f_1)}{\left(1 + (n-\xi_2)^2y\right)^{\frac{\eta}{2}}}\frac{\kappa_\eta (f_2)}{\left(1 + (m-\xi_2)^2y\right)^{\frac{\eta}{2}}} =\\
&=y^{\frac{1}{2}} \kappa_\eta (f_1)\kappa_\eta (f_2) \left( y^{-\eta} \left(\sum_{n\in \Z\smallsetminus\{k\}} \frac{1}{\left(y^{-1} + (n-\xi_2)^2\right)^{\frac{\eta}{2}} }\right)^2+ \right.\label{proof-L2.1-1}\\ 
&\left.\hspace{4cm}+ 2y^{-\frac{\eta}{2}}\frac{1}{\left(1+(k-\xi_2)^2y\right)^{\frac{\eta}{2}}}\sum_{n \in \Z\smallsetminus\{k\}} \frac{1}{\left(y^{-1} + (n-\xi_2)^2\right)^{\frac{\eta}{2}}}\right). \label{proof-L2.1-2}
%
\end{align}
Now note that the sum appearing in both \eqref{proof-L2.1-1} and \eqref{proof-L2.1-2} is bounded above by 
\begin{align}
\sum_{\ell=1}^\infty \frac{1}{(\ell-\theta)^\eta}+\sum_{j=1}^\infty\frac{1}{(j+\theta)^\eta}.\label{proof-L2.1-3}
\end{align} 
By setting $\theta=\pm\ha$ we can find an upper bound for \eqref{proof-L2.1-3} given by $2\zeta(\eta)(2^\eta-1)$.
We have therefore found an upper bound for  the absolute value of the sum \eqref{proof-L2.1-0}  restricted to $(n,m)\in\Z^2\smallsetminus \{(k,k)\}$ of the form 
\begin{align}
y^{\ha}\kappa_\eta(f_1)\kappa_\eta(f_2)\left(A_\eta y^{-\eta}+B_\eta y^{-\frac{\eta}{2}}\right),\label{upperbound-with-A_eta-and-B_eta}
\end{align}
where $A_\eta=(2\zeta(\eta)(2^\eta-1))^2$ and $B_\eta=4\zeta(\eta)(2^\eta-1)$. Since $y^{-\eta}\leq 2^{\frac{\eta}{2}}y^{-\frac{\eta}{2}}$ for $y\geq\ha$ and $B_\eta\leq 2^{\frac{\eta}{2}}A_\eta$, we obtain \eqref{L2.1-statement} with $C_\eta\geq2^{1+\frac{\eta}{2}}A_\eta$ and the result follows.
\end{proof}

\subsection{Tail bound in the smooth case}
We now consider the function $\Theta_{f_1}\overline{\Theta_{f_2}}:\GamG\to\C$ as \emph{two} possible complex-valued random variables by assuming that its argument is distributed either according to $\mu_{\GamG}^{\bn}$ 
or $\mu_{\GamG}$. In both cases, we use Lemma  \ref{L2.1} to accurately study the tails of these random variables.  

\subsubsection{Tail bound in the smooth rational case}
\begin{lemma}\label{L2.2}
Let $\eta>1$. Then, for every $f_1,f_2\in\mathcal{S}_\eta(\R)$ and  $R^2 > 2^\eta C_\eta \kappa_{\eta}(f_1)\kappa_{\eta}(f_2)$ we have 
\begin{align}
\mu_{\GamG}^{\bn}\left(|\Theta_{f_1} \overline{\Theta_{f_2}}| > R^2\right) = \frac{2D_{\mathrm{rat}}(f_1,f_2)}{\pi^2R^4}\left(1 + O_\eta\!\left((\kappa_\eta (f_1) \kappa_{\eta}(f_2))^{\eta} R^{-2\eta}\right)\right),\label{L2.2-statement}
\end{align}
where $C_\eta$ is as in Lemma \ref{L2.1}, 
\begin{align} 
D_{\mathrm{rat}}(f_1, f_2) = \int_0^{\pi}\left|(f_1)_\phi(0)\,(f_2)_\phi(0)\right|^2 \: \mathrm{d}\phi,\label{def-D_rat(f_1,f_2)}
\end{align}
and the implied constant  by $O_\eta$-notation in \eqref{L2.2-statement} can be taken to be $3 (2^\eta C_\eta)^\eta$. 
\end{lemma}
\begin{proof}
Recall the fundamental domain \eqref{fund-dom-Gamma-on-G} for the action of $\Gamma$ on $G$, and  define
\begin{align}
\mathcal{F}_T = \{(x + iy, \phi; \pmb{\xi}) \in \h \times [0, \pi) \times (-\tha,\tha]^2:\: x \in [-1/2, 1/2),\: y\geq T\}
\end{align}
so that we have the containments $\mathcal{F}_{1} \subseteq \mathcal{F}_{\Gamma} \subseteq \mathcal{F}_{\ha}$. Set $\widetilde{\kappa} = C_{\eta}\kappa_{\eta}(f_1)\kappa_{\eta}(f_2)$. 
We obtain, by Lemma \ref{L2.1},
\begin{align}
&\mu_{\GamG}^{\bn}\left\{(x + iy, \phi; \vecxi) \in \mathcal{F}_{\Gamma} : |\Theta_{f_1} \overline{\Theta_{f_2}}(x + iy, \phi; \vecxi)| > R^2\right\} \leq\\
&\leq \mu_{}^{\bn}\left\{(x + iy, \phi; \vecxi) \in \mathcal{F}_{\ha} : |\Theta_{f_1} \overline{\Theta_{f_2}}(x + iy, \phi; \vecxi)| > R^2\right\}\leq\\
&\leq \mu_{}^{\bn}\left\{(x + iy, \phi; \pmb{\xi}) \in \mathcal{F}_{\ha} : \left|y^{\frac{1}{2}}\,(f_1)_\phi(-\theta y^{\frac{1}{2}})\, (f_2)_\phi(-\theta y^{\frac{1}{2}})\right| +  \widetilde{\kappa} y^{-\frac{\eta - 1}{2}} > R^2\right\}.\label{proof-L2.2-0}
\end{align}
Since $f_1,f_2\in\mathcal{S}_\eta(\R)$, we have
\begin{align}
\left| (f_1)_\phi(-\theta y^{\frac{1}{2}})\,(f_2)_\phi(-\theta y^{\frac{1}{2}})\right| \leq \frac{\kappa_\eta (f_1) \kappa_{\eta}(f_2)}{\left(1 + \left|-\theta y^{\frac{1}{2}}\right|^2\right)^{\eta}} \leq \kappa_\eta (f_1) \kappa_{\eta}(f_2).
\end{align}
Therefore, the condition  
\begin{align}
\left|y^{\frac{1}{2}}(f_1)_\phi(-\theta y^{\frac{1}{2}})\,(f_2)_\phi(-\theta y^{\frac{1}{2}})\right| +  \widetilde{\kappa} y^{-\frac{\eta - 1}{2}} > R^2\label{proof-L2.2-1}
\end{align}
in \eqref{proof-L2.2-0} implies 
\begin{align}
\widetilde{\kappa}\left(y^{\frac{1}{2}} + y^{-\frac{\eta - 1}{2}}\right) > R^2.\label{proof-L2.2-2}
\end{align}
%
%
Furthermore, since $y \geq \frac{1}{2}$ for all points  $(x+iy,\phi;\vecxi)\in\mathcal{F}_{\ha}$, we have $y^{-\tfrac{\eta-1}{2} } \leq  2^{\frac{\eta}{2}}y^{\ha}$ and  \eqref{proof-L2.2-2} in turn implies
\begin{align}
y>\frac{R^4}{\widetilde{\kappa}^2 4^{\eta} }.
\end{align}
Therefore \eqref{proof-L2.2-1}  implies
\begin{align}
\left|y^{\frac{1}{2}}(f_1)_\phi(-\theta y^{\frac{1}{2}})\,(f_2)_\phi(-\theta y^{\frac{1}{2}})\right| +  \widetilde{\kappa} \left(\frac{R^4}{\widetilde{\kappa}^2 4^{\eta}}\right)^{\!-\frac{\eta - 1}{2}} > R^2.\label{proof-L2.2-3}
\end{align}
Let us introduce the shorthand
\begin{align}
I_T^{\bn}(\Lambda)=\mu_{}^{\bn}\left\{(x + iy, \phi; \vecxi) \in \mathcal{F}_{T} : \left| y^{\frac{1}{2}}\,(f_1)_\phi(-\theta y^{\frac{1}{2}})\,(f_2)_\phi(-\theta y^{\frac{1}{2}})\right| >\Lambda\right\}.\label{def-I_T^0}
\end{align}
Then, \eqref{proof-L2.2-0} and \eqref{proof-L2.2-3} yield the upper bound
\begin{align}
&\mu_{\GamG}^{\bn}
\left(|\Theta_{f_1} \overline{\Theta_{f_2}}|>R^2\right)
\leq\label{proof-L2.2-3a}
I_{\ha}^{\bn}\!\left(R^2- \widetilde{\kappa} \!\left(\tfrac{2^\eta \widetilde{ \kappa}}{R^2}\right)^{\eta - 1} \right).
\end{align}
Similarly, we obtain the lower bound
\begin{align}
&\mu_{\GamG}^{\bn}
\left(|\Theta_{f_1} \overline{\Theta_{f_2}}|>R^2\right)
\geq \label{proof-L2.2-4a}
I_{1}^{\bn}\left(R^2+ \widetilde{\kappa} \!\left(\tfrac{2^\eta \widetilde{\kappa}}{R^2}\right)^{\eta - 1} \right).
\end{align}

Using the fact that  $f_1,f_2 \in \mathcal{S}_{\eta}(\R)$, with $\eta > 1$, we see that for fixed $\theta\neq0$ 
we have  $\displaystyle\lim_{y\to\infty }y^\ha (f_1)_\phi(-\theta y^{\frac{1}{2}})\, (f_2)_\phi(-\theta y^{\frac{1}{2}})=0$. 
Thus, in order for the inequality in \eqref{def-I_T^0}
to hold for large $\Lambda$, we must have that $|\theta y^{\frac{1}{2}}|=O(1)$ as $y\to\infty$, i.e. $\theta=O(y^{-\ha})$. In particular, if $\theta$ is fixed, then $\theta=0$ (i.e. $\xi_2\in\Z$).
Recall \eqref{measure-mu0-on-G} and note that integrating the indicator of the event $\{| y^{\frac{1}{2}} (f_1)_\phi(-\theta y^{\ha})\,(f_2)_\phi(-\theta y^{\ha})| >\Lambda\}$ against the measure $\mu_{}^{\bn}$ on $\mathcal F_T$ forces $\vecxi\in\{\sve{0}{0}, \sve{1/2}{0}, \sve{0}{1/2}\}$. More precisely, $\xi_2=\theta=0$ with probability $\frac{2}{3}$ and $\xi_2=\theta=\tha$ with probability $\frac{1}{3}$. The terms with $\xi_2 = \theta = \tha$ do not contribute for $\Lambda$ large. Indeed, note that if $\xi_2=\theta=\tha$, then we have the estimates
\begin{align}
\frac{\left| y^{\frac{1}{2}}(f_1)_\phi(-\tha y^{\ha})\,(f_2)_\phi(-\tha y^{\ha})\right|}{\kappa_\eta(f_1)\kappa_\eta(f_2)}\leq \frac{y^{\ha}}{(1+\frac{1}{4} y)^{{\eta}}}< 4^\eta y^{\ha-\eta}\leq 2^{3\eta}
\end{align}
for $y\geq\ha$. Therefore, if we assume 
\begin{align}
\Lambda\geq 2^{3\eta} \kappa_\eta(f_1)\,\kappa_\eta(f_2),\label{proof-L2.2-6a0}
\end{align}
then 
we obtain
%
\begin{align}
I_T^{\bn}(\Lambda)  &= \frac{3}{\pi^2}\frac{2}{3}\int_{-\ha}^{\ha} \int_{0}^\pi \int_{y > \max\!\left\{T, \left|
(f_1)_\phi(0)\,(f_2)_\phi(0)
\right|^{-2}\Lambda^2\right\}} \frac{1}{y^2}  \de y\,\de\phi\,\de x.\label{proof-L2.2-6}
\end{align} 
Note that \eqref{proof-L2.2-6a0} implies that $\Lambda\geq \kappa_\eta(f_1)\,\kappa_\eta(f_2)$ and hence
$\left|(f_1)_\phi(0)\,(f_2)_\phi(0)\right|^2\leq\frac{\Lambda^2}{T}$ for $T\in\{\ha,1\}$. In this case \eqref{proof-L2.2-6} can be rewritten as
\begin{align}
I_{T}^{\bn}(\Lambda)&=\frac{2}{\pi^2}\int_{-\ha}^\ha \int_0^\pi\int_{|
(f_1)_\phi(0)\,(f_2)_\phi(0)|^{-2}\Lambda^2}^\infty\frac{1}{y^2}\de y\,\de\phi\,\de x=\\
&=\frac{2}{\pi^2}\frac{1}{\Lambda^2}\int_0^\pi |
(f_1)_\phi(0)\,(f_2)_\phi(0)|^2\de\phi.\label{proof-L2.2-7}
\end{align}
We now use the inequality (valid for $0 \leq u \leq \tfrac12$)
\begin{align}
\left| 1 - \frac{1}{(1 \mp u)^2} \right| \leq  6 u,
\end{align}
and hence
\begin{align}
\frac{1}{(1 \mp u)^2} = 1 + O(u)
\end{align}
where the implied constant can be taken to be $6$. Hence, with 
\begin{align}
\Lambda &= R^2\mp \widetilde{\kappa} \!\left(\frac{2^\eta \widetilde{\kappa}}{R^2}\right)^{\eta - 1},
\end{align}
we obtain
\begin{align}
&\frac{1}{\Lambda^2}=\frac{1}{\left(R^2\mp \widetilde{\kappa} \!\left(\frac{2^\eta \widetilde{\kappa}}{R^2}\right)^{\eta - 1} \right)^2}=\frac{1}{R^4} \frac{1}{\left(1 \mp  \!\left(\frac{2^{\eta-1} \widetilde{\kappa}}{R^2}\right)^{\eta} \right)^2} = \frac{1}{R^4} \left( 1+ O(2^{\eta(\eta-1)} \tilde \kappa^\eta R^{-2 \eta}) \right)\label{proof-L2.2-8}
\end{align}
provided
$2^{\eta(\eta-1)} \widetilde{\kappa}^\eta R^{-2\eta}< \tfrac12$. This is implied by the inequality 
\begin{align}
R^2>2^\eta \tilde \kappa.\label{proof-L2.2-9}
\end{align} Therefore, the  assumption \eqref{proof-L2.2-9} and the bound \eqref{proof-L2.2-8} allow us to write
\begin{align}
\frac{1}{\Lambda^2}=\frac{1}{R^4}\left( 1+O_\eta\!\left( \left(\kappa_\eta(f_1)\kappa_\eta(f_2)\right)^\eta R^{-2\eta}\right) \right) ,\label{proof-L2.2-10}
\end{align}
where the constant implied by the $O_\eta$-notation can be taken to be $3 (2^\eta C_\eta)^\eta$. Finally we notice that \eqref{proof-L2.2-9} also guarantees that  \eqref{proof-L2.2-6a0} holds. In fact, we obtain
\begin{align}
\Lambda&=R^2\mp  \widetilde{\kappa} \!\left(\frac{2^\eta \widetilde{\kappa}}{R^2}\right)^{\eta - 1}> \widetilde{\kappa} 
\left(2^{\eta}-1 \right) \geq \tilde \kappa.
\end{align}
Therefore, using the fact that $C_\eta \geq 2^{6\eta}$, we see that \eqref{proof-L2.2-9} implies \eqref{proof-L2.2-6a0} and hence in view of \eqref{proof-L2.2-7} and \eqref{proof-L2.2-10} we see
\begin{align}
I_T^\bn(\Lambda)=\frac{2}{\pi^2}D_{\mathrm{rat}}(f_1,f_2)\, R^{-4}\left(1 + O_\eta\!\left((\kappa_\eta (f_1) \kappa_{\eta}(f_2))^{\eta} R^{-2\eta}\right)\right)\label{proof-L2.2-12}
\end{align}
for $(T,\Lambda)$ either equal to $(\ha,R^2-  \widetilde{\kappa} \!\left(\frac{2^\eta \widetilde{\kappa}}{R^2}\right)^{\eta - 1})$ or to $(1,R^2+\, \widetilde{\kappa} \!\left(\frac{2^\eta \widetilde{\kappa}}{R^2}\right)^{\eta - 1})$. Finally, the bounds \eqref{proof-L2.2-3a} and \eqref{proof-L2.2-4a} imply \eqref{L2.2-statement}.
\end{proof}
Combining Theorem \ref{rat_lim_portmenteau} and Lemma \ref{L2.2}, we get the rational case of Theorem \ref{main-thm-4}.

\subsubsection{Tail bound in the smooth irrational case}

\begin{lemma}\label{L2.2_irrational}
Let $\eta>1$. Then for every  $f_1,f_2\in\mathcal{S}_\eta(\R)$ and $R^2 > 2^\eta C_\eta\kappa_{\eta}(f_1)\kappa_{\eta}(f_2)$ we have 
\begin{align}
\mu_{\GamG}\left(|\Theta_{f_1} \overline{\Theta_{f_2}}| > R^2\right) = \frac{2D_{\mathrm{irr}}(f_1,f_2)}{\pi^2R^6}\left(1 + O_\eta\!\left((\kappa_\eta (f_1) \kappa_{\eta}(f_2))^{\eta} R^{-2\eta}\right)\right),\label{L2.2-irr-statement}
\end{align}
where
\begin{align} 
D_{\mathrm{irr}}(f_1, f_2) = \int_0^{\pi}\int_{-\infty}^{\infty}\left|(f_1)_\phi(w)\,(f_2)_\phi(w)\right|^3 \: \de\phi\,\de w\label{def-D_irr(f_1,f_2)}
\end{align}
and the constant implied by the $O_\eta$-notation in \eqref{L2.2-irr-statement} can be taken to be $7 (2^\eta C_\eta)^\eta$.
\end{lemma}

\begin{proof}
Proceeding as in the proof of Lemma \ref{L2.2}, if we set
\begin{align}
I_T(\Lambda)=\mu\left\{(x + iy, \phi; \vecxi) \in \mathcal{F}_{T} : \left| y^{\frac{1}{2}}(f_1)_\phi(-\theta y^{\frac{1}{2}})\,(f_2)_\phi(-\theta y^{\frac{1}{2}})\right| >\Lambda\right\},
\end{align}
then we obtain the bounds
\begin{align}
I_{1}\!\left(R^2+ \widetilde{\kappa} \!\left(\tfrac{2^\eta \widetilde{\kappa}}{R^2}\right)^{\eta - 1} \right)\leq \mu_{\GamG}^{}
\left(|\Theta_{f_1} \overline{\Theta_{f_2}}|>R^2\right)
\leq\label{proof-L2.2-irr-1}
I_{\ha}\!\left(R^2-  \widetilde{\kappa} \!\left(\tfrac{2^\eta \widetilde{\kappa}}{R^2}\right)^{\eta - 1} \right).
\end{align}
If we assume that  
\begin{align}
\Lambda\geq 2 \kappa_\eta(f_1)\,\kappa_\eta(f_2),\label{proof-L2.2-irr-1a}
\end{align}
then 
$|(f_1)_\phi(-\theta y^{\ha})\,(f_2)_\phi(-\theta y^{\ha})|^2\leq\frac{\Lambda^2}{T}$ for $T\in\{\ha,1\}$ and we obtain
\begin{align}
I_T(\Lambda)=\frac{3}{\pi^2}\int_{-\ha}^{\ha}\int_{-\ha}^{\ha}\int_0^\pi
\iint_{\left\{-\ha<\theta\leq\ha,\:\: y>|
(f_1)_\phi(-\theta y^{\ha})\,(f_2)_\phi(-\theta y^{\ha})|^{-2}\Lambda^2\right\}} \frac{1}{y^2}\de y\,\de \theta\, \de\phi\, \de \xi_1\, \de x.\label{proof-L2.2-irr-2}
\end{align}
Now we use the change of variables  $w = -\theta y^{1/2}$ and notice that 
\begin{align}
&-\tha<\theta<0,\hspace{.3cm}y>|
(f_1)_\phi(-\theta y^{\ha})\,(f_2)_\phi(-\theta y^{\ha})|^{-2}\Lambda^2\\
&\iff 0<w<\infty,\hspace{.3cm} -|w| \left|
(f_1)_\phi(w)\,(f_2)_\phi(w)\right|\Lambda^{-1}<\theta<0
\end{align}
(since \eqref{proof-L2.2-irr-1a} implies that $-|w| \left|
(f_1)_\phi(w)\,(f_2)_\phi(w)\right|\Lambda^{-1}>-\ha$) and similarly
\begin{align}
&0<\theta<\tha,\hspace{.3cm}y>|
(f_1)_\phi(-\theta y^{\ha})\,(f_2)_\phi(-\theta y^{\ha})|^{-2}\Lambda^2\\
&\iff -\infty<w<0,\hspace{.3cm} 0<\theta<|w| \left|
(f_1)_\phi(w)\,(f_2)_\phi(w)\right|\Lambda^{-1}.
\end{align}
Therefore, since $\frac{1}{y^2}\de y=\frac{2\theta^2}{|w|^3}\de w$, \eqref{proof-L2.2-irr-2} becomes
\begin{align}
I_T(\Lambda)&=\frac{6}{\pi^2}\int_0^\pi \left(
\int_{0}^{\infty}\frac{1}{|w|^3}\!\!\int^{0}_{-|w||(f_1)_\phi(w)(f_2)_\phi](w)|\Lambda^{-1}}\!\!\!\!\!\!\!\!\!\!\!\!\theta^2\, \de \theta\, \de w
+
\int_{-\infty}^{0}\frac{1}{|w|^3}\!\!\int_{0}^{|w|(f_1)_\phi(w)(f_2)_\phi(w)|\Lambda^{-1}}\!\!\!\!\!\!\!\!\!\!\!\!\theta^2\,\de\theta\,\de w \!\right)\!\de\phi=\\
&=\frac{2}{\pi^2}\frac{1}{\Lambda^3}\int_0^{\pi}\int_{-\infty}^{\infty}\left|(f_1)_\phi(w)\,(f_2)_\phi(w)\right|^3\de w\,\de\phi.\label{proof-L2.2-irr-3}
\end{align}
Now we use the inequality (valid for $0\leq u \leq \tfrac12$)
\begin{equation}
\left| 1 - \frac{1}{(1\mp u)^3} \right| \leq 14 u
\end{equation}
and hence
\begin{equation}
\frac{1}{(1\mp u)^3} = 1 + O(u)
\end{equation}
where the implied constant can be taken to be $14$. Hence, with
\begin{equation}
\Lambda = R^2 \mp \tilde \kappa \left( \frac{2^\eta \tilde \kappa}{R^2} \right)^{\eta-1}
\end{equation}
we obtain
\begin{equation}
\frac{1}{\Lambda^3} = \frac{1}{R^6} (1+O( 2^{\eta(\eta-1)} \tilde\kappa^\eta R^{-2\eta}) ) \label{proof-L2.2-irr-5}
\end{equation}
provided $2^{\eta(\eta-1)} \tilde \kappa^\eta R^{-2\eta}< \tfrac12$. This is implied by the inequality
\begin{equation}
R^2 > 2^\eta \tilde\kappa. \label{proof-L2.2-irr-6}
\end{equation}
Therefore, assumption \eqref{proof-L2.2-irr-6} and the bound \eqref{proof-L2.2-irr-5} allow us to write
\begin{equation}
\frac{1}{\Lambda^3} = \frac{1}{R^6} \left( 1+ O_\eta( ( \kappa_\eta(f_1) \kappa_\eta(f_2) )^\eta R^{-2 \eta} ) \right), \label{proof-L2.2-irr-7}
\end{equation}
where the constant implied by the $O_\eta$-notation can be taken to be $7 (2^\eta C_\eta)^\eta$. As in the rational case, we notice that \eqref{proof-L2.2-irr-6} implies \eqref{proof-L2.2-irr-1a} and hence in view of \eqref{proof-L2.2-irr-3} and \eqref{proof-L2.2-irr-7} we obtain
\begin{align}
I_T(\Lambda)=\frac{2}{\pi^2}D_{\mathrm{irr}}(f_1,f_2)\, R^{-6}\left(1 + O_\eta\!\left((\kappa_\eta (f_1) \kappa_{\eta}(f_2))^{\eta} R^{-2\eta}\right)\right)
\end{align}
for $(T,\Lambda)$ either equal to $(\ha,R^2- \widetilde{\kappa} \!\left(\frac{2^\eta \widetilde{\kappa}}{R^2}\right)^{\eta - 1})$ or to $(1,R^2+  \widetilde{\kappa} \!\left(\frac{2^\eta \widetilde{\kappa}}{R^2}\right)^{\eta - 1})$. 
The bounds \eqref{proof-L2.2-irr-1} allow us to conclude the proof.
\end{proof}

Combining Theorem \ref{irr_limit_portmenteau} and Lemma \ref{L2.2_irrational}, we get the irrational case of Theorem \ref{main-thm-4}.

%% file: partition-of-unity-04.tex
\section{Dynamical approximation of the indicators}\label{section-approximation-of-indicators}
Let $\chi_s
=\mathbf{1}_{(0,s)}$ and $\chi=\chi_1$.
Our goal is  to accurately estimate $\mu_{\GamG}^{\bn}((\Theta_{\chi}\overline{\Theta_{\chi_b}})^{-1}(\mathcal{A}))$ and $\mu_{\GamG}((\Theta_{\chi}\overline{\Theta_{\chi_b}})^{-1}(\mathcal{A}))$  from \eqref{rat_limit_portmenteau-indicators-statement} and \eqref{irr_lim_portmenteau-indicators-statement} when $\mathcal A$ is the complement of a closed ball of large radius $R$. Since $\chi,\chi_b\notin\mathcal{S}_\eta(\R)$ with $\eta>1$, we cannot apply Lemmata \ref{L2.2} and \ref{L2.2_irrational} directly. Instead, we will apply them to  $\mu_{\GamG}^{\bn}((\Theta_{\chi^{(J)}}\overline{\Theta_{\chi_b^{(J)}}})^{-1}(B_R^c))$ and $\mu_{\GamG}^{}((\Theta_{\chi^{(J)}}\overline{\Theta_{\chi_b^{(J)}}})^{-1}(B_R^c))$, where $\chi^{(J)},\chi_b^{(J)}$ are smooth approximations of $\chi,\chi_b$ belonging in $\mathcal{S}_\eta(\R)$ with $\eta>1$. We now illustrate the construction of such approximations as in \cite{Cellarosi-Marklof}.

\subsection{Smooth trapezoidal approximation of indicators}\label{section:smooth-trapezoidal-approximation-of-indicators}
Let $a\leq b$ and $\varepsilon,\delta\geq0$. We introduce a piecewise quadratic function that is identically zero outside $(a-\varepsilon,b+\delta)$ and identically one within $(a,b)$. We set
%

\begin{align}\label{definition-trapezoid}
T_{a, b}^{\varepsilon, \delta} (w)= 
\begin{dcases*}
0 & $w \leq a - \varepsilon$ \\
\frac{2}{\varepsilon^2}(w - (a - \varepsilon))^2 & $a - \varepsilon < w \leq a - \frac{\varepsilon}{2}$\\
1 - \frac{2}{\varepsilon^2}(w - a)^2 & $a - \frac{\varepsilon}{2} < w \leq a$\\
1  & $a < w \leq b$ \\
1 - \frac{2}{\delta^2}(w - b)^2 & $b < w \leq b + \frac{\delta}{2}$\\
\frac{2}{\delta^2}(w - (b + \delta))^2 & $b + \frac{\delta}{2} < w \leq b + \delta$\\
0 & $w > b + \delta$.
\end{dcases*}
\end{align}
We have the following 
\begin{lemma}[\cite{Cellarosi-Marklof}, Lemma 3.1 therein]\label{lem-T_a,b^epsilon,delta-in-S_2}
Let $\varepsilon,\delta>0$. Then $T_{a,b}^{\varepsilon,\delta}\in\mathcal{S}_2(\R)$.
\end{lemma}

As a consequence, 
$T_{a,b}^{\varepsilon,\delta}\in\mathcal{S}_\eta(\R)$ for every $1<\eta\leq 2$. Using piecewise  functions consisting of higher degree polynomials, one could construct analogues of $T_{a,b}^{\varepsilon,\delta}$ that belong to $\mathcal{S}_\eta(\R)$ for larger $\eta$. This, however, is not necessary in our analysis since, as we shall see in Section \ref{section-final-proofs}, values of $\eta$ closer to 1 will yield better power savings in the error terms.  
Note that 
Lemma \eqref{lem-T_a,b^epsilon,delta-in-S_2} does not apply to $\chi_s$. 
Let $\Delta=T_{\frac{1}{3},\frac{1}{3}}^{\frac{1}{6},\frac{1}{3}}$.
%
%
We have the following partition of unity: 
\begin{align}\label{Chi=sums-of-Deltas-0}
\chi(w) = \sum_{j=0}^\infty \Delta(2^j w) + \sum_{j=0}^\infty \Delta(2^j(1-w)).
\end{align}
It is crucial to rewrite \eqref{Chi=sums-of-Deltas-0} dynamically, using the geodesic flow $\wtPhi$ from Section \ref{section:geodesic+horocycle+flows}. Using \eqref{def-R-a-tilde-y}, \eqref{def-W(0,y)}, and \eqref{def-tilde-R}, we get the identities
\begin{align}
[\tilde R(\wtPhi^{(-2 \log 2)j})f] (w) &= 2^{\frac{j}{2}} f(2^{j}w),\label{flow-2log2j-acting-on-f}\\
[\tilde R(\widetilde{\gamma}_4 \wtPhi^{(-2\log 2)j})f](w)&=2^{\frac{j}{2}}f(2^j(w-1)),
\end{align}
where $\widetilde{\gamma}_4$ is as in \eqref{def-gamma-4-tilde}. We then  rewrite  \eqref{Chi=sums-of-Deltas-0} as
\begin{align}
\chi(w) = \sum_{j=0}^\infty 2^{-\frac{j}{2}}[\tilde R(\wtPhi^{(-2 \log 2)j}) \Delta] (w) + \sum_{j=0}^\infty 2^{-\frac{j}{2}}[\tilde R(\widetilde\gamma_4\wtPhi^{(-2 \log 2)j}) \Delta_{-}](w).\label{Chi=sums-of-Deltas-1}
\end{align}
where $\Delta_{-}(w)=\Delta(-w)$. We observe that each sum in \eqref{Chi=sums-of-Deltas-1} is a weighted sum along the backward orbit of the geodesic flow. The flow renormalizes each sum, mapping one term to the next and providing the exponential weights. 
Furthermore, the action of the geodesic flow $\wtPhi$ allows us to obtain indicators $\chi_s$ from $\chi$ by rescaling, namely
\begin{align}
\chi_s = \sqrt{s}\tilde R(\wtPhi^{2\log s})\chi.\label{chi_s-from-chi_1}
\end{align} 
Therefore, using the group law \eqref{mult--univ-Jacobi}, we obtain
\begin{align}
\chi_s(w) &= \sqrt{s}\left(\sum_{j=0}^\infty 2^{-\frac{j}{2}}[\tilde R(\wtPhi^{2 \log (\frac{s}{2^{j}})}) \Delta](w) + \sum_{j=0}^\infty 2^{-\frac{j}{2}}[\tilde R(\widetilde\gamma_4 \wtPhi^{2 \log (\frac{s}{2^{j}})})\Delta_{-}](w)\right)\label{chi_s=sums-of-Deltas}\\
&=\chi_{s,L}(w)+\chi_{s,R}(w),\label{chi_s=sums-of-Deltas00}
\end{align}
where $\chi_{s,L}=T_{0,\frac{s}{3}}^{0,\frac{s}{3}}$ and $\chi_{s,R}=T_{\frac{2s}{3},s}^{\frac{s}{3},0}$. 
A natural way to construct an approximation $\chi_s^{(J)}$ of $\chi_s$ that is dynamically compatible with the action of the geodesic flow is to simply truncate the sums in  \eqref{chi_s=sums-of-Deltas} to $0\leq j\leq J-1$ for some $J\geq1$. This gives a smooth trapezoid function that equals $T_{\frac{s}{3\cdot 2^{J-1}},s-\frac{s}{3\cdot 2^{J-1}}}^{\frac{s}{6\cdot 2^{J-1}},\frac{s}{6\cdot 2^{J-1}}}$ and  can be written as
\begin{align}
\chi_s^{(J)}(w) &= \sqrt{s}\left(\sum_{j=0}^{J-1} 2^{-\frac{j}{2}}[\tilde R(\wtPhi^{2 \log (\frac{s}{2^{j}})}) \Delta](w) + \sum_{j=0}^{J-1} 2^{-\frac{j}{2}}[\tilde R(\widetilde\gamma_4 \wtPhi^{2 \log (\frac{s}{2^{j}})})\Delta_{-}](w)\right)\label{Chi-J=sums-of-Deltas-0}\\
&=\chi_{s,L}^{(J)}(w)+\chi_{s,R}^{(J)}(w),\label{chi_s^J=chi_s_L^J+chi_s_R^J-1}
\end{align}
where 
$\chi_{s,L}^{(J)}=T_{\frac{s}{3\cdot 2^{J-1}},\frac{s}{3}}^{\frac{s}{6\cdot 2^{J-1}},\frac{s}{3}}$ and $\chi_{s,R}^{(J)}=T_{\frac{2s}{3},s-\frac{s}{3\cdot 2^{J-1}}}^{\frac{s}{3},\frac{s}{6\cdot 2^{J-1}}}$, respectively. When $s=1$, we will drop $s$ from the notation and simply write $\chi$, $\chi_{L}$, $\chi_R$, $\chi^{(J)}$, $\chi_{L}^{(J)}$, and $\chi_{R}^{(J)}$.

\subsection{Defining $\Theta_{\chi_s}$ almost everywhere}\label{defining-Theta_chi-almost-everywhere}
By linearity, \eqref{Jacobi-theta-sum-1}, and \eqref{chi_s=sums-of-Deltas} we may write
\begin{align}
\Theta_{\chi_s}(\tg):=\sqrt{s}\left(\sum_{j=0}^{\infty}2^{-\frac{j}{2}}\Theta_{\Delta}(\tg \wtPhi^{2\log \left(\frac{s}{2^j}\right)})+\sum_{j=0}^{\infty}2^{-\frac{j}{2}}\Theta_{\Delta_-}(\tg \widetilde{\gamma}_4\wtPhi^{2\log \left(\frac{s}{2^j}\right)})\right).\label{def-Theta_chi_s}
\end{align}
Since $\chi_s \notin \mathcal{S}_\eta$ for $\eta>1$ this theta function is not defined for all $\tilde g \in \widetilde G$. On the other hand, for $J\geq1$ we have
\begin{equation}
\Theta_{\chi_{s}^{(J)}}(\tg)=\calL_{s}^{(J)}(\tg)+\calR_{s}^{(J)}(\tg\widetilde{\gamma}_4^s)\label{Theta-trapezoid_s=sum-of-L-and-R}
\end{equation}
where
\begin{align}
\calL_{s}^{(J)}(\tg)&= \Theta_{\chi^{(J)}_{s,L}}( \tilde g) = \sqrt{s}\sum_{j=0}^{J-1}2^{-\frac{j}{2}}\Theta_{\Delta}(\tg \wtPhi^{2\log \left(\frac{s}{2^j}\right)}),\\
\calR_{s}^{(J)}(\tg)&=  \Theta_{\chi^{(J)}_{s,L}(-\cdot)} (\tilde g) = \sqrt{s}\sum_{j=0}^{J-1}2^{-\frac{j}{2}}\Theta_{\Delta_{-}}(\tg  \wtPhi^{2 \log \left(\frac{s}{2^j}\right)}).
\end{align}
Note that these theta functions for the smoothed approximations are well defined on all of $\widetilde G$ by Lemma \ref{lem-T_a,b^epsilon,delta-in-S_2}. By splitting the summation in \eqref{def-Theta_chi_s} into blocks of length $J$ as in the proof of Lemma 3.14 in \cite{Cellarosi-Marklof}, we can write
\begin{align}
\Theta_{\chi_s}(\tg)=\sum_{k=0}^\infty 2^{-k\frac{J}{2}}\left( \calL_{s}^{(J)}(\tg\wtPhi^{-(2\log2)kJ})+\calR_{s}^{(J)}(\tg\widetilde{\gamma}_4\wtPhi^{-(2\log2)kJ})\right).\label{Theta_chi_s-as-sum-over-k}
\end{align}
The term $k=0$ in \eqref{Theta_chi_s-as-sum-over-k} corresponds to \eqref{Theta-trapezoid_s=sum-of-L-and-R} and will provide the main contribution to our estimates, while the other terms, which are obtained by the renormalising action of the backward geodesic flow, will contribute to the error term. 

As in Lemma \ref{L2.2} and Lemma \ref{L2.2_irrational}, we need to consider the push forward of the measures $\mu_{\GamG}^{\bn}$ and $\mu_{\GamG}$ via the function $\Theta_{\chi}\overline{\Theta_{\chi_b}}$, and hence it is crucial that latter is well defined on a  set of full $\mu_{\GamG}^{\bn}$ and $\mu_{\GamG}$-measure. When $s=1$, Theorem 3.10 in \cite{Cellarosi-Marklof} shows that there is a $\widetilde\Gamma$-invariant set $\widetilde W\subset \widetilde G$  (explicitly defined by means of  a Diophantine condition) of full $\widetilde\mu$-measure such that the series in \eqref{def-Theta_chi_s} are absolutely convergent for all $\tg\in\widetilde W$. Furthermore, the $\widetilde\mu$-a.e. defined function $\Theta_{\chi}(\tg)$ is $\tilde\Gamma$-invariant and hence we can write it as  $\Theta_{\chi}(\tilde\Gamma\tg)$.  All other values of $s>0$ can be handled in the same way using \eqref{chi_s-from-chi_1}. To describe the set $\widetilde W$ we use a set of coordinates on $\widetilde G$ different from $(x+iy,\phi;\vecxi,\zeta)\in\h\times\R\times\R^2\times\R$ on $\widetilde G$. We consider the flow $\widetilde\Phi^t$ together with its stable, unstable, and neutral manifolds. Letting
$\tilde n_+(v,\rho):=\left(\left[\sma{1}{v}{0}{1},0\right];\sve{\rho}{0},0\right)$ and  $\tilde n_-(u,\sigma):=\left(\left[\sma{1}{0}{u}{1},\arg(u\cdot+1)\right];\sve{0}{\sigma},0\right)$ parametrise the unstable and stable manifolds respectively, we can write any $\tg\in\widetilde G$ as
\begin{align}
\tg=\tilde n_+(v,\rho)\tilde n_-(u,\sigma)\wtPhi^t([I,0];\sve{0}{0},\tau)\label{dynamical-coordinates-tilde}
\end{align}
with $v,\rho,u,\sigma,t,\tau \in\R$. 
The set $\widetilde W$ of full $\widetilde\mu$-measure consists of all elements of the form \eqref{dynamical-coordinates-tilde} such that $v+\frac{1}{u}$ is Diophantine of type $\kappa>1$. Note that the projection $W=\pi(\widetilde W)\subset G$ via \eqref{projection-from-Gtilde-to-G} is of full $\mu$-measure and $\Gamma$-invariant.
If we write \eqref{dynamical-coordinates-tilde} with respect to our usual coordinates on $\widetilde \Gamma\backslash \widetilde G$, we obtain
\begin{align}
\tilde\Gamma n_+(v,\rho)n_-(u,\sigma)\tilde\Phi^t([I,0];\sve{0}{0},\tau)=\tilde\Gamma \left(n_x a_y k_\phi;\sve{\xi_1}{\xi_2},\zeta\right)
\end{align}
where $x=v+\frac{u}{u^2+e^{2t}}$, $y=\frac{e^t}{u^2+e^{2t}}$, $\xi_1=\rho+v\sigma$, $\xi_2=\sigma$, $\zeta=\tau$, and $\phi$ is determined by $\cos\phi=\frac{\sgn(u)e^{t}}{\sqrt{u^2+e^{2t}}}$ and $\sin\phi=\frac{\sgn(u)u}{\sqrt{u^2+e^{2t}}}$. In particular, since $\xi_1$ and $\xi_2$ do not depend on $u$, we have that for \emph{every} $\vecxi\in\R^2$ there exists a set $U_{\vecxi}\subset\tsltr$ of full $\mu_{\tsltr}$-measure such that $\Theta_{\chi_s}(x+i y,\phi;\vecxi,\zeta)$ is well-defined for every $(x+i y,\phi)\in U_{\vecxi}$, every $\zeta\in\R$, and every $s>0$. %
For $\vecxi\in\R^2$, we let $\widetilde W_\vecxi=\displaystyle \bigcup_{\widetilde \gamma\in\widetilde\Gamma}\widetilde\gamma(U_{\vecxi}\times\{\vecxi\}\times\R)$ and consider its projection $W_{\vecxi}:=\pi(\widetilde W_\vecxi)$, which is $\Gamma$-invariant 
since $\pi(\tilde\Gamma)=\Gamma$.
In particular, 
$W_{\bn}$ is $\Gamma$-invariant and 
 is of full $\mu^{\bn}$-measure and hence $\Gamma W_{\bn}$ is of full $\mu^{\bn}_{\GamG}$-measure. Therefore, since $\Theta_{\chi}(\tg)$ and $\Theta_{\chi_b}(\tg)$ are simultaneously defined (and $\tilde\Gamma$-invariant) for all  $\tg\in\tilde W_{\bn}$ and the product $\Theta_{\chi}(\tg)\overline{\Theta_{\chi_b}}(\tg)$ only depends on $g=\pi(\tg)$, then $\Theta_{\chi}\overline{\Theta_{\chi_b}}(\Gamma g)$  is defined for all $\Gamma g$ in the full $\mu_{\GamG}^{\bn}$-measure set  $\Gamma W_{\bn}$. 
\subsection{Expansion of the product $\Theta_{\chi}\overline{\Theta_{\chi_b}}$}\label{subsection-writing-product-of-Thetas}
Le us now use \eqref{Theta_chi_s-as-sum-over-k} to expand $\Theta_{\chi}\overline{\Theta_{\chi_b}}(g)$. 
Depending on whether we are considering the irrational case or the rational case, we consider $g\in W$ or $g\in W_{\bn}$ and we consider $\tg\in \tilde W$ or $\tg\in \tilde W_{\bn}$ such that $\pi(\tg)=g$. Then, we obtain
\begin{align}
\Theta_{\chi}\overline{\Theta_{\chi_b}}(g)&=
\Theta_{\chi^{(J)}}\overline{\Theta_{\chi_b^{(J)}}}(g)+\label{Theta_chi_aTheta_chi_b=Theta_chi_a^JTheta_chi_b^J+}\\
&~~~~+\sum_{k +\ell > 0} 2^{-\frac{(k + \ell)J}{2}}\calL_{1}^{(J)}(\tg \wtPhi^{-(2\log 2)k J})\overline{\calL_{b}^{(J)}(\tg \wtPhi^{-(2\log 2)\ell J})}+\label{Theta_chi_aTheta_chi_b_LL}\\
&~~~~+\sum_{k +\ell > 0}2^{-\frac{(k + \ell)J}{2}}\calL_{1}^{(J)}(\tg \wtPhi^{-(2\log 2)k J})\overline{\calR_{b}^{(J)}(\tg \tilde \gamma_4 \wtPhi^{-(2\log 2)\ell J})}+\label{Theta_chi_aTheta_chi_b_LR}\\
&~~~~+\sum_{k +\ell > 0} 2^{-\frac{(k + \ell)J}{2}}\calR_{1}^{(J)}(\tg \wtPhi^{-(2\log 2)k J})\overline{\calL_{b}^{(J)}(\tg \wtPhi^{-(2\log 2)\ell J})}+\label{Theta_chi_aTheta_chi_b_RL}\\
&~~~~+\sum_{k +\ell > 0}2^{-\frac{(k + \ell)J}{2}}\calR_{1}^{(J)}(\tg \wtPhi^{-(2\log 2)k J})\overline{\calR_{b}^{(J)}(\tg \tilde \gamma_4 \wtPhi^{-(2\log 2)\ell J})}.\label{Theta_chi_aTheta_chi_b_RR}
\end{align}
Note that each term in \eqref{Theta_chi_aTheta_chi_b_LL} only depends on $g=\pi(\tg)$ since we can once again use the renormalisation given by the geodesic flow to write 
\begin{align}
\calL_{1}^{(J)}(\tg \wtPhi^{-(2\log 2)k J})\overline{\calL_{b}^{(J)}(\tg \wtPhi^{-(2\log 2)\ell J})}&= \Theta_{\chi^{(J)}_{1,L}}(\tg \wtPhi^{-(2\log 2)k J})\overline{\Theta_{\chi^{(J)}_{b,L}}(\tg \wtPhi^{-(2\log 2)\ell J})}=\\
&= \Theta_{\chi^{(J,k)}_{1,L}}(\tg)\overline{\Theta_{\chi^{(J,\ell)}_{b,L}}(\tg)}=
\Theta_{\chi^{(J,k)}_{1,L}}\overline{\Theta_{\chi^{(J,\ell)}_{b,L}}} (g),
\end{align}
where $\chi^{(J, m)}_{s,L} = R(\wtPhi^{-(2\log 2)m J})\chi^{(J)}_{s,L}$ for $(s,m)\in\{(1,k),(b,\ell)\}$. A similar argument holds for all other summands in \eqref{Theta_chi_aTheta_chi_b_LR}--\eqref{Theta_chi_aTheta_chi_b_RR}. We now rewrite \eqref{Theta_chi_aTheta_chi_b=Theta_chi_a^JTheta_chi_b^J+}--\eqref{Theta_chi_aTheta_chi_b_RR} as
\begin{align}
\Theta_{\chi}\overline{\Theta_{\chi_b}}(g) = \: &\Theta_{\chi^{(J)}}\overline{\Theta_{\chi_b^{(J)}}}(g) + \sum_{k + \ell > 0}2^{-\frac{(k + \ell)J}{2}}\calLL^{(J)}_{k, \ell}(g) + \sum_{k + \ell > 0} 2^{-\frac{(k + \ell)J}{2}}\calLR^{(J)}_{k,\ell}(g) + \label{preparing-Theta_chi_aTheta_chi_b-for-union-bound-1}\\
 + &\sum_{k + \ell > 0} 2^{-\frac{(k + \ell)J}{2}}\calRL^{(J)}_{k,\ell}(g) + \sum_{k + \ell > 0} 2^{-\frac{(k + \ell)J}{2}}\calRR^{(J)}_{k,\ell}(g),\label{preparing-Theta_chi_aTheta_chi_b-for-union-bound-2}
\end{align}
with the obvious notations, where we have dropped the explicit dependence on 
$b$.

%% file: stationary-phase-lemma-9.tex
\section{Estimating $\kappa_\eta$ for $1<\eta\leq2$}\label{section-stationary-phase-and-estimate-of-kappa_eta}

Recall the function $\chi_{s,L}^{(J)}$  introduced in Section  \ref{section:smooth-trapezoidal-approximation-of-indicators} and used, along with its mirror image  $\chi_{s,L}^{(J)}(-\cdot)$, in Section \ref{defining-Theta_chi-almost-everywhere} to define $\calL_{s}^{(J)}$ and $\calR_{s}^{(J)}$. These functions are also implicitly used in \eqref{preparing-Theta_chi_aTheta_chi_b-for-union-bound-1}-\eqref{preparing-Theta_chi_aTheta_chi_b-for-union-bound-2}.
 It is crucial for us to estimate $\kappa_\eta(\chi_{s}^{(J)})$, $\kappa_\eta(\chi_{s,L}^{(J)}),$ and $\kappa_\eta(\chi_{s,L}^{(J)}(-\cdot))$, extending the case $\eta=2$ is considered in \cite{Cellarosi-Marklof}.
 
 In this section we will write $f(t,x)=[\mathscr{I}_t f](x)$, where $\mathscr{I}_t$ is the integral operator defined in \eqref{def-R-tilde-k-f}. 

\begin{lemma} The function $f(t,x)=[\mathscr{I}_t f](x)$ satisfies the Schr\"{o}dinger equation
\begin{align}
\begin{cases}\displaystyle\frac{i}{2\pi}\partial_t f(t,x)=\frac{1}{2}\left(-\frac{1}{4\pi^2}\frac{\de^2}{\de x^2}+x^2\right)f(t,x).\\
f(0,x)=f(x).
\end{cases}\label{quantum-harmonic-oscillator}
\end{align}
\end{lemma}
\begin{proof}
The fact that $f(0,x)=f(x)$ is part of the definition \eqref{def-R-tilde-k-f}. A direct computation yields
\begin{align}
\frac{i}{2\pi}\partial_t f(t,x)=&|\sin(t)|^{-\ha}\int_\R\e{\tha(x^2+w^2)\cot(t)-x w\csc(t)}\nonumber\\
&~~~~\left(-\tfrac{i}{4\pi}\cot(t)+\tha(x^2+w^2)\csc^{2}(t)-x w\cos(t)\csc^2(t)\right)f(w)\,\de w
\end{align}
and
\begin{align}
&\tha\left(-\tfrac{1}{4\pi^2}\tfrac{\de^2}{\de x^2}+x^2\right)f(t,x)=|\sin(t)|^{-\ha}\int_\R \e{\tha(x^2+w^2)\cot(t)-xw\csc(t)}\nonumber\\
&~~~~~~~~~~~~\left(\tha x^2\cot^2(t)+\tha w^2\csc^2(t)-xw\cos(t)\csc^2(t)-\tfrac{i\cot(t)}{4\pi}+\tha x^2\right)f(w)\,\de w.
\end{align}
The result then follows from the identity $1+\cot^2(t)=\csc^2(t)$.
\end{proof}

The following lemma can be seen as the consequence of a `coarse' stationary phase argument, and is very important in our bounds for $\kappa_\eta$.

\begin{lemma}[\cite{Cellarosi-Marklof},  Lemma 3.3 therein]\label{lem-coarse-stationary-phase-new}
Let $f$ be a compactly supported function with bounded variation. 
 Let $f(t,x)$ be the solution to \eqref{quantum-harmonic-oscillator} with initial condition $f$.  Then 
\begin{align}
|f(t,x)| \leq \|f\|_{\mathrm{sp}}\label{statement-lem-coarse-stationary-phase-new}
\end{align}
uniformly in $t$, where $\|f\|_{\mathrm{sp}}:= \max\!\left\{2\|f\|_1,3(\|f\|_\infty+V(f))\right\}$  and $V(f)$ is the total variation of $f$. In particular, if $f$ is bounded and piecewise $C^1$, then $V(f)=\|f'\|_1$. 
\end{lemma}

\begin{remark}\label{sp-norms-for-chi-functions}
We can  bound   the norm in  \eqref{statement-lem-coarse-stationary-phase-new} for the functions $\chi_s$, $\chi_s^{(J)}$, $\chi_{s,L}$, $\chi_{s,L}^{(J)}$ introduced in Section \ref{defining-Theta_chi-almost-everywhere} as follows. Since 
$\|T_{a,b}^{\varepsilon,\delta}\|_1=b-a+\tfrac{\varepsilon+\delta}{2}$, 
$\|T_{a,b}^{\varepsilon,\delta}\|_\infty=1$, and
$V(T_{a,b}^{\varepsilon,\delta})=2$, we get
\begin{align}
\|\chi_s\|_{\mathrm{sp}}&=\max\{2s, 9 \},\label{sp-norm-of-chi-s}\\
\|\chi_s^{(J)}\|_{\mathrm{sp}}&=\max\{2(1-2^{-J})s, 9 \}\leq \max\{2s, 9 \},\\
\|\chi_{s,L}\|_{\mathrm{sp}}&=\max\{s, 9 \},\label{sp-norm-of-chi-s,L}\\
\|\chi_{s,L}^{(J)}\|_{\mathrm{sp}}&=\max\{(1-2^{-J})s, 9\}\leq \max\{s, 9 \}.
\end{align}
\end{remark}


The following theorem will be used to estimate $\kappa_\eta(\chi_{s}^{(J)})$ and $\kappa_\eta(\chi_{s,L}^{(J)})$. It is in this theorem that we make the assumption that $1<\eta\leq 2$ for the first time.

\begin{theorem}\label{theorem-uniform-bound-f(t,x)-new}
Let $f$ be a  bounded, compactly supported $C^1$ function such that $\|f\|_1$, $\|f'\|_1$, and $\|f''\|_1$ are all finite.  Let $f(t,x)$ be the solution to \eqref{quantum-harmonic-oscillator} with initial condition $f$.  
Define $h^{p,q}(w)=w^p f^{(q)}(w)$.
Then for all $1<\eta\leq2$, $b\geq1$ and all $x,t\in\R$ we have 
\begin{align}
|f(t,x)|(1+|x|^2)^{\frac{\eta}{2}}
\leq 2^{\frac{\eta}{2}} \max&\left\{\|h^{0,0}\|_{\mathrm{sp}},\: b^{\eta-1}\left(\|h^{0,1}\|_1+\|h^{1,0}\|_{\mathrm{sp}}\right),\right.\label{statement-theorem-uniform-bound-f(t,x)-new}\\
&\left.~~ b^{\eta-2}\left(\|h^{1,1}\|_1+\|h^{0,2}\|_1+\|h^{0,0}\|_1+\|h^{2,0}\|_{\mathrm{sp}}\right)\right\}.\nonumber
\end{align}
\end{theorem}
\begin{proof}
We first use the bound coming from applying Lemma \ref{lem-coarse-stationary-phase-new} directly to $f$. 
For $|x|<1$ and all $t\in \R$ we get
\begin{align}
|(1+|x|^2)^{\frac{\eta}{2}}f(t,x)|\leq 2^{\frac{\eta}{2}} \|f\|_{\mathrm{sp}}.
\label{proof-theorem-uniform-bound-f(t,x)-0-new}
\end{align}
%
Let $b\geq 1$ and assume that $1\leq |x|<b$. Using integration by parts we get
\begin{align}
f(t,x)&= |\sin(t)|^{-\ha}\int_{\R} \e{\tha(x^2+w^2)\cot(t)} f(w)\,\e{-xw \csc(t)}\,\de w=\\
&= \left[- |\sin(t)|^{-\ha}\frac{\sin(t)}{2\pi i x} \e{\tha(x^2+w^2)\cot(t)}f(w)\right]_{-\infty}^\infty+ \label{proof-theorem-uniform-bound-f(t,x)-1-new}\\
&|\sin(t)|^{-\ha}\frac{\sin(t)}{2\pi i x}
\int_{\R}\e{\tha(x^2+w^2)\cot(t)-xw\csc(t)}
\left[f'(w)+2\pi i \cot(t)w f(w)\right]\de w\label{proof-theorem-uniform-bound-f(t,x)-2-new}.
\end{align}
The boundary term in \eqref{proof-theorem-uniform-bound-f(t,x)-1-new} is zero since $f$ is compactly supported. The integral in \eqref{proof-theorem-uniform-bound-f(t,x)-2-new}
splits into two pieces by linearity.
The first half is bounded in absolute value trivially by $\|f'\|_1$. To bound the absolute value of the second half of the integral (together with the prefactor $|\sin(t)|^{-\ha}$), we apply Lemma \ref{lem-coarse-stationary-phase-new} to the function $h^{1,0}$. 
We obtain
\begin{align}
\left|(1+|x|^2)^\frac{\eta}{2} f(t,x)\right|\leq\frac{(1+|x|^2)^\frac{\eta}{2}}{|x|}
\left(\|h^{0,1}\|_1+ \|h^{1,0}\|_{\mathrm{sp}}\right).
\end{align}
Since $|x|\geq1$, we can write $(1+|x|^2)^{\frac{\eta}{2}}|x|^{-1}\leq 2^{\frac{\eta}{2}}|x|^{\eta-1}$.  Therefore, since $\eta>1$ and  $1\leq |x|<b$,  we have
\begin{align}
\left|(1+|x|^2)^\frac{\eta}{2} f(t,x)\right|\leq 2^{\frac{\eta}{2}} b^{\eta-1}(\|h^{0,1}\|_1+\|h^{1,0}\|_{\mathrm{sp}}),
\label{proof-theorem-uniform-bound-f(t,x)-2a-new}
\end{align}
for all in $t\in\R$. 
In the remaining range $|x|\geq b$ we start from \eqref{proof-theorem-uniform-bound-f(t,x)-1-new}--\eqref{proof-theorem-uniform-bound-f(t,x)-2-new} and perform one more integration by parts. We get
\begin{align}
&f(t,x)=\\
&=|\sin(t)|^{-\ha}\left(\frac{\sin(t)}{2\pi i x}\right)^2\int_{\R}\e{\tha(x^2+w^2)\cot(t)-xw\csc(t)}\times\label{proof-theorem-uniform-bound-f(t,x)-3-new}\\
&~~~~\times\left(4\pi i \cot(t)wf'(w)+f''(w)+2\pi i\cot(t)f(w)-4\pi^2\cot^2(t) w^2f(w)\right)
\,\de w\label{proof-theorem-uniform-bound-f(t,x)-4-new}.
\end{align}
We split  \eqref{proof-theorem-uniform-bound-f(t,x)-3-new}--\eqref{proof-theorem-uniform-bound-f(t,x)-4-new} into four integrals by linearity. The first one yields the bound $\frac{1}{\pi} |\sin(t)|^\ha|\cos(t)|\|h^{1,1}\|_1$, the second the bound $\frac{1}{4\pi^2}|\sin(t)|^{\frac{3}{2}}\|f''\|_1$, and the third the bound $\frac{1}{2\pi}|\sin(t)|^{\ha}|\cos(t)|\|f\|_1$. The fourth integral is dealt with by applying Lemma \ref{lem-coarse-stationary-phase-new} to the function $h^{2,0}$ in order to avoid the possible divergence for $t$ close to an integer multiple of $\pi$. We thus obtain
\begin{align}
\left|(1+|x|^2)^\frac{\eta}{2} f(t,x)\right|\leq\frac{(1+|x|^2)^{\frac{\eta}{2}}}{|x|^2}\left(\|h^{1,1}\|_1+\|h^{0,2}\|_1+\|h^{0,0}\|_1+\|h^{2,0}\|_{\mathrm{sp}}\right).
\end{align}
Since $|x|\geq b\geq1$ and $\eta\leq2$, we have
\begin{align}
\left|(1+|x|^2)^\frac{\eta}{2} f(t,x)\right|\leq 2^{\frac{\eta}{2}} b^{\eta-2}\left(\|h^{1,1}\|_1+\|h^{0,2}\|_1+\|h^{0,0}\|_1+\|h^{2,0}\|_{\mathrm{sp}}\right)\label{proof-theorem-uniform-bound-f(t,x)-5-new}
\end{align}
for all $t\in\R$. 
Combining \eqref{proof-theorem-uniform-bound-f(t,x)-0-new}, \eqref{proof-theorem-uniform-bound-f(t,x)-2a-new} and  \eqref{proof-theorem-uniform-bound-f(t,x)-5-new},  we obtain \eqref{statement-theorem-uniform-bound-f(t,x)-new}.
\end{proof}

We now extend Lemma 3.1 in \cite{Cellarosi-Marklof}, in which the case $\eta=2$ is considered, by proving the following 
\begin{prop}\label{prop-key-estimate-kappa_eta(chi_s^J)}
Fix $s\geq1$, and $J\geq1$ 
as above. Then for $1<\eta\leq2$ we have 
\begin{align}
\kappa_\eta(\chi_{s}^{(J)})\leq K(s)\,  2^{(\eta-1)J}\hspace{1cm}\mbox{and}\hspace{1cm}\kappa_\eta(\chi_{s,L}^{(J)}(\pm\cdot))\leq K_{L}(s)\,   2^{(\eta-1)J},\label{statement-prop-key-estimate-kappa_eta(chi_s^J)}
\end{align}
where $K(s)$ and $K_{L}(s)$ are given explicitly in  \eqref{piecewise-formula-C3-new} and \eqref{piecewise-formula-C3L-new} respectively. 
\end{prop}
\begin{proof}
Let $f=\chi_{s}^{(J)}=T_{\frac{s}{3\cdot 2^{J-1}},s-\frac{s}{3\cdot 2^{J-1}}}^{\frac{s}{6\cdot 2^{J-1}},\frac{s}{6\cdot 2^{J-1}}}$. Let us estimate the norms featured in \eqref{statement-theorem-uniform-bound-f(t,x)-new}. Using the notation $h^{p,q}$ as in Theorem \ref{theorem-uniform-bound-f(t,x)-new} and the definition \eqref{definition-trapezoid} directly (see also Remark \ref{sp-norms-for-chi-functions}), we obtain 
$\|h^{0,0}\|_{\mathrm{sp}}\leq \max\{2s,9\}$,
$\|(h^{0,0})'\|_1=\|h^{0,1}\|_1=2$,  
$\|h^{1,0}\|_1=\frac{1}{2}s^2(1-2^{-J})\leq\frac{s^2}{2}$, 
$\|h^{1,0}\|_{\infty}\leq s$,
$\|h^{1,1}\|_1=s$,
$\|(h^{1,0})'\|_1\leq\|h^{0,0}\|_1+\|h^{1,1}\|_1\leq 2s$,
$\|h^{1,0}\|_{\mathrm{sp}}\leq\max\{s^2,9s\}$,
$\|h^{0,2}\|_1=\frac{24\cdot 2^J}{s}\leq 24\cdot 2^J$,
$\|h^{2,0}\|_1\leq\frac{s^3}{3}$,
$\|h^{2,0}\|_\infty\leq s^2$,
$\|h^{2,1}\|_1\leq s^2$,
$\|(h^{2,0})'\|_1\leq2\|h^{1,0}\|_1+\|h^{2,1}\|_1\leq2 s^2$, and 
$\|h^{2,0}\|_{\mathrm{sp}}\leq \max\{\frac{2}{3}s^3, 9 s^2\}$.
 Using the definition \eqref{def-kappa_eta(f)} and Theorem \ref{theorem-uniform-bound-f(t,x)-new}  with $b=
2^{J}$, 
we obtain
\begin{align}
&\sup_{\phi,w}|f(\phi,w)|(1+|w|^2)^{\frac{\eta}{2}}\leq\\
&\leq2^{\frac{\eta}{2}}\max\left\{\max\{2s,9\}, \textcolor{white}{2^{2^{2^{2}}}}\right.
\\
&~~~~~~~~~~~~~~~~~~~~
2^{(\eta-1)J}\left(2+\max\{s^2,9s\}
\right),\\
&~~~~~~~~~~~~~~~~~~~\left.
2^{(\eta-2)J}\left(s+24\cdot 2^J+s+\max\{\tfrac{2}{3}s^3,9s^2\}\right)\right\}\leq\\
&\leq2^{\frac{\eta}{2}} 2^{(\eta-1)J} \max\left\{2s,9,2+\max\left\{9s,s^2\right\},s+24+\max\left\{\tfrac{9}{2}s^2,\tfrac{1}{3}s^3\right\}\right\}\leq\\
&\leq K(s)\, 
2^{(\eta-1)J},\label{where-C_3-appears-1}
\end{align}
where
\begin{align}
K(s)=48+2s+\max\left\{9s^2,\tfrac{2}{3}s^3\right\}.
\label{piecewise-formula-C3-new}
\end{align}
We can repeat the same argument letting $f=\chi_{s,L}^{(J)}=T_{\frac{s}{3\cdot 2^{J-1}},\frac{s}{3}}^{\frac{s}{6\cdot 2^{J-1}},\frac{s}{3}}$ instead. Also in this case the norms featured in \eqref{statement-theorem-uniform-bound-f(t,x)-new} can be estimated directly using \eqref{definition-trapezoid}, although the computations are rather tedious. 
We have
$\|h^{0,0}\|_{\mathrm{sp}}\leq  \max\{s, 9 \}$,
$\|(h^{0,0})'\|_1=\|h^{0,1}\|_1=2$,  
$\|h^{1,0}\|_1=\frac{55}{432}s^2 (1-2^{-2J})\leq\frac{55}{432}s^2$, 
$\|h^{1,0}\|_{\infty}\leq \frac{2}{3}s$,
$\|h^{1,1}\|_1=\frac{1}{2}s(1+2^{-J})\leq \frac{3}{4}s$,
$\|(h^{1,0})'\|_1\leq\|h^{0,0}\|_1+\|h^{1,1}\|_1\leq \frac{5}{4}s$,
$\|h^{1,0}\|_{\mathrm{sp}}\leq\max\{\frac{55}{216}s^2,\frac{23}{4}s\}$,
$\|h^{0,2}\|_1=\frac{12}{s}(1+2^J)\leq 18\cdot 2^J$,
$\|h^{2,0}\|_1=\frac{19}{432}s^3(1-2^{3J})\leq\frac{19}{432}s^3$,
$\|h^{2,0}\|_\infty\leq \frac{4}{9}s^2$,
$\|h^{2,1}\|_1=\frac{55}{216}s^2(1+2^{-2J})\leq\frac{275}{864}s^2$,
$\|(h^{2,0})'\|_1\leq2\|h^{1,0}\|_1+\|h^{2,1}\|_1\leq \frac{55}{96}s^2$, and 
$\|h^{2,0}\|_{\mathrm{sp}}\leq \max\{\frac{19}{216}s^3,\frac{293}{96}s^2\}$.
%
%
%
%
%
%
%
%
Arguing as above, we obtain 

\begin{align}
&\sup_{\phi,w}|f(\phi,w)|(1+|w|^2)^{\frac{\eta}{2}}\leq\\
&\textcolor{black}{\leq2^{\frac{\eta}{2}}\max\left\{\max\{s,9\}, \textcolor{white}{2^{2^{2^{2}}}}\right.}\\
&\textcolor{black}{~~~~~~~~~~~~~~~~~~~~
2^{(\eta-1)J}\left(2+\max\{\tfrac{55}{216}s^2,\tfrac{23}{4}s\}
\right),}\\
&\textcolor{black}{~~~~~~~~~~~~~~~~~~~\left.
2^{(\eta-2)J}\left(\tfrac{3}{4}s+18\cdot 2^J+\tfrac{s}{2}+\max\{\tfrac{19}{216}s^3,\tfrac{293}{96}s^2\}\right)\right\}\leq}\\
&\textcolor{black}{\leq2^{\frac{\eta}{2}} 2^{(\eta-1)J} \max\left\{s,9,2+\max\{\tfrac{23}{4}s,\tfrac{55}{216}s^2\},18+\tfrac{5}{8}s+\max\{\tfrac{293}{192}s^2,\tfrac{19}{432}s^3\}\right\}\leq}\\
&\leq K_{L}(s)\, 
2^{(\eta-1)J},\label{where-C_3L-appears-1}
\end{align}
where
\begin{align}
K_{L}(s)=36+\tfrac{5}{4}s+\max\left\{\tfrac{293}{96}s^2,\tfrac{19}{216}s^3\right\}.
\label{piecewise-formula-C3L-new}
\end{align}
The argument for $\chi_{s,L}^{(J)}(-\cdot)$ is the same as for $\chi_{s,L}^{(J)}$.
\end{proof}

\begin{remark}\label{remark-ratio-C3-C3L}
An elementary but tedious analysis shows that for all $s\geq1$ we have
\begin{align}
2^{\ha}\leq\frac{K(s)}{K_{L}(s)}\leq2^3.\label{inequalities-ratio-C3andC3L}
\end{align}
\end{remark}

\subsection{Divergence of expected displacement for harmonic quantum oscillator with piecewise initial data}

 We illustrate here an application of Theorem \ref{theorem-uniform-bound-f(t,x)-new}. The equation \eqref{quantum-harmonic-oscillator} describes a 1-dimensional quantum harmonic oscillator. 
If we interpret $x\mapsto |f(t,x)|^2$ as a probability density, then we can consider the $m$-th moment of the displacement of the quantum particle at time $t$, i.e. $\langle x^m\rangle=\displaystyle\int_\R x^m |f(t,x)|^2\,\de x$. It is well known that the regularity of the initial condition $f$ determines the decay of $f(t,\cdot)$. If, for instance, $f\in\mathcal{S}(\R)$, then $f(t,\cdot)\in\mathcal{S}(\R)$  and hence all moments $\langle x^m\rangle$ are finite. 
However, if the initial condition is only piecewise $C^1$ with finite jumps (such as $f=\chi$), then $\langle x\rangle$ diverges because $|f(t,x)|^2$ decays too slowly as $|x|\to\infty$. In Theorem \ref{theorem-uniform-bound-f(t,x)} below we consider smooth approximations to such initial data, and we  precisely quantify the rate at which the divergence in the expected displacement occurs as the approximation increases.

 For the rest of this section,  $f_\epsilon$ will not denote $\mathscr F_\epsilon f$, and we will instead use it to denote a sequence of functions approaching $f$ in the suitable sense as $\epsilon\to0^+$.

\begin{defin}\label{def-convergence-with-bounded-gradient}
A sequence $(f_\epsilon)_{1\geq\epsilon>0}$ of $C^1$ functions converging in $L^1(\R)$ to a piecewise $C^1$ function $f_0\in C^1_{\mathrm{pw}}$ as $\epsilon\to0^+$ is said to \emph{converge with bounded gradient} if there exists a constant $C>0$ such that for all $\epsilon>0$ we have
$\|f_\epsilon\|_1\leq C$, 
$\|f_\epsilon'\|_1\leq C$, and
$\|f_\epsilon''\|_1\leq C\epsilon^{-1}$.
That is to say, the derivative of the functions $f_\epsilon$ will diverge as $\epsilon\to0$ to approximate the piecewise function $f_0$, but such divergence must occur in proportionally small intervals.
\end{defin}
 This definition can be motivated as follows. Let $f_0\in C^1_{\mathrm{pw}}$ be piecewise $C^1$ with finitely many jump discontinuities at the points $x_j$, and construct a sequence $f_\epsilon(x)=f_0(x)+\sum_{j}g_j((x-x_j)/ 
\epsilon)$ where the $g_j$'s are compactly supported functions in $[-1,1]$ chosen in order to make $f_\epsilon$ differentiable. By direct computation one can check that $f_\epsilon\to f_0$ in $L^1$ with bounded gradient as $\epsilon\to0^+$.
\begin{theorem}\label{theorem-uniform-bound-f(t,x)}
Let $f_0\in C^1_{\mathrm{pw}}$  be a bounded, compactly supported function with finitely many jump discontinuities. Let $(f_\epsilon)_{\epsilon>0}$ be a sequence of bounded, compactly supported $C^1$ functions converging to $f_0$ with bounded gradient. Let $f_\epsilon(t,x)$ be the solution to \eqref{quantum-harmonic-oscillator} with initial condition $f_\epsilon$. Then for all $1<\eta \leq 2$ and all sufficiently small $\epsilon$ there exists a constant $C'>0$ such that 
\begin{align}
|f_\epsilon(t,x)|\leq C'(1+|x|^2)^{-\frac{\eta}{2}}\epsilon^{-(\eta-1)}\label{statement-theorem-uniform-bound-f(t,x)}
\end{align}
uniformly for all $t$.
\end{theorem}
To illustrate the above theorem, fix $1<\eta\leq2$. As in \eqref{def-kappa_eta(f)}-\eqref{def-S_eta}, we think of $\eta$ as a regularity parameter. Observe that \eqref{statement-theorem-uniform-bound-f(t,x)} means that $f_\epsilon\in\mathcal{S}_\eta(\R)$ for  every $\epsilon>0$ and provides a bound $\kappa_\eta(f_\epsilon)=O(\epsilon^{1-\eta})$. 
This also  implies that for  every fixed $\epsilon>0$ the probability density $x\mapsto|f_\epsilon(t,x)|^2$ decays sufficiently fast as $|x|\to\infty$ so that the expected displacement $\int_\R x|f_\epsilon(t,x)|^2\de x$ is finite (and in fact uniformly bounded for all times). On the other hand, as mentioned above, the integral $\int_\R x|f_0(t,x)|^2\de x$ diverges for all $t>0$ due to the lack of regularity of $f_0$. These two contrasting behaviours are reconciled by the theorem: the smaller $\epsilon>0$ is, the closer $f_\epsilon$ is to $f_0$, and \eqref{statement-theorem-uniform-bound-f(t,x)} shows how the uniform bound for the expected displacement blows up to infinity as $\epsilon\to0^+$. 


\begin{proof}[Proof of Theorem \ref{theorem-uniform-bound-f(t,x)}] In order to apply Theorem \ref{theorem-uniform-bound-f(t,x)-new} to the function $f=f_\epsilon$, we need to consider the various norms featured in \eqref{statement-theorem-uniform-bound-f(t,x)-new}. We have
\begin{equation}
\| h^{0,0} \|_{\mathrm{sp}} =  \max \{ 2\| f_\epsilon \|_1, 3(\|f_\epsilon\|_\infty + \|f_\epsilon'\|_1) \} \leq K_1
\end{equation}
with $K_1 =3(M+C)$, $C$ as in Definition \ref{def-convergence-with-bounded-gradient} and $M = \sup_{0< \epsilon<1} \|f_\epsilon\|_{\infty}$. Also,
\begin{equation}
\| h^{0,1}\|_1 + \|h^{1,0} \|_{\mathrm{sp}} \leq C + \max \{2 AC , 3(AM +  C + AC) \} \leq K_2 
\end{equation}
where $K_2 = 4C +3A(C+M)$ and $A$ is such that $\operatorname{supp} f_\epsilon \subset [-A,A]$ for all $0<\epsilon<1$. Finally,
\begin{align}
&\| h^{1,1} \|_1 + \|h^{0,2}\|_1 + \|h^{0,0}\|_1 + \|h^{2,0}\|_{\mathrm{sp}} \\
&\leq A C + C \epsilon^{-1} + C + \max \{2 A^2 C,3( A^2 M + 2A C + A^2 C)\} \leq \epsilon^{-1} K_3
\end{align}
with $K_3 = 2 C + 3 A^2 (C+M) + 7 AC$. Using $b=\epsilon^{-1}$ in Theorem \ref{theorem-uniform-bound-f(t,x)-new}, we obtain Theorem \ref{statement-theorem-uniform-bound-f(t,x)} with $C'=2^{\frac{\eta}{2}} \max \{ K_1, K_2, K_3\}$.
\end{proof}

%% file: tail-approximations-10.tex
\section{Tail approximations}\label{section_tail_approximations}
Recall \eqref{preparing-Theta_chi_aTheta_chi_b-for-union-bound-1}--\eqref{preparing-Theta_chi_aTheta_chi_b-for-union-bound-2}, in which we have  approximated $\Theta_{\chi}\overline{\Theta_{\chi_b}}(g)$ by
 $\Theta_{\chi^{(J)}}\overline{\Theta_{\chi_b^{(J)}}}(g)$. Let  $\mu=\mu_{\GamG}^{\bn}$ in the rational case and $\mu=\mu_{\GamG}$ in the irrational case. In this section:
 \begin{itemize}
 \item[(i)] we first estimate how much the tail probability $\mu\left(|\Theta_{\chi}\overline{\Theta_{\chi_b}}|>R^2\right)$
 differs from $\mu\left(|\Theta_{\chi^{(J)}} \overline{\Theta_{\chi^{(J)}_b}}| >R^2(1\pm\delta)^2\right)$ for suitable choices of parameters $J\geq4$ and $\delta>0$ (see Lemmata \ref{smooth_approx_rat} and \ref{smooth_approx_irr}),
 \item[(ii)] then we give an asymptotic for $\mu\!\left(|\Theta_{\chi^{(J)}} \overline{\Theta_{\chi^{(J)}_b}}| >R^2(1\pm\delta)^2\right)$ which includes explicit error terms (see Lemmata \ref{lemma-final-approximations-rat} and \ref{lemma-final-approximations-irr}).  
 \end{itemize} 


\subsection{Rational case}
\begin{lemma}\label{smooth_approx_rat}
Fix $1<\eta \leq 2$, $b\geq1$, and let $K_{L}(b)$ as in Proposition \ref{prop-key-estimate-kappa_eta(chi_s^J)}. Then for all $J\geq4$, $R^2 > 64 \, C_\eta K_L(b)^2 2^{2(\eta - 1)J}$ and $2^{-\tfrac{J-1}{2}} <\delta< \tfrac12$ we have 
\begin{align}
\mu_{\GamG}^{\bn}\left(|\Theta_{\chi}\overline{\Theta_{\chi_b}}|>R^2\right)-\mu_{\GamG}^{\bn}\left(|\Theta_{\chi^{(J)}} \overline{\Theta_{\chi^{(J)}_b}}| >R^2(1-\delta)^2\right) &= O\!\left(\frac{1}{R^4 \delta^2 2^{J}}\right) \text{ and }\label{statement-smooth_approx_rat-1}\\
\mu_{\GamG}^{\bn}\left(|\Theta_{\chi^{(J)}} \overline{\Theta_{\chi^{(J)}_b}}| >R^2(1+\delta)^2\right) - \mu_{\GamG}^{\bn}\left(|\Theta_{\chi}\overline{\Theta_{\chi_b}}|>R^2\right) &= O\!\left(\frac{1}{R^4 \delta^2 2^{J}}\right),\label{statement-smooth_approx_rat-2}
\end{align}
where the constants implied by the $O$-notation in \eqref{statement-smooth_approx_rat-1}-\eqref{statement-smooth_approx_rat-2} can both be taken to be $ 2^{17} \| \chi_{b,L} \|_{\mathrm{sp}}^4$ (see \eqref{chi_s=sums-of-Deltas00}, Lemma \ref{lem-coarse-stationary-phase-new}, and \eqref{sp-norm-of-chi-s,L}). 
\end{lemma}

\begin{proof}
Set $S = \{(k,\ell)\in \Z_{\geq0}^2:\:\:(k,\ell)\neq(0,0)\}$. 
Note that $1 = (1 - \delta)^2 + 4\frac{(1-\delta)^2}{4}\displaystyle\sum_{(k,\ell) \in S} \delta^{k +\ell}$. We use \eqref{preparing-Theta_chi_aTheta_chi_b-for-union-bound-1}--\eqref{preparing-Theta_chi_aTheta_chi_b-for-union-bound-2} and a   union bound to obtain
\begin{align}
&\mu_{\GamG}^{\bn}\left(|\Theta_{\chi}\overline{\Theta_{\chi_b}}|>R^2\right) - \mu_{\GamG}^{\bn}\left(|\Theta_{\chi^{(J)}}\overline{\Theta_{\chi_b^{(J)}}}| > R^2 (1 - \delta)^2\right) \leq \label{rational-union-bound-0}\\
&\leq \: \sum_{(k,\ell) \in S}\mu_{\GamG}^{\bn}\left(2^{-{(k + \ell)J\over 2}}|\calLL^{(J)}_{k, \ell}| > R^2 \tfrac{(1-\delta)^2}{4}\delta^{k + \ell}\right)+ \label{rational-union-bound-1}\\
&+\sum_{(k,\ell) \in S}\mu_{\GamG}^{\bn}\left(2^{-{(k + \ell)J\over 2}}|\calLR^{(J)}_{k, \ell}| > R^2 \tfrac{(1-\delta)^2}{4}\delta^{k + \ell}\right) +\label{rational-union-bound-2}\\
& +\sum_{(k,\ell) \in S}\mu_{\GamG}^{\bn}\left(2^{-{(k + \ell)J\over 2}}|\calRL^{(J)}_{k, \ell}| > R^2 \tfrac{(1-\delta)^2}{4}\delta^{k + \ell}\right) +\label{rational-union-bound-3}\\
&+ \sum_{(k,\ell) \in S}\mu_{\GamG}^{\bn}\left(2^{-{(k + \ell)J\over 2}}|\calRR^{(J)}_{k, \ell}| > R^2 \tfrac{(1-\delta)^2}{4}\delta^{k + \ell}\right). \label{rational-union-bound-4}
\end{align}
We focus on estimating the sum in \eqref{rational-union-bound-1}, since \eqref{rational-union-bound-2}--\eqref{rational-union-bound-4} are similar.  First observe that another union bound and the invariance of the measure imply, for every $(k,\ell)\in S$, 
\begin{align}
&\mu_{\GamG}^{\bn}\left(2^{-{(k + \ell)J\over 2}}|\calLL^{(J)}_{k, \ell}| > R^2 \tfrac{(1-\delta)^2}{4}\delta^{k + \ell}\right)\leq\\
&\leq\mu_{\GamG}^{\bn}\left\{g\in\mathcal F_\Gamma:\: \left|\Theta_{\chi_{L}^{(J)}}(\tg\wtPhi^{-(2\log2)k J})\right|>R\tfrac{1-\delta}{2}(\delta 2^{\frac{J}{2}})^{\frac{k+\ell}{2}}\right\}+\label{proof-lemma-smooth_approx_rat-1m1}\\
&~~~~+\mu_{\GamG}^{\bn}\left\{g\in\mathcal F_\Gamma:\: \left|\Theta_{\chi_{b,L}^{(J)}}(\tg\wtPhi^{-(2\log2)\ell J})\right|>R\tfrac{1-\delta}{2}(\delta 2^{\frac{J}{2}})^{\frac{k+\ell}{2}}\right\}=\label{proof-lemma-smooth_approx_rat-1m2}\\
&=\mu_{\GamG}^{\bn}\left\{g\in\mathcal F_\Gamma:\: \left|\Theta_{\chi_{L}^{(J)}}(\tg)\right|>R\tfrac{1-\delta}{2}(\delta 2^{\frac{J}{2}})^{\frac{k+\ell}{2}}\right\}+\label{proof-lemma-smooth_approx_rat-1}\\
&~~~~+\mu_{\GamG}^{\bn}\left\{g\in\mathcal F_\Gamma:\: \left|\Theta_{\chi_{b,L}^{(J)}}(\tg)\right|>R\tfrac{1-\delta}{2}(\delta 2^{\frac{J}{2}})^{\frac{k+\ell}{2}}\right\}.\label{proof-lemma-smooth_approx_rat-2}
\end{align}
In \eqref{proof-lemma-smooth_approx_rat-1m1}--\eqref{proof-lemma-smooth_approx_rat-2} use $\tg$ such that $\pi(\tg)=g$ as discussed in Section \ref{subsection-writing-product-of-Thetas}. To apply Lemma \ref{L2.2} to \eqref{proof-lemma-smooth_approx_rat-1} and \eqref{proof-lemma-smooth_approx_rat-2} we need
\begin{equation}
R \tfrac{1-\delta}{2} ( \delta 2^{\frac{J}{2}} )^{\tfrac{k+\ell}{2}} > 2^{\tfrac{\eta}{2}} C_\eta^{\ha} \max \{ \kappa_\eta(\chi_{L}^{(J)}),  \kappa_\eta(\chi_{b,L}^{(J)}) \} . \label{proof-lemma-smooth_approx_rat-3}
\end{equation}
By Proposition \ref{prop-key-estimate-kappa_eta(chi_s^J)} we can write
$\kappa_\eta(\chi_{s,L}^{(J)})\leq K_{L}(s)\, 2^{(\eta-1)J}$ for $s\in\{1,b\}$ and thus, since $\eta\leq2$ and $\delta < \ha$, \eqref{proof-lemma-smooth_approx_rat-3} is implied by the inequality
\begin{equation}
R > 8 \,C_\eta^{\ha}  \, K_L(s) \, 2^{(\eta-1) J}.
\end{equation}
Under this assumption, Lemma \ref{L2.2} gives
\begin{align}
\mu_{\GamG}^{\bn}\left\{g\in\mathcal F_\Gamma:\: \left|\Theta_{\chi_{L}^{(J)}}(\tg)\right|>R\tfrac{1-\delta}{2}(\delta 2^{\frac{J}{2}})^{\frac{k+\ell}{2}}\right\}=\frac{\frac{2}{\pi^2}D_{\mathrm{rat}}(\chi_{L}^{(J)},\chi_{L}^{(J)})}{R^4\tfrac{(1-\delta)^4}{16}(\delta^2 2^J)^{k+\ell}}(1+\mathscr{E}),\label{proof-lemma-smooth_approx_rat-4}
\end{align}
where
$|\mathscr{E}|\leq 3 (2^\eta C_\eta)^\eta \kappa_\eta(\chi_{L}^{(J)})^{2\eta}\left(R\tfrac{1-\delta}{2}(\delta 2^{\frac{J}{2}})^{\frac{k+\ell}{2}}\right)^{-2\eta}$. By \eqref{proof-lemma-smooth_approx_rat-3} we obtain $|\mathscr{E}| \leq 3$, and by Lemma \ref{lem-coarse-stationary-phase-new} we have
\begin{align}
D_{\mathrm{rat}}(\chi_{L}^{(J)},\chi_{L}^{(J)})=\int_0^\pi\left|\left(\chi_{L}^{(J)}\right)_\phi(0)\right|^4\de\phi \leq \pi   \|\chi_{L}^{(J)}\|_{\mathrm{sp}}^4.
\end{align}
Since $\delta<\tfrac12$ we have $(1-\delta)^4 > 2^{-4}$, and thus \eqref{proof-lemma-smooth_approx_rat-4} can be bounded above by
\begin{equation}
\frac{ 2^{11} \| \chi_L^{(J)} \|_\mathrm{sp}^4}{\pi R^4 (\delta^2 2^J)^{k+\ell} }.
\end{equation}
An identical argument allows us to bound the quantity in \eqref{proof-lemma-smooth_approx_rat-2} from above by 
\begin{equation}
\frac{ 2^{11} \| \chi_{b,L}^{(J)} \|_\mathrm{sp}^4}{\pi R^4 (\delta^2 2^J)^{k+\ell} }.
\end{equation}
Now observe that for $q>2$ we have $\sum_{(k,\ell)\in S}q^{-(k+\ell)}=\frac{1}{q}\!\left(\frac{2q^2-q}{(q-1)^2}\right)\leq\frac{6}{q}$. This bound can be used for $q=\delta^2 2^{J}$ since $\delta>2^{-\frac{J-1}{2}}$ by hypothesis. Therefore,  since $\|\chi_L^{(J)}\|_{\mathrm{sp}} <\|\chi_{b,L}^{(J)}\|_{\mathrm{sp}}$, the sum in \eqref{rational-union-bound-1} is bounded above by
\begin{equation}
\frac{2^{13} \| \chi_{b,L}^{(J)} \|_{\mathrm{sp}}^4 }{  R^4 \delta^2 2^J}.
\end{equation}
We can obtain the same upper bound for each of the sums in \eqref{rational-union-bound-2}--\eqref{rational-union-bound-4} and therefore obtain \eqref{statement-smooth_approx_rat-1} with an implied constant equal to $2^{15} \|\chi_{b,L}^{(J)}\|_{\mathrm{sp}}^4$.
To prove \eqref{statement-smooth_approx_rat-2}, we use a slightly different union bound, arising from the identity 
\begin{align}
(1+\delta)^2&=1+(1-\delta)^2\left(\sum_{k\geq0}\sum_{\ell\geq1}\delta^{k+\ell}+\sum_{k\geq1}\sum_{\ell\geq0}\delta^{k+\ell}+\sum_{k\geq1}\sum_{\ell\geq1}\delta^{k+\ell}\right)=\\
&=1+4\frac{(1-\delta)^2}{4}\sum_{(k,\ell)\in S}c_{k,\ell}\delta^{k+\ell},
\end{align}
where $c_{k,\ell}=1$ if $k=0$ or $\ell=0$ and $c_{k,\ell}=3$ otherwise. Now \eqref{preparing-Theta_chi_aTheta_chi_b-for-union-bound-1}--\eqref{preparing-Theta_chi_aTheta_chi_b-for-union-bound-2} and a union bound yield
\begin{align}
&\mu_{\GamG}^{\bn}\left(|\Theta_{\chi}\overline{\Theta_{\chi_b}}|>R^2 (1+\delta)^2\right) - \mu_{\GamG}^{\bn}\left(|\Theta_{\chi^{(J)}}\overline{\Theta_{\chi_b^{(J)}}}| > R^2\right) \leq \label{rational-union-bound-5}\\
&\leq \: \sum_{(k,\ell) \in S}\mu_{\GamG}^{\bn}\left(2^{-{(k + \ell)J\over 2}}|\calLL^{(J)}_{k, \ell}| > c_{k,\ell}R^2 \tfrac{(1-\delta)^2}{4}\delta^{k + \ell}\right)+ \label{rational-union-bound-6}\\
&+\sum_{(k,\ell) \in S}\mu_{\GamG}^{\bn}\left(2^{-{(k + \ell)J\over 2}}|\calLR^{(J)}_{k, \ell}| > c_{k,\ell}R^2 \tfrac{(1-\delta)^2}{4}\delta^{k + \ell}\right) +\label{rational-union-bound-7}\\
& +\sum_{(k,\ell) \in S}\mu_{\GamG}^{\bn}\left(2^{-{(k + \ell)J\over 2}}|\calRL^{(J)}_{k, \ell}| > c_{k,\ell}R^2 \tfrac{(1-\delta)^2}{4}\delta^{k + \ell}\right) +\label{rational-union-bound-8}\\
&+ \sum_{(k,\ell) \in S}\mu_{\GamG}^{\bn}\left(2^{-{(k + \ell)J\over 2}}|\calRR^{(J)}_{k, \ell}| > c_{k,\ell}R^2 \tfrac{(1-\delta)^2}{4}\delta^{k + \ell}\right). \label{rational-union-bound-9}
\end{align}
Arguing as before, the sum in \eqref{rational-union-bound-6} can be bounded above by
\begin{align}
\frac{2^{15} \| \chi_{b,L}^{(J)} \|_{\mathrm{sp}}^4 }{ R^4\delta^2 2^J},
\end{align}
and similarly for \eqref{rational-union-bound-7}--\eqref{rational-union-bound-9}. Using $\|\chi_{b,L}^{(J)}\|_{\mathrm{sp}} \leq \|\chi_{b,L}\|_{\mathrm{sp}}$ (see Remark \ref{sp-norms-for-chi-functions}) we obtain \eqref{statement-smooth_approx_rat-2} with an implied constant equal to $2^{17}  \| \chi_{b,L}\|_{\mathrm{sp}}^4 $.
\end{proof}

\begin{lemma}\label{lemma-final-approximations-rat}
Fix $1<\eta\leq2$, $b\geq1$, and let $K(b)$ be as in Proposition \ref{prop-key-estimate-kappa_eta(chi_s^J)}. Then for all $J \geq 1$, $R^2> 16 \, C_\eta K(b)^2 2^{2(\eta-1)J}$ and $0<\delta<\ha$ we have
\begin{align}
&\mu_{\GamG}^{\bn}\left\{ |\Theta_{\chi^{(J)}}\overline{\Theta_{\chi_{b}^{(J)}}}|>R^2(1\pm\delta)^2\right\}=\nonumber\\
&=\frac{2D_{\mathrm{rat}}(\chi,\chi_b)}{\pi^2 R^4}\left(1+O(\delta)\right)\left(1+O_\eta(2^{2\eta(\eta-1)J}R^{-2\eta})\right)\left(1+O(2^{-J})\right), \label{lemma-final-approximations-rat-statement}
\end{align}
where the constants implied by the O-notations can be taken to be $30$,   $2^{10} C_\eta^2 K(b)^4$, and 
$2^5 b \| \chi_b\|_{\mathrm{sp}}^3$ (see \eqref{sp-norm-of-chi-s}) respectively.
\end{lemma}
\begin{proof}
In order to apply Lemma \ref{L2.2} we require
\begin{equation}
R^2(1\pm \delta)^2 > 2^\eta C_\eta \kappa_\eta(\chi^{(J)})\kappa_\eta(\chi_{b}^{(J)}) \label{proof-lemma-final-approximations-rat-0}.
\end{equation}
By Proposition \ref{prop-key-estimate-kappa_eta(chi_s^J)} we have
$\kappa_\eta(\chi^{(J)})\kappa_\eta(\chi_{b}^{(J)})\leq K(1)K(b) 2^{2(\eta-1)J}$. Since $\eta\leq2$, $\delta < \tfrac12$ and $K(1) \leq K(b)$, \eqref{proof-lemma-final-approximations-rat-0} is implied by the inequality
\begin{equation}
R^2 >   16 \, C_\eta  K(b)^2 2^{2(\eta-1) J}.\label{inequality-after-using-K1<Kb-rat}
\end{equation}
We obtain
\begin{align}
\mu_{\GamG}^{\bn}\left\{ |\Theta_{\chi^{(J)}}\overline{\Theta_{\chi_{b}^{(J)}}}|>R^2(1\pm\delta)^2\right\}=\frac{\frac{2}{\pi^2}D_{\mathrm{rat}}(\chi^{(J)},\chi_{b}^{(J)})}{R^4(1\pm\delta)^4}\left(1+\mathscr{E}_1\right),\label{proof-lemma-final-approximations-rat-1}
\end{align}
where
\begin{align}
|\mathscr{E}_1| &\leq 3 (2^\eta C_\eta  \kappa_\eta(\chi^{(J)})\kappa_\eta(\chi_{b}^{(J)}))^\eta \left(R(1\pm\delta)\right)^{-2\eta} \\
&\leq 3 (2^{\eta+2} C_\eta K(b)^2)^\eta \, 2^{2\eta(\eta-1) J} R^{-2\eta} \\
&\leq 2^{10} C_\eta^2 K(b)^4\,  2^{2\eta(\eta-1) J} R^{-2\eta}.\label{proof-lemma-final-approximations-rat-2}
\end{align} 
This yields the second implied constant in the statement of the Lemma. For the first implied constant, note that for $0<\delta<\ha$ we have 
\begin{align}
\frac{1}{(1\pm\delta)^4}\leq\frac{1}{(1-\delta)^4}=(1+\mathscr{E}_2), \hspace{.5cm}|\mathscr{E}_2|\leq30 \delta.\label{proof-lemma-final-approximations-rat-3}
\end{align}
The third implied constant is more laborious. We have the estimate
\begin{align}
&\left|D_{\text{rat}}(\chi, \chi_b) - D_{\text{rat}}(\chi^{(J)}, \chi_b^{(J)})\right| \leq\\
&\leq \int_{0}^{\pi}\left|\left|\chi_{\phi}(0)\,(\chi_{b})_{\phi}(0)\right|^2 - \left|(\chi^{(J)})_{\phi}(0)\,(\chi_b^{(J)})_{\phi}(0)\right|^2 \right|\: \de\phi\leq.
\end{align}
By the triangle inequality this is bounded above by
\begin{align}
&\int_{0}^{\pi}\left|\left|\chi_{\phi}(0)\,(\chi_{b})_{\phi}(0)\right|^2 - \left|(\chi^{(J)})_{\phi}(0)\,(\chi_b)_{\phi}(0)\right|^2 \right|\: \de\phi\, +\\
&~~~+ \int_{0}^{\pi}\left|\left|(\chi^{(J)})_{\phi}(0)\,(\chi_{b})_{\phi}(0)\right|^2 - \left|(\chi^{(J)})_{\phi}(0)\,(\chi_b^{(J)})_{\phi}(0)\right|^2 \right|\: \de\phi\leq.
\end{align}
Now apply Lemma \ref{lem-coarse-stationary-phase-new} to $| (\chi_b)_\phi(0)|$ in the first term, and $| (\chi^{(J)})_\phi(0)|$ in the second term to obtain
\begin{align}
\left|D_{\text{rat}}(\chi, \chi_b) - D_{\text{rat}}(\chi^{(J)}, \chi_b^{(J)})\right|  &\leq \|\chi_b\|_{\mathrm{sp}}^2 \int_{0}^{\pi}\left|\left|\chi_{\phi}(0)\right|^2 - \left|(\chi^{(J)})_{\phi}(0)\right|^2\right|\de\phi \\
&~~~+ \|\chi^{(J)}\|_{\mathrm{sp}}^2\int_{0}^{\pi}\left|\left|(\chi_{b})_{\phi}(0)\right|^2 - \left|(\chi_{b}^{(J)})_{\phi}(0)\right|^2\right|\de\phi.
\end{align}
Factorising the integrands, applying Lemma \ref{lem-coarse-stationary-phase-new} twice more, and using the reverse triangle inequality gives
\begin{align}
&\left|D_{\text{rat}}(\chi, \chi_b) - D_{\text{rat}}(\chi^{(J)}, \chi_b^{(J)})\right|\leq\\
&\leq \|\chi_b\|_{\mathrm{sp}}^2 ( \|\chi\|_{\mathrm{sp}} + \| \chi^{(J)}\|_{\mathrm{sp}} )  \int_{0}^{\pi}\left|\chi_{\phi}(0) - (\chi^{(J)})_{\phi}(0)\right|\,\de\phi \label{proof-lemma-final-approximations-rat-5}\\
&~~~+ \|\chi^{(J)}\|_{\mathrm{sp}}^2 (\|\chi_b\|_{\mathrm{sp}}+\|\chi_b^{(J)}\|_{\mathrm{sp}} ) \int_{0}^{\pi}\left|(\chi_{b})_{\phi}(0) - (\chi_{b}^{(J)})_{\phi}(0)\right|\,\de\phi.\label{proof-lemma-final-approximations-rat-6}
\end{align}
Furthermore, since $\|\chi_s-\chi_s^{(J)}\|_1=s2^{-J}$, note that 
\begin{align}
\int_0^\pi \left| (\chi_s)_\phi(0)-(\chi_s^{(J)})_\phi(0)\right|\de \phi&=\int_0^\pi\left| (\chi_s-\chi_s^{(J)})_\phi(0)\right|\de\phi\leq\\
&\leq\int_{0}^{\pi}\frac{\de \phi}{\sqrt{\sin(\phi)}}\int_{\R}\left|(\chi_s-\chi_s^{(J)})(v)\right|\de v=
\frac{\sqrt{2}\pi^{\frac{3}{2}}s}{\Gamma(\frac{3}{4})^2}2^{-J}.\label{proof-lemma-final-approximations-rat-7}
\end{align}
By Remark \ref{sp-norms-for-chi-functions}, each of the norms is bounded above by $\| \chi_b \|_{\mathrm{sp}}$ and we obtain
\begin{equation}
\left|D_{\text{rat}}(\chi, \chi_b) - D_{\text{rat}}(\chi^{(J)}, \chi_b^{(J)})\right|\leq 2^5 \|\chi_b\|_{\mathrm{sp}}^3 b \, 2^{-J}.\label{proof-lemma-final-approximations-rat-8}
\end{equation}
Combining \eqref{proof-lemma-final-approximations-rat-1} with \eqref{proof-lemma-final-approximations-rat-2}, \eqref{proof-lemma-final-approximations-rat-3} and \eqref{proof-lemma-final-approximations-rat-8} gives the result.
\end{proof}

\subsection{Irrational case}
\begin{lemma}\label{smooth_approx_irr}
Fix $1<\eta\leq2$, $b\geq1$ and let $K_L(b)$ be as in Proposition \ref{prop-key-estimate-kappa_eta(chi_s^J)}. Then for all $J\geq 4$, $R^2 > 64 C_\eta K_L(b)^2 2^{2(\eta-1)J}$ and $2^{-\frac{J-1}{2}}<\delta<\ha$ we have 
\begin{align}
\mu_{\GamG}^{}\left(|\Theta_{\chi}\overline{\Theta_{\chi_b}}|>R^2\right)-\mu\left(|\Theta_{\chi^{(J)}} \overline{\Theta_{\chi^{(J)}_b}}| >R^2(1-\delta)^2\right) &= O\!\left(\frac{1}{R^6 \delta^3 2^{\frac{3}{2}J}}\right) \text{ and }\label{statement-smooth_approx_irr-1}\\
\mu_{\GamG}^{}\left(|\Theta_{\chi^{(J)}} \overline{\Theta_{\chi^{(J)}_b}}| >R^2(1+\delta)^2\right) - \mu\left(|\Theta_{\chi}\overline{\Theta_{\chi_b}}|>R^2\right) &= O\!\left(\frac{1}{R^6 \delta^3 2^{\frac{3}{2}J}}\right),\label{statement-smooth_approx_irr-2}
\end{align}
where the constants implied by the $O$-notation in \eqref{statement-smooth_approx_rat-1}-\eqref{statement-smooth_approx_rat-2} can both be taken to be $2^{22} \|\chi_{b,L}\|_{\mathrm{sp}}^5$ (see \eqref{sp-norm-of-chi-s,L}).
\end{lemma}

\begin{proof}
We follow the same strategy as the proof of Lemma \ref{smooth_approx_rat} and we omit several details.
Set $S = \{(k,\ell)\in \Z_{\geq0}^2:\:\:(k,\ell)\neq(0,0)\}$. 
We use \eqref{preparing-Theta_chi_aTheta_chi_b-for-union-bound-1}--\eqref{preparing-Theta_chi_aTheta_chi_b-for-union-bound-2} and a   union bound to obtain
\begin{align}
&\mu_{\GamG}^{}\left(|\Theta_{\chi}\overline{\Theta_{\chi_b}}|>R^2\right) - \mu_{\GamG}^{}\left(|\Theta_{\chi^{(J)}}\overline{\Theta_{\chi_b^{(J)}}}| > R^2 (1 - \delta)^2\right) \leq \label{irrational-union-bound-0}\\
&\leq \: \sum_{(k,\ell) \in S}\mu_{\GamG}^{}\left(2^{-{(k + \ell)J\over 2}}|\calLL^{(J)}_{k, \ell}| > R^2 \tfrac{(1-\delta)^2}{4}\delta^{k + \ell}\right)+ \label{irrational-union-bound-1}\\
&+\sum_{(k,\ell) \in S}\mu_{\GamG}^{}\left(2^{-{(k + \ell)J\over 2}}|\calLR^{(J)}_{k, \ell}| > R^2 \tfrac{(1-\delta)^2}{4}\delta^{k + \ell}\right) +\label{irrational-union-bound-2}\\
& +\sum_{(k,\ell) \in S}\mu_{\GamG}^{}\left(2^{-{(k + \ell)J\over 2}}|\calRL^{(J)}_{k, \ell}| > R^2 \tfrac{(1-\delta)^2}{4}\delta^{k + \ell}\right) +\label{irrational-union-bound-3}\\
&+ \sum_{(k,\ell) \in S}\mu_{\GamG}^{}\left(2^{-{(k + \ell)J\over 2}}|\calRR^{(J)}_{k, \ell}| > R^2 \tfrac{(1-\delta)^2}{4}\delta^{k + \ell}\right). \label{irrational-union-bound-4}
\end{align}
Each term of the sum in \eqref{irrational-union-bound-1} can be estimated as
\begin{align}
&\mu_{\GamG}^{}\left(2^{-{(k + \ell)J\over 2}}|\calLL^{(J)}_{k, \ell}| > R^2 \tfrac{(1-\delta)^2}{4}\delta^{k + \ell}\right)\leq\\
&
\leq\mu_{\GamG}^{}\left\{g\in\mathcal F_\Gamma:\: \left|\Theta_{\chi_{L}^{(J)}}(\tg)\right|>R\tfrac{1-\delta}{2}(\delta 2^{\frac{J}{2}})^{\frac{k+\ell}{2}}\right\}+\label{proof-lemma-smooth_approx_irr-1}\\
&~~~~+\mu_{\GamG}^{}\left\{g\in\mathcal F_\Gamma:\: \left|\Theta_{\chi_{b,L}^{(J)}}(\tg)\right|>R\tfrac{1-\delta}{2}(\delta 2^{\frac{J}{2}})^{\frac{k+\ell}{2}}\right\}.\label{proof-lemma-smooth_approx_irr-2}
\end{align}
In \eqref{proof-lemma-smooth_approx_irr-1}--\eqref{proof-lemma-smooth_approx_irr-2} we use $\tg$ such that $\pi(\tg)=g$ as discussed in Section \ref{subsection-writing-product-of-Thetas}. To apply Lemma \ref{L2.2_irrational} to each of these terms we need
\begin{equation}
R \tfrac{1-\delta}{2} (\delta 2^{\tfrac{J}{2}} )^{\tfrac{k+\ell}{2}} > 2^{\tfrac{\eta}{2}} C_\eta^{\ha} \max \{ \kappa_\eta(\chi_L^{(J)}),\kappa_\eta(\chi_{b,L}^{(J)}) \}. 
\end{equation}
By Proposition \ref{prop-key-estimate-kappa_eta(chi_s^J)}, and since $\eta\leq 2$ this is implied by the assumption
\begin{equation}
R> 8 C_\eta^{\ha} K_L(b)\, 2^{(\eta-1)J}.
\end{equation}
Under this assumption, Lemma \ref{L2.2_irrational} gives
\begin{align}
\mu_{\GamG}^{}\left\{g\in\mathcal F_\Gamma:\: \left|\Theta_{\chi_{L}^{(J)}}(\tg)\right|>R\tfrac{1-\delta}{2}(\delta 2^{\frac{J}{2}})^{\frac{k+\ell}{2}}\right\}=\frac{\frac{2}{\pi^2}D_{\mathrm{irr}}(\chi_{L}^{(J)},\chi_{L}^{(J)})}{R^6\tfrac{(1-\delta)^6}{64}(\delta^3 2^{\frac{3}{2}J})^{k+\ell}}(1+\mathscr{E}),\label{proof-lemma-smooth_approx_irr-4}
\end{align}
where
$|\mathscr{E}|\leq 7
$. Note that $\delta<\ha$ implies $(1-\delta)^{-6}\leq64$. Additionally,  Lemma \ref{lem-coarse-stationary-phase-new} and the unitarity of $f\mapsto f_\phi$ yields
\begin{align}
D_{\mathrm{irr}}(\chi_{L}^{(J)},\chi_{L}^{(J)})&=\int_0^\pi\int_{-\infty}^{\infty}\left|\left(\chi_{L}^{(J)}\right)_\phi\!(w)\right|^6\de w\,\de\phi\leq 
\|\chi_{L}^{(J)}\|_{\mathrm{sp}}^4\int_{0}^\pi\left\|\ \!\!\left(\chi_L^{(J)}\right)_\phi\right\|_2^2\,\de\phi=\\
&=\pi\|\chi_L^{(J)}\|_{\mathrm{sp}}^4 \, \|\chi_L^{(J)}\|_2^2 .
\end{align}
Since $|\chi_L^{(J)}(w)| \leq 1$, we have $\|\chi_L^{(J)}\|_2^2 \leq \|\chi_L^{(J)} \|_1 \leq \frac12 \| \chi_L^{(J)}\|_\mathrm{sp}$, and thus
\begin{align}
\mu_{\GamG}^{}\left\{g\in\mathcal F_\Gamma:\: \left|\Theta_{\chi_{L}^{(J)}}(\tg)\right|>R\tfrac{1-\delta}{2}(\delta 2^{\frac{J}{2}})^{\frac{k+\ell}{2}}\right\} < \frac{2^{14} \|\chi_L^{(J)}\|_{\mathrm{sp}}^5}{ R^6  ( \delta^3 2^{\frac{3}{2} J} )^{k+\ell} }.
\end{align}
An identical argument allows us to bound the quantity in \eqref{proof-lemma-smooth_approx_irr-2} from above by  $\frac{2^{14} \|\chi_{b,L}^{(J)}\|_{\mathrm{sp}}^5}{ R^6  ( \delta^3 2^{\frac32 J} )^{k+\ell} }$. Now observe that for $q>2^{\frac{3}{2}}$ we have $\sum_{(k,\ell)\in S}q^{-(k+\ell)} \leq\frac{4}{q}$. This bound can be used for $q=\delta^3 2^{\frac{3}{2}J}$ since $\delta>2^{-\frac{J-1}{2}}$ by hypothesis. Therefore, the sum in \eqref{irrational-union-bound-1} is bounded above by
\begin{equation}
\frac{2^{15} \|\chi_{b,L}^{(J)}\|_{\mathrm{sp}}^5}{R^6} \sum_{(k,\ell)\in S} \frac{1}{  ( \delta^3 2^{\frac{3}{2} J} )^{k+\ell} } \leq \frac{2^{17} \|\chi_{b,L}^{(J)}\|_{\mathrm{sp}}^5}{R^6 \delta^3 2^{\frac{3}{2} J} }.
\end{equation}
We can obtain the same upper bound for each of the sums in \eqref{irrational-union-bound-2}--\eqref{irrational-union-bound-4} and therefore obtain \eqref{statement-smooth_approx_irr-1} with an implied constant equal to $2^{19} \|\chi_{b,L}^{(J)}\|_{\mathrm{sp}}^5$. To prove \eqref{statement-smooth_approx_irr-2}, we can use a slightly different union bound as in the proof of \eqref{statement-smooth_approx_rat-2}. Combined with the upper bound $\|\chi_{b,L}^{(J)}\|_\mathrm{sp} \leq \|\chi_{b,L} \|_{\mathrm{sp}}$, this yields \eqref{statement-smooth_approx_irr-2} with an implied constant equal to $2^{21} \|\chi_{b,L}\|_{\mathrm{sp}}^5$.
\end{proof}

\begin{lemma}\label{lemma-final-approximations-irr}
Fix $1<\eta\leq2$, $b\geq1$, and let $K(b)$ be as in Proposition \ref{prop-key-estimate-kappa_eta(chi_s^J)}. Suppose that $R^2> 16 C_\eta K(b)^2 2^{2(\eta-1) J}$ and $0<\delta<\ha$. Then 
\begin{align}
&\mu_{\GamG}^{}\left\{ |\Theta_{\chi^{(J)}}\overline{\Theta_{\chi_{b}^{(J)}}}|>R^2(1\pm\delta)^2\right\}=\nonumber\\
&=\frac{2D_{\mathrm{irr}}(\chi,\chi_b)}{\pi^2 R^6}\left(1+O(\delta)\right)\left(1+O_\eta(2^{2\eta(\eta-1)J}R^{-2\eta} )\right)\left(1+O(2^{-\frac{J}{2}})\right), \label{lemma-final-approximations-irr-statement}
\end{align}
where the constants implied by the O-notations can be taken to be $126$,   $2^{11} C_\eta^2 K(b)^4$, and $8 \sqrt{b} \|\chi_b \|_{\mathrm{sp}}^5$ (see \eqref{sp-norm-of-chi-s}) respectively.
\end{lemma}
\begin{proof}
We follow the strategy of Lemma 3.16 in \cite{Cellarosi-Marklof}, but keep track of all implied constants.The argument is similar to the proof of Lemma \ref{lemma-final-approximations-rat} and omit several details. In order to apply Lemma \ref{L2.2_irrational} we require
\begin{equation}
R^2(1\pm\delta)^2 > 2^\eta C_\eta \kappa_\eta(\chi^{(J)}) \kappa_\eta (\chi_b^{(J)}).
\end{equation}
Applying Proposition \ref{prop-key-estimate-kappa_eta(chi_s^J)}, and using $\eta \leq 2$ and $\delta < \tfrac12$, this is implied by the assumption
\begin{equation}
R^2 > 16 C_\eta K(b)^2 2^{2(\eta-1) J}.
\end{equation}
Applying Lemma \ref{L2.2_irrational} we obtain
\begin{align}
\mu_{\GamG}^{}\left\{ |\Theta_{\chi^{(J)}}\overline{\Theta_{\chi_{b}^{(J)}}}|>R^2(1\pm\delta)^2\right\}=\frac{\frac{2}{\pi^2}D_{\mathrm{irr}}(\chi^{(J)},\chi_{b}^{(J)})}{R^6(1\pm\delta)^6}\left(1+\mathscr{E}_1\right),\label{proof-lemma-final-approximations-irr-1}
\end{align}
where
\begin{align}
|\mathscr{E}_1|&\leq 7(2^\eta C_\eta \kappa_\eta(\chi^{(J)})\kappa_\eta(\chi_{b}^{(J)}))^\eta \left(R(1\pm\delta)\right)^{-2\eta} \\
&\leq 7(2^{\eta+2}C_\eta K(b)^2)^{\eta} 2^{2\eta(\eta-1)J} R^{-2\eta}\\
&\leq 2^{11}C_\eta^2 K(b)^4\label{proof-lemma-final-approximations-irr-2}
\end{align} 
This yields the second implied constant in the statement of the Lemma. For the first implied constant, note that for $0<\delta<\ha$ we have 
\begin{align}
\frac{1}{(1\pm\delta)^6}\leq\frac{1}{(1-\delta)^6}=(1+\mathscr{E}_2), \hspace{.5cm}|\mathscr{E}_2|\leq126 \delta.\label{proof-lemma-final-approximations-irr-3}
\end{align}
For the third implied constant, note that
\begin{align}
&\left|D_{\text{irr}}(\chi, \chi_b) - D_{\text{irr}}(\chi^{(J)}, \chi_b^{(J)})\right|\leq\\
&\leq\int_0^{\pi}\left|\int_{-\infty}^{\infty} \left(|\chi_\phi(w)(\chi_b)_\phi(w)|^3-|(\chi^{(J)})_\phi(w)(\chi_b^{(J)})_\phi|^3\right)\de w\right|\de\phi=\\
&=\int_{0}^{\pi}\left|\left\|(\chi_\phi(\chi_b)_\phi)^3\right\|_1-\left\|((\chi^{(J)})_\phi(\chi_b^{(J)})_\phi)^3 \right\|_1\right|\de\phi\leq
\\
&\leq\int_{0}^{\pi}\left\| (\chi_\phi(\chi_b)_\phi)^3-((\chi^{(J)})_\phi(\chi_b^{(J)})_\phi)^3\right\|_1\de\phi,\label{proof-lemma-final-approximations-irr-5a}
\end{align}
where in \eqref{proof-lemma-final-approximations-irr-5a} we have used the reverse triangle inequality for the $1$-norm. 
Now, using the identity $x^3-y^3=(x-y)(x^2+xy+y^2)$ and the bound 
\begin{align}
&\left|(\chi_\phi(w)(\chi_b)_\phi(w))^2+\chi_\phi(w)(\chi_b)_\phi(w)(\chi^{(J)})_\phi(w)(\chi_b^{(J)})_\phi(w)+((\chi^{(J)})_\phi(w)(\chi_b^{(J)})_\phi(w))^2\right|\leq\nonumber\\
&\leq 3\|\chi_b \|_{\mathrm{sp}}^4
\end{align}
coming from Lemma \ref{lem-coarse-stationary-phase-new} and Remark \ref{sp-norms-for-chi-functions}, we see that \eqref{proof-lemma-final-approximations-irr-5a} implies
\begin{align}
&\left|D_{\text{irr}}(\chi, \chi_b) - D_{\text{irr}}(\chi^{(J)}, \chi_b^{(J)})\right|\leq\\
&\leq 3 \|\chi_b \|_{\mathrm{sp}}^4\int_{0}^{\pi}\left\| \chi_\phi(\chi_b)_\phi-(\chi^{(J)})_\phi(\chi_b^{(J)})_\phi\right\|_1\de\phi.\label{proof-lemma-final-approximations-irr-5b}
\end{align}
Now we add and subtract $\chi_\phi(\chi_b^{(J)})_\phi$ in \eqref{proof-lemma-final-approximations-irr-5b} and use the triangle inequality, followed by the Cauchy-Schwarz inequality, to obtain
\begin{align}
&\left|D_{\text{irr}}(\chi, \chi_b) - D_{\text{irr}}(\chi^{(J)}, \chi_b^{(J)})\right|\leq\\
&\leq 3 \|\chi_b \|_{\mathrm{sp}}^4 \int_{0}^{\pi}\left(\|\chi_\phi\|_2\|(\chi_b)_\phi-(\chi_b^{(J)})_\phi\|_2+\|(\chi_b^{(J)})_\phi\|_2\|\chi_\phi-(\chi^{(J)})_\phi \|_2\right)\de\phi.
\end{align}
The unitarity of the Shale-Weil representation and the bounds $\|\chi_\phi\|_2^2 \leq \frac12 \|\chi_b\|_\mathrm{sp}$, $\|(\chi_b^{(J)})_\phi\|_2^2 \leq  \frac12 \|\chi_b\|_{\mathrm{sp}}$, and  $\|\chi_s-\chi_s^{(J)}\|_2 \leq 3^{-\frac12} 2^{-\frac{J}{2}}\sqrt{s}$ imply\footnote{Note that $\|\chi_s-\chi_s^{(J)}\|_2$ is erroneously claimed to be $O(2^{-J})$ in \cite{Cellarosi-Marklof} but this does not affect 
the proof of Theorem 3.13 therein. In the same way, if we replaced  $O(2^{-\frac{J}{2}})$ in \eqref{lemma-final-approximations-irr-statement}  by $O(2^{-J})$, we would  get the same optimal power saving in Theorem \ref{thm-with-eta-dependent-exponent-irr}, see $\beta^*$ in \eqref{proof-main-thm-3-irr-3}.
}
\begin{align}
&\left|D_{\text{irr}}(\chi, \chi_b) - D_{\text{irr}}(\chi^{(J)}, \chi_b^{(J)})\right|\leq 8 \sqrt{b} \|\chi_b \|_{\mathrm{sp}}^5\, 2^{-\frac{J}{2}} \label{proof-lemma-final-approximations-irr-4b}.
\end{align}
Combining \eqref{proof-lemma-final-approximations-irr-1} with \eqref{proof-lemma-final-approximations-irr-2}, \eqref{proof-lemma-final-approximations-irr-3} and \eqref{proof-lemma-final-approximations-irr-4b} gives the result.
\end{proof}

\section{Proof of the main theorems}\label{section-final-proofs}
In this section, we combine the tail bounds in Lemmata \ref{L2.2} and \ref{L2.2_irrational} with the tail approximations from Section 
\ref{section_tail_approximations} to prove Theorems \ref{thm-with-eta-dependent-exponent-rat} and \ref{thm-with-eta-dependent-exponent-irr}, which are analogue to Theorem \ref{main-thm-4} but hold for $(f_1,f_2)=(\chi,\chi_b)$ and have explicit, $\eta$-dependent power savings, namely $O_\eta(R^{-\frac{2\eta}{6\eta(\eta-1)+1}})$ and $O_\eta(R^{-\frac{6\eta}{16\eta(\eta-1)+3}})$, respectively. 
We then take $\eta\to1^+$ to obtain Theorems \ref{thm-rat-epsilon-explicit} and \ref{thm-irr-epsilon-explicit}, which are explicit versions of Theorems \ref{main-thm-3-rat} and \ref{main-thm-3-irr}, respectively. 
\subsection{Rational case}
\import{./}{a-computation-03.tex}

\begin{theorem}\label{thm-with-eta-dependent-exponent-rat}
Let $\lambda$ be a probability measure on $\R$ which is absolutely continuous with respect to the Lebesgue measure. Let $1<\eta\leq2$ and $b\geq1$. Set $R_0^{\mathrm{rat}}=R_0^{\mathrm{rat}}(b,\eta)$ as in \eqref{def-R_0(b,eta)}. Then for all $R>R_0^{\mathrm{rat}}$ we have
\begin{align}
\lim_{N\to\infty }\lambda\!\left\{x\in\R: \left|\tfrac{1}{N}S_{N}(x;0,0)\overline{S_{\lfloor b N\rfloor}(x;0,0)}\right|>R^2\right\}=\frac{2D_{\mathrm{rat}}(\chi,\chi_b)}{\pi^2 R^4}\left(1\!+\!O_{\eta}\!\left(R^{-\frac{2\eta}{6\eta(\eta-1)+1}}\right)\!\right)\!,\label{statement-thm-with-eta-dependent-exponent-rat}
\end{align}
where the constant implied by the $O_{\eta}$-notation can be taken to be as in \eqref{def-P_eta}
\end{theorem}
\begin{proof}
Corollary \ref{key-tail-limit-theorem-rat} gives
\begin{align}
\lim_{N\to\infty }\lambda\left\{x\in\R:\: \left|\tfrac{1}{N}S_N(x;0,0)\overline{S_{\lfloor b N\rfloor}(x;0,0)}\right|>R^2\right\}=\mu_{\GamG}^{\bn}\left(|\Theta_{\chi}\overline{\Theta_{\chi_b}}|>R^2\right).\label{proof-main-thm-3-rat-0}
\end{align}
Let us now write $J=\alpha\log_2(R)$ and $\delta=R^{-\beta}$, for some $\alpha,\beta>0$. 
We then rewrite the assumptions of Lemma \ref{smooth_approx_rat} as 
\begin{align}
R^2>64  K_L(b)^2 C_\eta R^{2(\eta-1)\alpha},\hspace{.3cm}\sqrt{2}R^{-\frac{\alpha}{2}}<R^{-\beta}<\tha,\hspace{.3cm}, R^\alpha\geq 16,\label{rewrite-assumptions-lemma-smooth_approx_rat}
\end{align}
and those of Lemma \ref{lemma-final-approximations-rat}  as
\begin{align}
R^2> 16  K(b)^2 C_\eta R^{2(\eta-1)\alpha},\hspace{.3cm}0<R^{-\beta}<\tha,\hspace{.3cm}R^{\alpha}\geq 2. \label{rewrite-assumptions-lemma-final-approximations-rat}
\end{align}
In order for the first inequality of  \eqref{rewrite-assumptions-lemma-smooth_approx_rat} %
to hold for all sufficiently large $R$, we need  $(\eta-1)\alpha<1$. The same is true for the first inequality of \eqref{rewrite-assumptions-lemma-final-approximations-rat}. Similarly, the second inequality of \eqref{rewrite-assumptions-lemma-smooth_approx_rat} implies $\beta<\frac{\alpha}{2}$. 
By the lower bound in Remark \ref{remark-ratio-C3-C3L}, 
%
we have
$K_{L}(b)^2\leq \ha K(b)^2$.
Therefore if we assume that 
\begin{align}
R^2>32 K(b)^2 C_\eta R^{2(\eta-1)\alpha},\hspace{.3cm}R^\alpha>16,\hspace{.3cm} R^\beta>2,\hspace{.3cm} R^{\frac{\alpha}{2}-\beta}>\sqrt{2},\label{first-3-conditions-on-R}
\end{align} then \eqref{rewrite-assumptions-lemma-smooth_approx_rat} and \eqref{rewrite-assumptions-lemma-final-approximations-rat} hold and we can apply both  Lemma \ref{smooth_approx_rat} and Lemma \ref{lemma-final-approximations-rat}. We obtain
%
\begin{align}
&\mu_{\GamG}^{\bn}\left(|\Theta_{\chi}\overline{\Theta_{\chi_b}}|>R^2\right)=\\
&=\frac{2D_{\mathrm{rat}}(\chi,\chi_b)}{\pi^2R^4}\left(1+O(R^{-\beta})\right)\left(1+O_\eta(R^{-2\eta+2\eta(\eta-1)\alpha})\right)\left(1+O(R^{-\alpha})\right)+O\!\left(R^{-\alpha+2\beta-4}\right),\label{proof-main-thm-3-rat-1}
\end{align} 
where $D_{\mathrm{rat}}(\chi,\chi_b)$ is given explicitly in \eqref{computation-I(a,b)-statement}.
Now we choose $R$ sufficiently large so  that, in addition to satisfying \eqref{first-3-conditions-on-R},  each of the first three $O$-terms in \eqref{proof-main-thm-3-rat-1} is bounded above in absolute value by $1$. This means, recalling the implicit constants in 
 Lemma \ref{lemma-final-approximations-rat}, 
\begin{align}
R>&
\max\left\{(32 K(b)^2 C_\eta)^{\frac{1}{2-2(\eta-1)\alpha}},\: 16^{\frac{1}{\alpha}},\:2^{\frac{1}{\beta}},\:2^{\frac{1}{\alpha-2\beta}},\:30^{\frac{1}{\beta}},\: \right.\label{rat-first-part-of-R>max}\\
&\left.~~~~~~~~(32 K(b)^2 C_\eta )^{\frac{1}{\eta-\eta(\eta-1)\alpha}},\: (32b\|\chi_b\|_{\mathrm{sp}}^3)^\frac{1}{\alpha}\right\}.\label{rat-second-part-of-R>max}
\end{align}
Observing that $\frac{1}{2-2(\eta-1)\alpha}\leq \frac{1}{\eta-\eta(\eta-1)\alpha}$ for $1<\eta\leq2$ and $0<\alpha<\frac{1}{\eta-1}$ and that, by \eqref{sp-norm-of-chi-s}, $32 b\|\chi_b\|_{\mathrm{sp}}^3\geq32\cdot 9^2\geq 16$, we can rewrite the above condition on $R$ as
\begin{align}
R>R_0^{\mathrm{rat}}(\alpha,\beta;b,\eta):=\max\!\left\{(32 K(b)^2 C_\eta )^{\frac{1}{\eta-\eta(\eta-1)\alpha}},\:30^{\frac{1}{\beta}},\: 2^{\frac{1}{\alpha-2\beta}},\:(32b\|\chi_b\|_{\mathrm{sp}}^3)^\frac{1}{\alpha}\right\}.\label{proof-main-thm-3-rat-condition-on-R-2}
\end{align}
 %
Under the assumption \eqref{proof-main-thm-3-rat-condition-on-R-2}, combining the contribution of the error terms in \eqref{proof-main-thm-3-rat-1}, we can write\footnote{The transition from \eqref{proof-main-thm-3-rat-1} to \eqref{proof-main-thm-3-rat-2} uses the fact that $A(1+B)(1+C)(1+D)+E\leq A(1+4 B+C+2 D+\frac{E}{A})$ provided $A,E>0$ and $0<B,C,D<1$.}
\begin{align}
\mu_{\GamG}^{\bn}\left(|\Theta_{\chi}\overline{\Theta_{\chi_b}}|>R^2\right)&=\frac{2D_{\mathrm{rat}}(\chi,\chi_b)}{\pi^2R^4}\left(1+O_\eta\!\left(R^{-\beta}+R^{-2\eta+2\eta(\eta-1)\alpha}+R^{-\alpha}+R^{-\alpha+2\beta}\right)\right),
\label{proof-main-thm-3-rat-2}
\end{align}
where the implied constant in the $O_\eta$-notation in \eqref{proof-main-thm-3-rat-2} can be taken as
\begin{align}
\max\left\{120,\, 2^{10}K(b)^4 C_\eta^2,\, 2^6b\|\chi_b\|_{\mathrm{sp}}^3,\, \frac{2^{17}\|\chi_{b,L}\|_{\mathrm{sp}}^4}{\tfrac{2}{\pi^2}D_{\mathrm{rat}}(\chi,\chi_b)}\right\}.
\label{def-B_eta-0}
\end{align}
Using the trivial bounds $K(b)\geq 2^5$ and $C_\eta\geq2^6$, we see that $2^{10}K(b)^4C_\eta^2\geq 2^{42}\geq 120$ and hence the first term in \ref{def-B_eta-0} can be dropped.  Using the explicit formul\ae\, \eqref{piecewise-formula-C3-new}  and \eqref{sp-norm-of-chi-s} for $K(b)$ and $\|\chi_b\|_{\mathrm{sp}}$ respectively, we can verify that for $b\geq 1$ the bound  $K(b)^4\geq 2^{14} b\|\chi_b\|_{\mathrm{sp}}^3$ holds. Therefore $2^{10}K(b)^4C_\eta^2\geq 2^{22} K(b)^4\geq 2^{36}b\|\chi_b\|_{\mathrm{sp}}^3$ and hence the third term in \ref{def-B_eta-0} can be dropped. Finally, using the explicit formula \eqref{sp-norm-of-chi-s,L} for $\|\chi_{b,L}\|_{\mathrm{sp}}$, we can check that  the bound $K(b)^4\geq 2^{10}\|\chi_{b,L}\|_{\mathrm{sp}}^4$ holds for every $b\geq1$. 
This fact, along with the lower bound $D_{\mathrm{rat}}(\chi,\chi_b)\geq 2\log(2)$ (see Theorem \ref{thm-computation-constant-D_rat(chi,chi_b)}), yield 
  $2^{10}K(b)^4C_\eta^2\geq 2^{22}K(b)^4\geq 2^{32}\|\chi_{b,L}\|_{\mathrm{sp}}^4\geq\frac{\pi^2 2^{16}}{2\log(2)}\|\chi_{b,L}\|_{\mathrm{sp}}^4\geq\pi^2 2^{16}\frac{\|\chi_{b,L}\|_{\mathrm{sp}}^4}{D_{\mathrm{rat}(\chi,\chi_b)}}$ and hence the last term in \ref{def-B_eta-0} can also be dropped. 
We have shown that the implied constant in the $O_\eta$-notation in \eqref{proof-main-thm-3-rat-2} can be taken as
\begin{align}
P_0^{\mathrm{rat}}(b,\eta)=2^{10}K(b)^4 C_\eta^2. \label{def-P_eta-0}
\end{align}

The best power saving in \eqref{proof-main-thm-3-rat-2} comes from choosing $(\alpha^*,\beta^*)$ in the region 
$\{(\alpha,\beta):\:0<\alpha<\frac{1}{\eta-1},\: 0<\beta<\frac{\alpha}{2}\}$
 in order to maximise $m(\alpha,\beta)=\min\{\alpha,\beta,2\eta-2\eta(\eta-1)\alpha,\alpha-2\beta\}=\min\{\beta,2\eta-2\eta(\eta-1)\alpha,\alpha-2\beta\}$. That is
\begin{align}
(\alpha^*, \beta^*)=\left(\frac{6\eta}{6\eta(\eta-1)+1},\frac{2\eta}{6\eta(\eta-1)+1}\right),\label{proof-main-thm-3-rat-3}
\end{align}
at which $m(\alpha^*,\beta^*)=\beta^*$.
Therefore, combining \eqref{proof-main-thm-3-rat-0}, \eqref{proof-main-thm-3-rat-1}, \eqref{proof-main-thm-3-rat-condition-on-R-2}, \eqref{proof-main-thm-3-rat-2}, and \eqref{proof-main-thm-3-rat-3}, we obtain that for every $R>R_0^{\mathrm{rat}}(\alpha^*,\beta^*;b,\eta)$ we have
\begin{align}
\lim_{N\to\infty}\!\lambda\!\left\{x\in\R: \left|\tfrac{1}{N}S_N(x;0,0)\overline{S_{\lfloor b N\rfloor}(x;0,0)}\right|\!>R^2\right\}
=\frac{2D_{\mathrm{rat}}(\chi,\chi_b)}{\pi^2R^4}\left(1+O_\eta\!\left(R^{-\frac{2\eta}{6\eta(\eta-1)+1}}\right)\!\right)\!,\label{proof-main-thm-3-rat-4}
\end{align}
where the implicit constant in the $O_\eta$-notation in \eqref{proof-main-thm-3-rat-4} can be taken as $4P_0^{\mathrm{rat}}(b,\eta)$.
Note that 
\begin{align}
\frac{1}{\eta-\eta(\eta-1)\alpha^*}&=\frac{1}{\eta}+6(\eta-1),\label{proof-main-thm-3-rat-4a}\\
\frac{1}{\beta^*}=\frac{1}{\alpha^*-2\beta^*}
&=\frac{1}{2\eta}+3(\eta-1)
,\label{proof-main-thm-3-rat-4b}\\
\frac{1}{\alpha^*}&=\frac{1}{6\eta}+(\eta-1)
.\label{proof-main-thm-3-rat-4c}
\end{align}
We claim that 
\begin{align}
R_0^{\mathrm{rat}}(\alpha^*,\beta^*,\eta)=(32 K(b)^2 C_\eta)^{\frac{1}{\eta}+6(\eta-1)}.\label{claim-on-R'_0(alphastar,betastar,eta)} 
\end{align} 
From \eqref{proof-main-thm-3-rat-4b} we immediately have  $30^{\frac{1}{\beta^*}}\geq 2^{\frac{1}{\alpha^*-2\beta^*}}$. Then we see by \eqref{proof-main-thm-3-rat-4a}-\eqref{proof-main-thm-3-rat-4b} that 
\begin{align}
(32 K(b)^2 C_\eta )^{\frac{1}{\eta-\eta(\eta-1)\alpha^*}}=(2^{10}K(b)^4C_\eta^2)^{\frac{1}{2\eta}+3(\eta-1)}\geq 30^{\frac{1}{2\eta}+3(\eta-1)}=30^{\frac{1}{\beta^*}}
\end{align}
because $2^{10}K(b)^4C_\eta^2\geq 2^{10}K(1)^4 2^{12}\geq 2^{45}\geq30$. Finally,  \eqref{proof-main-thm-3-rat-4c} and the explicit bound $K(b)^2\geq 2^{10}b^{\frac{1}{6}}\|\chi_b\|_{\mathrm{sp}}^{\ha}$ imply
\begin{align}
(32 K(b)^2 C_\eta )^{\frac{1}{\eta-\eta(\eta-1)\alpha^*}}=(32K(b)^2C_\eta)^{\frac{1}{\eta}+6(\eta-1)}\geq (2^{\frac{5}{6}}b^{\frac{1}{6}}\|\chi_b\|_{\mathrm{sp}}^{\ha})^{\frac{1}{\eta}+6(\eta-1)}=(32b\|\chi_b\|_{\mathrm{sp}}^3)^{\frac{1}{\alpha^*}}
\end{align}
since
$2^5 K(b)^2 C_\eta\geq 2^{11}K(b)^2\geq 2^{21}b^{\frac{1}{6}}\|\chi_b\|_{\mathrm{sp}}^{\ha}\geq 2^{\frac{5}{6}}b^{\frac{1}{6}}\|\chi_b\|_{\mathrm{sp}}^{\ha}$. This concludes the proof of the claim \eqref{claim-on-R'_0(alphastar,betastar,eta)}.
 Now, 
  if we set 
\begin{align}
R_0^{\mathrm{rat}}(b,\eta)&:=(32 K(b)^2C_\eta)^{\frac{1}{\eta}+6(\eta-1)}
,\label{def-R_0(b,eta)}
\end{align}
then 
 for every $R>R_0^{\mathrm{rat}}(b,\eta)$ we have \eqref{proof-main-thm-3-rat-4} with an implicit constant of 
\begin{align}
P^{\mathrm{rat}}(b,\eta):=2^{12} K(b)^4 C_\eta^2
\label{def-P_eta}
\end{align}
and the theorem is proven.
\end{proof}

\begin{theorem}[Explicit version of Theorem \ref{main-thm-3-rat}]\label{thm-rat-epsilon-explicit}
Let $\lambda$ be a probability measure on $\R$ which is absolutely continuous with respect to the Lebesgue measure. Let $b\geq1$ and  $0<\varepsilon\leq1$, and set
\begin{align}
R^{\mathrm{rat}}(b,\varepsilon)=2^{38}K(b)^4\varepsilon^{-2}
\label{Rrat(b,eps)}
\end{align}
where $K(b)$ is given in \eqref{piecewise-formula-C3-new}.
Then for all $R>R^{\mathrm{rat}}(b,\varepsilon)$ we have
\begin{align}
\lim_{N\to\infty}\lambda\!\left\{x\in\R:\: \left|\tfrac{1}{N}S_N(x;0,0)\overline{S_{\lfloor b N\rfloor}(x;0,
0)}\right|>R^2\right\}=
\frac{2D_{\mathrm{rat}}(b)}{\pi^2 R^4}\left(1+O_\varepsilon\!\left(R^{-2+\varepsilon}\right)\right),
\label{statement-thm-rat-epsilon-explicit}
\end{align}
where  $D_{\mathrm{rat}}(b)=D_{\mathrm{rat}}(\chi,\chi_b)$ is given in \eqref{computation-I(a,b)-statement}, and the constant implied by the $O_{\varepsilon}$-notation in  \eqref{statement-thm-rat-epsilon-explicit} can be taken to be 
\begin{align}
2^{42}K(b)^4\varepsilon^{-4}.\label{implied-const-explicit-rat}
\end{align}
\end{theorem}
\begin{proof}
Note that the exponent in \eqref{statement-thm-with-eta-dependent-exponent-rat} satisfies $$\displaystyle\lim_{\eta\to1^+}\frac{2\eta}{6\eta(\eta-1)+1}=2^-.$$ This means that the closer $\eta>1$ is to $1$, the better the power saving in Theorem \ref{thm-with-eta-dependent-exponent-rat} is.
%
%
%
If we set $\frac{2\eta}{6\eta(\eta-1)+1}=2-\varepsilon$ for some $\varepsilon>0$, then $\varepsilon\mapsto \eta(\varepsilon)=\frac{7-3\varepsilon+\sqrt{3\varepsilon^2-18\varepsilon+25}}{6(2-\varepsilon)}$ maps $(0,1]$ to $(1,\tfrac{4+\sqrt{10}}{6}]$ monotonically. Considering the Taylor expansion of $\eta(\varepsilon)$ near $\varepsilon=0$ we find the lower bound $\eta(\varepsilon)>1+\frac{\varepsilon}{10}$.

We aim to find a simple upper bound for $R_0^{\mathrm{rat}}(b,\eta(\varepsilon))$. If the function  $\eta\mapsto R_0^{\mathrm{rat}}(b,\eta)$ were decreasing, then we could get an upper bound by replacing $R_0^{\mathrm{rat}}(b,\eta(\varepsilon))$ with $R_0^{\mathrm{rat}}(b,1+\frac{\varepsilon}{10})$. However, we only have monotonicity near $\eta=1$ and not on the whole interval $(1,\tfrac{4+\sqrt{10}}{6}]$. To address this issue, let us find an upper bound for $R_0^{\mathrm{rat}}(b,\eta)$ which is decreasing as a function of $\eta$ in the desired range. For the rest of this proof we will assume that $1<\eta\leq\tfrac{4+\sqrt{10}}{6}$. We have
%
\begin{align}
\zeta(\eta)\leq \frac{\sqrt{2}}{\eta-1}\label{upper-bound-zeta(eta)-with-sqrt2-constant}
\end{align}
(see Corollary \ref{cor-upper-and-lower-bounds-for-zeta5/4} in Appendix \ref{appendix-Zeta}). Recall that $C_\eta=2^{6\eta}\zeta(\eta)^2$. We use \eqref{def-R_0(b,eta)}, \eqref{upper-bound-zeta(eta)-with-sqrt2-constant}, the fact that $\frac{1}{\eta}+6(\eta-1)\leq2$, and the bound $(\eta-1)^{-2(\frac{1}{\eta}+6(\eta-1))}\leq2^5(\eta-1)^{-2}$ to obtain
\begin{align}
R_0^{\mathrm{rat}}(b,\eta)\leq \left(\frac{2^{6+6\eta}K(b)^2}{(\eta-1)^2}\right)^{\frac{1}{\eta}+6(\eta-1)}\leq \frac{2^{20+2\sqrt{10}}K(b)^4}{(\eta-1)^{2(\frac{1}{\eta}+6(\eta-1))}}\leq\frac{2^{25+2\sqrt{10}}K(b)^4}{(\eta-1)^2}.\label{monotonic-upperbound-of-R_0^rat(b,eta)}
\end{align}
Therefore, since the right-most upper bound in \eqref{monotonic-upperbound-of-R_0^rat(b,eta)} is decreasing in $\eta$, 
\begin{align}
R_0^{\mathrm{rat}}(b,\eta(\varepsilon))\leq\frac{2^{25+2\sqrt{10}}K(b)^{4}}{(\frac{\varepsilon}{10})^2}\leq 2^{38}K(b)^4\varepsilon^{-2}=:R^{\mathrm{rat}}(b,\varepsilon).
\end{align}
Therefore, if $R>R^{\mathrm{rat}}(b,\varepsilon)$, then by Theorem \ref{thm-with-eta-dependent-exponent-rat} we obtain
\begin{align}
\lim_{N\to\infty }\lambda \left\{x\in\R: \left|\tfrac{1}{N}S_{N}(x;0,0)\overline{S_{\lfloor b N\rfloor}(x;0,0)}\right|>R^2\right\}=\frac{2D_{\mathrm{rat}}(\chi,\chi_b)}{\pi^2 R^4}\left(1+O_{\varepsilon}\!\left(R^{-2+\varepsilon}\right)\right)\!.\label{statement-thm-with-epsilon-dependent-exponent-rat}
\end{align}
The constant implied by the $O_\varepsilon$-notation in \eqref{statement-thm-with-epsilon-dependent-exponent-rat} is given by $P^{\mathrm{rat}}(b,\eta(\varepsilon))$ using \eqref{def-P_eta} and can be bounded above as follows:
\begin{align}
P^{\mathrm{rat}}(b,\eta(\varepsilon))\leq2^{22+2\sqrt{10}}K(b)^4(\tfrac{\varepsilon}{10})^{-4}\leq 2^{42}K(b)^4\varepsilon^{-4}.
\end{align}
 \end{proof}

By taking $b=1$ in Theorem \ref{thm-rat-epsilon-explicit} and bounding $K(1)^4=59^4\leq 2^{24}$, we obtain the following
\begin{cor}[Explicit version of Theorem \ref{main-thm-1}]\label{cor-explicit-version-thm-1}
Let $\lambda$ be a probability measure on $\R$ which is absolutely continuous with respect to the Lebesgue measure. Let $0<\varepsilon\leq1$ and set $R^{\mathrm{rat}}(\varepsilon)=2^{62}\varepsilon^{-2}$. Then for all $R>R^{\mathrm{rat}}(\varepsilon)$ we have
\begin{align}
\lim_{N\to\infty}\lambda\!\left\{x\in\R:\: \left|\tfrac{1}{\sqrt N}S_N(x;0,0)\right|>R\right\}=
\frac{4\log 2}{\pi^2 R^4}\left(1+O_\varepsilon\!\left(R^{-2+\varepsilon}\right)\right),
\label{statement-thm-rat-epsilon-explicit-b=1}
\end{align}
the constant implied by the $O_{\varepsilon}$-notation in  \eqref{statement-thm-rat-epsilon-explicit-b=1} can be taken to be 
$2^{66}\varepsilon^{-4}$.
\end{cor}

\subsection{Irrational case}
\begin{theorem}\label{thm-with-eta-dependent-exponent-irr}
Let $\lambda$ be a probability measure on $\R$ which is absolutely continuous with respect to the Lebesgue measure. Let $(c,\alpha)\notin\Q^2$, $1<\eta\leq2$, and $b\geq1$. Set $R_0^{\mathrm{irr}}=R_0^{\mathrm{irr}}(b,\eta)$ as in \eqref{def-R_0(b,eta)-irr}. Then for all $R>R_0^{\mathrm{irr}}$ we have
\begin{align}
\lim_{N\to\infty }\lambda\!\left\{x\in\R\!: \left|\tfrac{1}{N}S_{N}(x;c,\alpha)\overline{S_{\lfloor b N\rfloor}(x;c,\alpha)}\right|>R^2\right\}\!=\!\frac{2D_{\mathrm{irr}}(\chi,\chi_b)}{\pi^2R^6}\left(1\!+\!O_{\eta}\!\left(R^{-\frac{6\eta}{16\eta(\eta-1)+3}}\right)\!\right)\!,\label{statement-thm-with-eta-dependent-exponent-irr}
\end{align}
where
\begin{align}
D_{\mathrm{irr}}(\chi, \chi_b) = \int_{-\infty}^{\infty}\int_0^{\pi}\left|\chi_\phi(w)\,(\chi_b)_\phi(w)\right|^3 \: \de\phi\,\de w\label{D_irr-def}
\end{align}
and the constant implied by the $O_{\eta}$-notation can be taken to be as in \eqref{def-P_eta-irr}
\end{theorem}
\begin{proof}
Corollary \ref{key-tail-limit-theorem-irr} gives
\begin{align}
\lim_{N\to\infty }\lambda\left\{x\in\R:\: \left|\tfrac{1}{N}S_N(x;c,\alpha)\overline{S_{\lfloor b N\rfloor}(x;c,\alpha)}\right|>R^2\right\}=\mu_{\GamG}^{}\left(|\Theta_{\chi}\overline{\Theta_{\chi_b}}|>R^2\right).\label{proof-main-thm-3-irr-0}
\end{align}
We argue as in the proof of Theorem \ref{thm-with-eta-dependent-exponent-rat} and write $J=\alpha\log_2(R)$, and $\delta=R^{-\beta}$, for some $\alpha,\beta>0$. 
%
%
The assumptions of Lemmata \ref{smooth_approx_irr} and  \ref{lemma-final-approximations-irr} read as \eqref{rewrite-assumptions-lemma-smooth_approx_rat} and  \eqref{rewrite-assumptions-lemma-final-approximations-rat} respectively, i.e. they are the same as in the rational case. Therefore, if we assume \eqref{first-3-conditions-on-R}, then we can apply both Lemma \ref{smooth_approx_irr} and Lemma \ref{lemma-final-approximations-irr}. 
\begin{align}
&\mu_{\GamG}^{}\left(|\Theta_{\chi}\overline{\Theta_{\chi_b}}|>R^2\right)=\\
&=\frac{2D_{\mathrm{irr}}(\chi,\chi_b)}{\pi^2R^6}\left(1+O(R^{-\beta})\right)\left(1+O_\eta(R^{-2\eta+2\eta(\eta-1)\alpha})\right)\left(1+O(R^{-\frac{\alpha}{2}})\right)+O\!\left(R^{-\frac{3}{2}\alpha+3\beta-6}\right).\label{proof-main-thm-3-irr-1}
\end{align} 
Now we choose $R$ sufficiently large so  that, in addition to satisfying \eqref{first-3-conditions-on-R}
,  each of the first three $O$-terms in \eqref{proof-main-thm-3-irr-1} is bounded above in absolute value by $1$. This means, recalling the implicit constants in 
 Lemma \ref{lemma-final-approximations-irr}, 
\begin{align}
R>&
\max\left\{(32 K(b)^2 C_\eta)^{\frac{1}{2-2(\eta-1)\alpha}},\: 16^{\frac{1}{\alpha}},\:2^{\frac{1}{\beta}},\:2^{\frac{1}{\alpha-2\beta}},\:126^{\frac{1}{\beta}},\: \right.\\
&\left.~~~~~~~~~~~(2^{\frac{11}{2}}K(b)^2C_\eta)^{\frac{1}{\eta-\eta(\eta-1)\alpha}}, (2^6 b\|\chi_b\|_{\mathrm{sp}}^{10})^{\frac{1}{\alpha}} \right\}\\
&=\max\left\{(2^{\frac{11}{2}}K(b)^2C_\eta)^{\frac{1}{\eta-\eta(\eta-1)\alpha}}, 126^{\frac{1}{\beta}}, 2^{\frac{1}{\alpha-2\beta}}, (2^6 b\|\chi_b\|_{\mathrm{sp}}^{10})^{\frac{1}{\alpha}} \right\}=:R_0^{\mathrm{irr}}(\alpha,\beta;b,\eta).
\label{proof-main-thm-3-irr-condition-on-R-2}
\end{align}
Under the assumption $R>R_0^{\mathrm{irr}}(\alpha,\beta;b,\eta)$,  
combining the contribution of the error terms in \eqref{proof-main-thm-3-irr-1}, we can write 
\begin{align}
\mu_{\GamG}^{}\left(|\Theta_{\chi}\overline{\Theta_{\chi_b}}|>R^2\right)&=\frac{2D_{\mathrm{irr}}(\chi,\chi_b)}{\pi^2 R^6}\left(1+O_\eta\!\left(R^{-\beta}+R^{-2\eta+2\eta(\eta-1)\alpha}+R^{-\frac{\alpha}{2}}+R^{-\frac{3}{2}\alpha+3\beta}\right)\right),
\label{proof-main-thm-3-irr-2}
\end{align}
where the implied constant in the $O_\eta$-notation in \eqref{proof-main-thm-3-rat-2} can be taken as
\begin{align}
\max\left\{504,\: 2^{11}K(b)^4 C_\eta^2,\: 16\sqrt{b}\|\chi_b\|_{\mathrm{sp}}^5,\: \frac{2^{22} \|\chi_{b,L}\|_{\mathrm{sp}}^5}{\frac{2}{\pi^2}D_{\mathrm{irr}}(\chi,\chi_b)}\right\}.
\label{def-B_eta-0-irr}
\end{align}
Since $2^{11}K(b)^4C_\eta^2\geq2^{46}\geq 504$, the first term in 
\eqref{def-B_eta-0-irr} can be dropped. Since $K(b)^4\geq 2^7 \sqrt{b}\|\chi_b\|_{\mathrm{sp}}^5$ for every $b\geq1$, we have $2^{11}K(b)^4C_\eta^2\geq2^{23}K(b)^4\geq 2^{30}\sqrt{b}\|\chi_b\|_{\mathrm{sp}}^5$ and the third term in \eqref{def-B_eta-0-irr} can be dropped. 
We conjecture that $D_{\mathrm{irr}}(\chi,\chi_b)\geq D_{\mathrm{irr}}(\chi,\chi)=3$. If that were the case, then we could use the bound $K(b)^4\geq 2^7 \|\chi_{b,L}\|_{\mathrm{sp}}^5$ to obtain $2^{11}K(b)^4C_\eta^2\geq2^{23}K(b)^4\geq 2^{30}\|\chi_{b,L}\|_{\mathrm{sp}}^5\geq\frac{\pi^2 2^{21}}{3}\|\chi_{b,L}\|_{\mathrm{sp}}^5\geq\pi^2 2^{21}\frac{\|\chi_{b,L}\|_{\mathrm{sp}}^5}{D_{\mathrm{irr}}(\chi,\chi_b)}$ and we would also drop the last term of \eqref{def-B_eta-0-irr}. However, since we do not have an explicit lower bound for $D_{\mathrm{irr}}(\chi,\chi_b)$, we avoid the comparison between the second and the last term of \eqref{def-B_eta-0-irr}. Therefore, we take the implied constant in the $O_\eta$-notation in \eqref{proof-main-thm-3-irr-2} 
as
%
\begin{align}
P_0^{\mathrm{irr}}(b,\eta)=\max\!\left\{2^{11}K(b)^4 C_\eta^2,\:\frac{\pi^2 2^{21}\|\chi_{b,L}\|_{\mathrm{sp}}^5}{D_{\mathrm{irr}}(\chi,\chi_b)}\right\}. \label{def-P_eta-irr-0}
\end{align}
Note that, refining the argument above, the second term of \eqref{def-P_eta-irr-0} can be dropped as long as we have the lower bound $D_{\mathrm{irr}}(\chi,\chi_b)\geq2^{-9}\pi^2$.
The best power saving in \eqref{proof-main-thm-3-irr-2} comes from choosing $(\alpha^*,\beta^*)$ in the region 
$\{(\alpha,\beta):\:0<\alpha<\frac{1}{\eta-1},\: 0<\beta<\frac{\alpha}{2}\}$
 in order to maximise $m(\alpha,\beta)=\min\{\frac{\alpha}{2},\beta,2\eta-2\eta(\eta-1)\alpha,\frac{3}{2}\alpha-3\beta\}=\min\{\beta,2\eta-2\eta(\eta-1)\alpha,\frac{3}{2}\alpha-3\beta\}$. That is
 \begin{align}
(\alpha^*, \beta^*)=\left(\frac{16\eta}{16\eta(\eta-1)+3},\frac{6\eta}{16\eta(\eta-1)+3}\right)\label{proof-main-thm-3-irr-3}
\end{align}
at which $m(\alpha^*,\beta^*)=\beta^*$. Therefore, combining \eqref{proof-main-thm-3-irr-0}, \eqref{proof-main-thm-3-irr-1}, \eqref{proof-main-thm-3-irr-condition-on-R-2}, \eqref{proof-main-thm-3-irr-2}, and \eqref{proof-main-thm-3-irr-3}, we obtain that for every $R>R_0^{\mathrm{irr}}(\alpha^*,\beta^*,\eta)$ we have
\begin{align}
\lim_{N\to\infty}\!\lambda\!\left\{x\in\R: \left|\tfrac{1}{N}S_N(x;c,\alpha)\overline{S_{\lfloor b N\rfloor}(x;c,\alpha)}\right|\!>R^2\right\}
=\frac{2D_{\mathrm{irr}}(\chi,\chi_b)}{\pi^2R^6}\left(1+O_\eta\!\left(R^{-\frac{6\eta}{16\eta(\eta-1)+3}}\right)\!\right)\!,\label{proof-main-thm-3-irr-4}
\end{align}
where the implicit constant in the $O_\eta$-notation in \eqref{proof-main-thm-3-irr-4} can be taken as $4P_0^{\mathrm{irr}}(b,\eta)$. 
Note that
\begin{align}
\frac{1}{\eta-\eta(\eta-1)\alpha^*}&=\frac{1}{\eta}+\frac{16}{3}(\eta-1),\label{proof-main-thm-3-irr-4a}\\
\frac{1}{\beta^*}&=\frac{1}{2\eta}+\frac{8}{3}(\eta-1)
,\label{proof-main-thm-3-irr-4b}\\
\frac{1}{\alpha^*-2\beta^*}&=\frac{3}{4\eta}+4(\eta-1)
,\label{proof-main-thm-3-irr-4c}\\
\frac{1}{\alpha^*}&=\frac{3}{16\eta}+(\eta-1)
.\label{proof-main-thm-3-irr-4d}
\end{align}
We claim that 
\begin{align}
R_0^{\mathrm{irr}}(\alpha^*,\beta^*;b,\eta)=(2^{\frac{11}{2}}K(b)^2C_\eta)^{\frac{1}{\eta}+\frac{16}{3}(\eta-1)}.\label{claim-on-R'_0(alphastar,betastar,eta)-irr} 
\end{align} 
To see this, first use \eqref{proof-main-thm-3-irr-4b}-\eqref{proof-main-thm-3-irr-4c} to see that
$126^{\frac{1}{\beta^*}}=(126^{\frac{2}{3}})^{\frac{3}{4\eta}+4(\eta-1)}\geq 2^{\frac{3}{4\eta}+4(\eta-1)}=2^{\frac{1}{\alpha^*-2\beta^*}}$. Then \eqref{proof-main-thm-3-irr-4a}-\eqref{proof-main-thm-3-irr-4b} to  obtain
\begin{align}
(2^{\frac{11}{2}}K(b)^2C_\eta)^{\frac{1}{\eta-\eta(\eta-1)\alpha^*}}=(2^{11}K(b)^4C_\eta^2)^{\frac{1}{2\eta}+\frac{8}{3}(\eta-1)}\geq126^{\frac{1}{2\eta}+\frac{8}{3}(\eta-1)}\geq 126^{\frac{1}{\beta^*}}
\end{align}
because $2^{11}K(b)^4C_\eta^2\geq 2^{46}\geq 126$. Finally, \eqref{proof-main-thm-3-irr-4d} and  the explicit bound  $K(b)^2\geq2^5 b^{\frac{3}{16}}\|\chi_b\|_{\mathrm{sp}}^{\frac{15}{8}}$ yield
\begin{align}
(2^{\frac{11}{2}}K(b)^2C_\eta)^{\frac{1}{\eta-\eta(\eta-1)\alpha^*}}=(2^{\frac{11}{2}}K(b)^2C_\eta)^{\frac{1}{\eta}+\frac{16}{3}(\eta-1)}\\\geq(2^{\frac{9}{8}}b^{\frac{3}{16}}\|\chi_b\|_{\mathrm{sp}}^{\frac{15}{8}})^{\frac{1}{\eta}+\frac{16}{3}(\eta-1)}=(2^6 b\|\chi_b\|_{\mathrm{sp}}^{10})^{\frac{1}{\alpha^*}} 
\end{align}
since $2^{\frac{11}{2}}K(b)^2C_\eta\geq2^{\frac{23}{2}}K(b)^2\geq2^{\frac{33}{2}} b^{\frac{3}{16}}\|\chi_b\|_{\mathrm{sp}}^{\frac{15}{8}}\geq 2^{\frac{9}{8}} b^{\frac{3}{16}}\|\chi_b\|_{\mathrm{sp}}^{\frac{15}{8}}$. This concludes the proof of the claim \eqref{claim-on-R'_0(alphastar,betastar,eta)-irr}. 
%
Therefore, if we set
\begin{align}
R_0^{\mathrm{irr}}(b,\eta)=(2^{\frac{11}{2}}K(b)^2C_\eta)^{\frac{1}{\eta}+\frac{16}{3}(\eta-1)},\label{def-R_0(b,eta)-irr}
\end{align}
then for every $R>R_0^{\mathrm{irr}}(b,\eta)$ we have \eqref{proof-main-thm-3-irr-4} with an implicit constant of 
\begin{align}
P^{\mathrm{irr}}(b,\eta):=\max\!\left\{2^{13}K(b)^4 C_\eta^2,\:\frac{\pi^2 2^{24}\|\chi_{b,L}\|_{\mathrm{sp}}^5}{D_{\mathrm{irr}}(\chi,\chi_b)}\right\}.\label{def-P_eta-irr}
\end{align}
and the theorem is proven.
\end{proof}

\begin{theorem}[Explicit version of Theorem \ref{main-thm-3-irr}]\label{thm-irr-epsilon-explicit} 
Let $\lambda$ be a probability measure on $\R$ which is absolutely continuous with respect to the Lebesgue measure.   Let $(c,\alpha)\notin \Q^2$,  $b\geq1$ and  $0<\varepsilon\leq1$. Set
\begin{align}
R^{\mathrm{irr}}(b,\varepsilon)=2^{39}K(b)^4\varepsilon^{-2}
\label{Rirr(b,eps)}
\end{align}
where $K(b)$ is given in \eqref{piecewise-formula-C3-new}.
Then for all $R>R_b^{\mathrm{irr}}(\varepsilon)$ we have
\begin{align}
\lim_{N\to\infty}\lambda\!\left\{x\in\R:\: \left|\tfrac{1}{N}S_N(x;c,\alpha)\overline{S_{\lfloor b N\rfloor}(x;c,
\alpha)}\right|>R^2\right\}=
\frac{2D_{\mathrm{irr}}(b)}{\pi^2 R^6}\left(1+O_\varepsilon\!\left(R^{-2+\varepsilon}\right)\right),\label{statement-thm-irr-epsilon-explicit}
\end{align}
where  $D_{\mathrm{irr}}(b)=D_{\mathrm{irr}}(\chi,\chi_b)$ is given in \eqref{D_irr-def}, 
and the constant implied by the $O_{\varepsilon}$-notation in  \eqref{statement-thm-irr-epsilon-explicit} can be taken to be
\begin{align}
\max\left\{2^{41}K(b)^4\varepsilon^{-4},\:\ 2^{28}\max\{b^5,3^{10}\}D_{\mathrm{irr}}(b)^{-1}\right\}.
\label{implied-constant-irr-epsilon-version}
\end{align}
\end{theorem}
\begin{proof}
Note that the exponent in \eqref{statement-thm-with-eta-dependent-exponent-irr} satisfies $$\displaystyle\lim_{\eta\to1^+}\frac{6\eta}{16\eta(\eta-1)+3}=2^-.$$ This means that the closer $\eta>1$ is to $1$, the better the power saving in Theorem \ref{thm-with-eta-dependent-exponent-irr} is.
If we set $\frac{6\eta}{16\eta(\eta-1)+3}=2-\varepsilon$, then $\varepsilon\mapsto\eta(\varepsilon)=\frac{19-8\varepsilon+\sqrt{16\varepsilon^2-112\varepsilon+169}}{16(2-\varepsilon)}$ maps $(0,1]$ to $(1,\frac{11+\sqrt{73}}{16}]$ monotonically. Hence we assume that $1<\eta\leq\frac{11+\sqrt{73}}{16}$. We aim to employ the lower bound $\eta(\varepsilon)>1+\frac{3}{26}\varepsilon$. Note that, by Corollary \ref{cor-upper-and-lower-bounds-for-zeta5/4}, \eqref{upper-bound-zeta(eta)-with-sqrt2-constant} holds in this case. We use \eqref{def-R_0(b,eta)-irr}, along with the fact that $\frac{1}{\eta}+\frac{16}{3}(\eta-1)\leq2$ and the bound $(\eta-1)^{-2(\frac{1}{\eta}+\frac{16}{3}(\eta-1))}\leq 2^5(s-1)^{-2}$ to obtain
\begin{align}
R_0^{\mathrm{irr}}(b,\eta)\leq \left(\frac{2^{\frac{13}{2}+6\eta}K(b)^2}{(\eta-1)^2}\right)^{\frac{1}{\eta}+\frac{16}{3}(\eta-1)}\leq\frac{2^{\frac{85}{4}+\frac{3}{4}\sqrt{73}}K(b)^4}{(\eta-1)^{2(\frac{1}{\eta}+\frac{16}{3}(\eta-1))}}\leq\frac{2^{\frac{105}{4}+\frac{3}{4}\sqrt{73}}K(b)^4}{(\eta-1)^2}.
\end{align}
By monotonicity we have
\begin{align}
R_0^{\mathrm{irr}}(b,\eta(\varepsilon))\leq\frac{2^{\frac{105}{4}+\frac{3}{4}\sqrt{73}}K(b)^4}{(\frac{3}{26}\varepsilon)^2}\leq2^{39}K(b)^4\varepsilon^{-2}=:R^{\mathrm{irr}}(b,\varepsilon).
\end{align}
Therefore, if $R>R^{\mathrm{irr}}(b,\varepsilon)$, then by Theorem \ref{thm-with-eta-dependent-exponent-irr} we obtain
\begin{align}
\lim_{N\to\infty }\lambda\!\left\{x\in\R\!: \left|\tfrac{1}{N}S_{N}(x;c,\alpha)\overline{S_{\lfloor b N\rfloor}(x;c,\alpha)}\right|>R^2\right\}\!=\!\frac{2D_{\mathrm{irr}}(\chi,\chi_b)}{\pi^2R^6}\left(1\!+\!O_{\eta}\!\left(R^{-2+\varepsilon}\right)\!\right)
\end{align}
where the constant implied by the $O_\varepsilon$-notation is given by $P^{\mathrm{irr}}(b,\eta(\varepsilon))$ using \eqref{def-P_eta-irr}. 
%
%
%
The bounds
$2^{11}K(b)^4C_{\eta(\varepsilon)}^2\leq 
2^{\frac{85}{4}+\frac{3}{4}\sqrt{73}}K(b)^4 (\eta(\varepsilon)-1)^{-4}\leq 2^{41}K(b)^4 \varepsilon^{-4}$ and \eqref{sp-norm-of-chi-s,L} imply that the constant implied by the $O_\varepsilon$-notation in \eqref{statement-thm-irr-epsilon-explicit}
can be taken as \eqref{implied-constant-irr-epsilon-version}
%
\end{proof}

If we take $b=1$ in Theorem \ref{thm-irr-epsilon-explicit} and, 
use the fact that $D_{\mathrm{irr}}(1)=3$ (see (3.121) in \cite{Cellarosi-Marklof}), we obtain the following
\begin{cor}[Explicit version of Theorem \ref{main-thm-2}]\label{cor-explicit-version-thm-2}
Let $\lambda$ be a probability measure on $\R$ which is absolutely continuous with respect to the Lebesgue measure. Let $(c,\alpha)\notin \Q^2$ and  $0<\varepsilon\leq1$. Set $R^{\mathrm{irr}}(\varepsilon)=2^{63}\varepsilon^{-2}
$. Then for all $R>R^{\mathrm{irr}}(\varepsilon)$ we have
\begin{align}
\lim_{N\to\infty}\lambda\!\left\{x\in\R:\: \left|\tfrac{1}{\sqrt N}S_N(x;c,\alpha)\right|>R\right\}=
\frac{6}{\pi^2 R^6}\left(1+O_\varepsilon\!\left(R^{-2+\varepsilon}\right)\right),
\label{statement-thm-irr-epsilon-explicit-b=1}
\end{align}
the constant implied by the $O_{\varepsilon}$-notation in  \eqref{statement-thm-irr-epsilon-explicit-b=1} can be taken to be  $2^{65}\varepsilon^{-4}$.
\end{cor}
\subsection{Final discussion on the explicit constants}\label{section-improve-constants}
As pointed out in the Introduction, we put a lot of effort into finding explicit, simple-looking constants for the implied $O_\varepsilon$-notation in Theorems \ref{main-thm-1}-\ref{main-thm-2} and in their generalisations \ref{main-thm-3-rat}-\ref{main-thm-3-irr}. These constants are given in Corollaries \ref{cor-explicit-version-thm-1} and \ref{cor-explicit-version-thm-2}, and in  
\eqref{Rrat(b,eps)}, \eqref{implied-const-explicit-rat}, \eqref{Rirr(b,eps)}, \eqref{implied-constant-irr-epsilon-version}.
We do not claim that the dependence of these constants upon $b$ and $\varepsilon$ is optimal. One can in principle make these constants smaller by improving \eqref{def-R_0(b,eta)}, \eqref{def-P_eta}, \eqref{def-R_0(b,eta)-irr}, and \eqref{def-P_eta-irr}. This task would  involve fairly elementary estimates, but would also be extremely lengthy and we purposefully avoided it. 
For the enthusiastic reader who seeks to optimise the constants, we list below some steps in which improvements can be made, and highlight some of the effects of such potential changes. 
\begin{enumerate}
\item To improve the constant $C_\eta=2^{6\eta}\zeta(\eta)^2$ in Lemma \ref{L2.1}: 
\begin{enumerate}
\item Bound \eqref{proof-L2.1-3} above by $\zeta(\eta)(2^\eta-1)+\zeta(\eta,\frac{3}{2})$, where $\zeta(\eta,a)$ is the Hurwitz zeta function $\zeta(\eta,a)=\sum_{n=0}^\infty(n+a)^{-\eta}$.
\item Replace $A_\eta$ by  $A'_\eta=(\zeta(\eta)(2^\eta-1)+\zeta(\eta,\frac{3}{2}))^2$ and $B_\eta$ by $B'_\eta=2(\zeta(\eta)(2^\eta-1)+\zeta(\eta,\frac{3}{2}))$ in \eqref{upperbound-with-A_eta-and-B_eta}. These new constants satisfy the inequality $B'_\eta\leq \frac{7}{20} A'_\eta$.
\item Replace $C_\eta$ by $A'_\eta(2^\frac{\eta}{2}+\frac{7}{20})$. Near $\eta=1$, the improved constant would be smaller than $C_\eta$ by at least a factor of $9$.  
\end{enumerate} 
\item With the revised definition of $C_\eta$ discussed above, the statement of Lemma  \ref{L2.2} would still holds, but its proof would need to be changed. In fact, since $C_\eta$ enters the definition of $\tilde\kappa$,  
one has to show that 
\eqref{proof-L2.2-9} still implies \eqref{proof-L2.2-6a0}. Albeit true, this would be more laborious to prove.  Similarly, Lemma \ref{L2.2_irrational} would still hold. Its proof would not require major changes as it would still be easy to see that \eqref{proof-L2.2-irr-6} implies \eqref{proof-L2.2-irr-1a}. 
\item The lower bound $C_\eta\geq2^6$, which we use in the proof of Theorems \ref{thm-with-eta-dependent-exponent-rat} and \ref{thm-with-eta-dependent-exponent-irr}, would still hold but its proof would require some effort.
\item The constants $K(s)$ and $K_L(s)$ in the key Proposition \ref{prop-key-estimate-kappa_eta(chi_s^J)} could be improved by avoiding the bound $2^{\frac{\eta}{2}}\leq2$ for $1<\eta\leq2$ in \eqref{where-C_3-appears-1} and  \eqref{where-C_3L-appears-1}. This would complicate the statements and the arguments in Sections \ref{section_tail_approximations}-\ref{section-final-proofs} by introducing dependence of such constants upon $\eta$. 
\item In the proof of Lemma \ref{lemma-final-approximations-rat}, in addition to improving $C_\eta$ as discussed in point 1 above, we could avoid the use of the inequality $K(1)\leq K(b)$ and replace  \eqref{inequality-after-using-K1<Kb-rat} with $R^2>16\cdot C_\eta K(1)K(b)2^{2(\eta-1)J}=16\cdot59\cdot C_\eta K(b)2^{2(\eta-1)J}$. The improvement would be especially beneficial for large values of $b$, since $K(b)\sim\frac{2}{3}b^3$ as $b\to\infty$. Moreover, restricting $1<\eta\leq \eta_0$ (see point 8a below), one could replace the factor $16$ with $2^{2+\eta_0}$. In the same way, we could lower the implied constant $2^{10}C_\eta^2 K(b)^4$ by replacing \eqref{proof-lemma-final-approximations-rat-2}
with $3\cdot 2^{\eta_0^2+2\eta_0}C_\eta^{\eta_0}K(b)^{2\eta_0} 2^{2\eta(\eta-1) J} R^{-2\eta}$. Similar considerations hold for Lemma \ref{lemma-final-approximations-irr}.
\item In Lemma \ref{smooth_approx_rat} the power 2 of $K_L(b)$ could not be reduced due to \eqref{proof-lemma-smooth_approx_rat-3}.  On the other hand, assuming $1<\eta\leq \eta_0$ (see point 8a below), the factor $64$ in the assumption of the lemma could be lowered to $2^{\frac{7}{2}+\eta_0}$. The same holds for Lemma \ref{smooth_approx_irr}.
\item If we carried out the improvements described in points 5-6 above, then we would need to modify the proof of Theorem \ref{thm-with-eta-dependent-exponent-rat}. Specifically, we would replace  the first inequality of \eqref{rewrite-assumptions-lemma-smooth_approx_rat} by $R^2>2^{\frac{7}{2}+\eta_0}K_L(b)^2C_\eta R^{2(\eta-1)\alpha}$ and the first inequality of \eqref{rewrite-assumptions-lemma-final-approximations-rat} by $R^2>59\cdot2^{2+\eta_0}K(b)C_\eta 2^{2(\eta-1)\alpha}$. Consequently, because of the different powers of $K_L(b)$ and $K(b)$, Remark \ref{remark-ratio-C3-C3L} could not be used and the sufficient condition \eqref{first-3-conditions-on-R} for \eqref{rewrite-assumptions-lemma-smooth_approx_rat}-\eqref{rewrite-assumptions-lemma-final-approximations-rat}
 would have to be changed. In turn, this change would affect \eqref{rat-first-part-of-R>max}, \eqref{proof-main-thm-3-rat-condition-on-R-2}, \eqref{claim-on-R'_0(alphastar,betastar,eta)}, and \eqref{def-R_0(b,eta)}. 
 The lowering of the implied constants $2^{10}C_\eta^2 K(b)^4$ in Lemma \ref{lemma-final-approximations-rat} would  affect \eqref{rat-second-part-of-R>max}, \eqref{proof-main-thm-3-rat-condition-on-R-2},  \eqref{claim-on-R'_0(alphastar,betastar,eta)}, \eqref{def-R_0(b,eta)}, but also \eqref{def-B_eta-0}, \eqref{def-P_eta-0}, and  \eqref{def-P_eta}. Ultimately the term $(K(b)^2C_\eta)^{\frac{1}{\eta}+6(\eta-1)}$ in $R_0^{\mathrm{rat}}(n,\eta)$ could be replaced by its square root. The term $K(b)^{4}C_\eta^2$  in \eqref{def-B_eta-0} could be replaced by $K(b)^\eta C_\eta^\eta$ but \eqref{def-P_eta-0} would have to be written as the maximum of two terms (one depending on $b$ and $\eta$ and one only dependent on $b$). 
 Similar considerations apply to Theorem \ref{thm-with-eta-dependent-exponent-irr}, ultimately affecting \eqref{def-R_0(b,eta)-irr} and \eqref{def-P_eta-irr}.
 
 \item Reducing the interval $0<\varepsilon\leq 1$ in Theorem \ref{thm-rat-epsilon-explicit} to $0<\varepsilon\leq\varepsilon_0$ would improve the constant \eqref{Rrat(b,eps)}  (and hence the constant $R^{\mathrm{rat}}(\varepsilon)$ in Corollary \ref{cor-explicit-version-thm-1}) 
 as follows:
\begin{enumerate}
\item It would reduce the interval $1<\eta\leq\frac{4+\sqrt{10}}{6}$ to $1<\eta\leq\eta_0=\frac{7-3\varepsilon_0+\sqrt{3\varepsilon_0^2-18\varepsilon_0+25}}{6(2-\varepsilon_0)}$ in the proof of the theorem. 
\item The constant $\sqrt{2}$ in Corollary \ref{cor-upper-and-lower-bounds-for-zeta5/4} (and hence in \eqref{upper-bound-zeta(eta)-with-sqrt2-constant}) could  be lowered and hence improve the first inequality in \eqref{monotonic-upperbound-of-R_0^rat(b,eta)}. 
\item For sufficiently small $\varepsilon_0$ we would gain the monotonicity of the function $\eta\mapsto R_0^{\mathrm{rat}}(b,\eta)$ (or its improved version, as described in point 7 above) and of  the second term in \eqref{monotonic-upperbound-of-R_0^rat(b,eta)} (i.e. the upper bound of \eqref{def-R_0(b,eta)}  obtained using \eqref{upper-bound-zeta(eta)-with-sqrt2-constant}) on the interval $(1,\eta_0]$. 
\item Even if we did not alter \eqref{def-R_0(b,eta)}, the bound $\frac{1}{\eta}+6(\eta-1)\leq2$, used in \eqref{monotonic-upperbound-of-R_0^rat(b,eta)}, could be replaced by $\frac{1}{\eta}+6(\eta-1)\leq \frac{2}{2-\varepsilon_0}$ and consequently reduce the exponent to which $K(b)^2$ is raised from 2 to slightly above 1. Therefore $K(b)^4$ in \eqref{Rrat(b,eps)} could be replaced by $K(b)^{\frac{4}{2-\varepsilon_0}}$. The dependence of \eqref{Rrat(b,eps)} upon $\varepsilon$ would be unchanged.  If we did alter \eqref{def-R_0(b,eta)} as described in points 5-7 above, then the dependence of \eqref{Rrat(b,eps)} upon $b$ and $\varepsilon$ would be $K(b)^{\frac{2}{2-\varepsilon_0}}\varepsilon^{-1}$ (up to a multiplicative constant). 
\end{enumerate}
\item 
The improvement of the constant \eqref{implied-const-explicit-rat} is a little more subtle. The reason is that, as outlined in point 7 above,  \eqref{def-P_eta-0} would be written as the maximum of two terms. The term involving $b$ and $\eta$ would include $K(b)$ raised to a power slightly above $1$ (improving the power $4$ in \eqref{implied-const-explicit-rat}) and $\varepsilon$ raised to a power slightly below $-2$ (improving the power $-4$). We could improve the constant $2^{66}\varepsilon^{-4}$ in Corollary \ref{cor-explicit-version-thm-1} accordingly.

\item 
Arguing as in points 8-9, we could also improve the constants \eqref{Rirr(b,eps)} and \eqref{implied-constant-irr-epsilon-version} in   
Theorem \ref{thm-irr-epsilon-explicit} and those in Corollary \ref{cor-explicit-version-thm-2}.
\end{enumerate}

%% file: a-computation-03.tex
It is convenient to express $D_{\mathrm{rat}}(\chi,\chi_b)$ from  Lemma \ref{lemma-final-approximations-rat} in closed form.   
Using the definition \eqref{def-D_rat(f_1,f_2)}  we have 
%
\begin{align}
D_{\mathrm{rat}}(\chi,\chi_b)
=\displaystyle\int_{0}^\pi |\chi_\phi(0)|^2\, |(\chi_b)_\phi(0)|^2\, d\phi.\label{computation-b-1}
\end{align}
\begin{theorem}\label{thm-computation-constant-D_rat(chi,chi_b)}
Let $b\geq1$.  Then
\begin{align}
D_{\mathrm{rat}}(\chi,\chi_b)=\begin{cases}
2\log2&\mbox{if $b=1$;}\\
2b\,\mathrm{coth}^{-1}(b)+\frac{1}{2}\log(b^2-1)+\frac{b^2}{2}\log(1-\frac{1}{b^2})&\mbox{if $b>1$.}
\end{cases}
\label{computation-I(a,b)-statement}
\end{align}
\end{theorem}
\begin{proof}
Note that for $0<\phi<\pi$ we have 
\begin{align} 
|f_\phi(0)|^2=\frac{1}{\sin(\phi)} \left| \,\int_{-\infty}^\infty \e{
v^2 \frac{\cot(\phi)}{2}
}
f(v)dv
\right|^2
\end{align}
and the change of variables $x=\frac{\cot(\phi)}{2}$ in \eqref{computation-b-1} yields
\begin{align}
D_{\mathrm{rat}}(\chi,\chi_b)=2\int_{-\infty}^\infty \left|\int_0^1 \e{v^2 x}dv\right|^2 \left|\int_0^b \e{w^2 x}dw\right|^2dx.\label{computation-b-2_0}
\end{align}
By performing the  three successive changes of variables $w'=w/b$, $x'=bx$, and $v'=v/b$, we obtain
\begin{align}
D_{\mathrm{rat}}(\chi,\chi_b)=b^2D_{\mathrm{rat}}(\chi,\chi_{1/b})
\label{computation-b-2}
\end{align}
%
Therefore it is enough to compute $D_{\mathrm{rat}}(\chi,\chi_c)$ with $0<c\leq1$. Observing that the integrand is even with respect to the variable $x$, we can rewrite it 
as
\begin{align}
D_{\mathrm{rat}}(\chi,\chi_c)=4\int_0^\infty p(x)\,q\!\left(\tfrac{1}{x}\right)\frac{dx}{x},\label{computation-b-4}
\end{align}
where
\begin{align}
p(x)=\displaystyle\int_0^1\int_0^1 \e{x(v_1^2-v_2^2)}dv_1dv_2 \hspace{.3cm}\mbox{and}\hspace{.3cm}q(x)=\displaystyle\frac{1}{x}\int_0^c\int_0^c\e{\frac{1}{x}(w_1^2-w_2^2)}dw_1dw_2.
\end{align}
The integral in \eqref{computation-b-4} can be seen as a Mellin convolution (evaluated at $z=1$, see also \eqref{Mellin-inversion-formula}) over the multiplicative group $\R_{>0}$. 
Let us compute the Mellin transforms of $p$ and $q$ respectively. Let $0<\Re(s)<1$. We use Fubini's theorem to change the order of integration and  integrate with respect to $x$ first and  then  with respect to the two remaining variables ($dv_1\,dv_2$ and $dw_1\,dw_2$ respectively). We obtain
\begin{align}
P(s)&=\displaystyle\int_0^\infty x^{s-1}p(x)dx=\int_0^1\int_0^1\int_0^\infty x^{s-1}\e{x(v_1^2-v_2^2)}dx\,dv_1\,dv_2=\\
&=\displaystyle\frac{\pi^{\tfrac{3}{2}-s}\csc\!\left(\frac{\pi s}{2}\right)}{2^{s+2}(1-s)\,\Gamma\!\left(\tfrac{3}{2}-s\right)}
\end{align}
and
\begin{align}
Q(s)&=\displaystyle\int_0^\infty x^{s-1}q(x)dx=\int_0^c\int_0^c\int_0^\infty x^{s-2}\e{x^{-1}(w_1^2-w_2^2)}dx\,dw_1\,dw_2=\\
&=\displaystyle\frac{c^{2s}\pi^{\tfrac{1}{2}+s}\,\sec\!\left(\frac{\pi s}{2}\right)}{2^{3-s}s\,\Gamma\!\left(\frac{1}{2}+s\right)}.
\end{align}
We use the double angle formula for sine and the identity $\Gamma\!\left(\tfrac{3}{2}-s\right)\Gamma\!\left(\tfrac{1}{2}+s\right)=\frac{1}{2}\pi(1-2s)\sec(\pi s)$ to simplify
\begin{align}
P(s)Q(s)&\displaystyle=
\frac{\pi c^{2s}}{8}\frac{\cot(\pi s)}{s(s-1)(2s-1)}.
\end{align}
Now we apply the Mellin inversion formula 
\begin{align}
\int_0^\infty p(x)\,q\!\left(\tfrac{z}{x}\right)\frac{dx}{x}=\int_{\frac{1}{2}-i\infty}^{\frac{1}{2}+i\infty} P(s)Q(s)\,z^{-s}ds\label{Mellin-inversion-formula}
\end{align}
in the particular case $z=1$. Note that for $s\to\frac{1}{2}$ we have $\cot(\pi s)=-\pi\left(s-\frac{1}{2}\right)+O((s-\frac{1}{2})^3)$ and therefore $s=\frac{1}{2}$ is not a pole of $P(s)Q(s)$. The poles of $\cot(\pi s)$ are the integers and $\mathrm{Res}(\cot(\pi s), s=n)=\frac{1}{\pi}$ for every $n\in \Z$.
Therefore 
\begin{align}
\mathrm{Res}\!\left(P(s)Q(s),s=n\right)=\begin{cases}\displaystyle\frac{-3c^2+2c^2\log(c)}{8} &\mbox{if $n=1$,}\\ \\ \displaystyle\frac{c^{2n}}{8 n(n-1)(2n-1)}=\frac{c^{2n}}{8}\left(\frac{1}{n}-\frac{2}{n-\frac{1}{2}}+\frac{1}{n-1}\right)&\mbox{if $n\geq2$.}\end{cases}
\end{align}
By the residue's theorem, we have
\begin{align}
D_{\mathrm{rat}}(\chi,\chi_c)=&\frac{4}{2\pi i}\int_{\frac{1}{2}-i\infty}^{\frac{1}{2}+i\infty}P(s)Q(s)ds=-4\sum_{n\geq1}\mathrm{Res}\left(P(s)Q(s),s=n\right)=\\
=&-\frac{1}{2}\left(-3c^2+2c^2\log(c)+\sum_{n=2}^\infty c^{2n}\left(\frac{1}{n}-\frac{2}{n-\frac{1}{2}}+\frac{1}{n-1}\right) \right).
\end{align}
Now, if $c=1$, we can use the digamma function  $\psi(x)=\frac{\Gamma'(x)}{\Gamma(x)}$, which satisfies $\psi(x+1)=\psi(x)+1/x$, to write 
\begin{align}
\sum_{n=2}^N\left(\frac{1}{n}+\frac{1}{n-1}\right)=2\psi(N)+2\gamma-1+\frac{1}{N}\hspace{.4cm}\mbox{and}\hspace{.4cm}\sum_{n=2}^N\frac{1}{n-\frac{1}{2}}=\psi\!\left(N+\frac{1}{2}\right)-\psi\!\left(\frac{3}{2}\right).
\end{align}
Since $\psi\!\left(N+\frac{1}{2}\right)-\psi(N)=\frac{1}{2N}+O(N^{-2})$ as $N\to\infty$ and  $\psi\!\left(\frac{1}{2}\right)=-2\log 2-\gamma$, we have
\begin{align}
D_{\mathrm{rat}}(\chi,\chi)=&-\frac{1}{2}\left(-3+\lim_{N\to\infty}\sum_{n=2}^N \left(\frac{1}{n}-\frac{2}{n-\frac{1}{2}}+\frac{1}{n-1}\right) \right)=\\
=&-\frac{1}{2}\left(-4+2\gamma+2\psi\!\left(\frac{3}{2}\right)\right)=2\log 2.
\end{align}
This proves the first part of \eqref{computation-I(a,b)-statement}.
On the other hand, for $0<c<1$, we can use the series  
\begin{align}
-\log(1 - x)=\sum_{k=1}^\infty \frac{x^k}{k}\hspace{.3cm}\mbox{and}\hspace{.3cm}\mathrm{tanh}^{-1}(x)=\sum_{k=0}^\infty\frac{x^{2k+1}}{2k+1}
\end{align}
to simplify
\begin{align}
&D_{\mathrm{rat}}(\chi,\chi_c)=-\frac{1}{2}\left(-3c^2+2c^2\log(c)+\sum_{n=2}^\infty c^{2n}\left(\frac{1}{n}-\frac{2}{n-\frac{1}{2}}+\frac{1}{n-1}\right) \right)=\\
&=-\frac{1}{2}\left(-3c^2+2c^2\log(c) -c^2-\log(1-c^2)+4c^2-4c\,\mathrm{tanh}^{-1}(c)-c^2\log(1-c^2) \right)=\\
&=\,2c\,\mathrm{tanh}^{-1}(c)-c^2\log(c)+\frac{1+c^2}{2}\log(1-c^2).\label{computation-b-3}
\end{align}
Using \eqref{computation-b-3} with $c=b^{-1}$ and \eqref{computation-b-2}, we obtain  the second part of \eqref{computation-I(a,b)-statement}. 
\end{proof}
\begin{remark}
The function $b\mapsto D_{\mathrm{rat}}(\chi,\chi_b)$ is once but not twice continuously differentiable at $b=1$. Its  first derivative at $1$  equals $2\log 2$, while  its second derivative  equals $\log(1 - b^{-2})$ for $b>1$.
Moreover, $ D_{\mathrm{rat}}(\chi,\chi_b)=\frac{3 + 2 \log b}{2}+O(b^{-2})$ as $b\to\infty$.
\end{remark}

%% file: appendix-zeta-01.tex
\appendix

\section{Bounds for the zeta function}\label{appendix-Zeta}
The constant $C_\eta=2^{6\eta}\zeta(\eta)^2$ in Lemmata \ref{L2.1} and \ref{L2.2_irrational} involves the Riemann zeta function and is used throughout Sections 4, 7, and 8 of the paper. In the proofs of Theorems \ref{thm-rat-epsilon-explicit} and \ref{thm-irr-epsilon-explicit} we use the upper bound \eqref{upper-bound-zeta(eta)-with-sqrt2-constant}, which we prove in this appendix.

%
For $\eta>1$ let $\zeta_{\mathrm{a}}(\eta)=\sum_{n=1}^\infty (-1)^{n-1}n^{-\eta}$ be the alternating zeta function (also known as the Dirichlet eta function).
\begin{lemma}\label{lem-upper-bound-zeta}
For $\eta>1$ we have 
\begin{align}
\zeta_{\mathrm{a}}(\eta)\left(\frac{1}{\log(2)(\eta-1)}+\frac{1}{2}\right)\leq \zeta(\eta)\leq\zeta_{\mathrm{a}}(\eta)\left(\frac{1}{\log(2)(\eta-1)}+\frac{1}{2}+\frac{\log(2)}{12}(\eta-1)\right).\label{statement-lem-upper-bound-zeta}
\end{align}
\end{lemma}
\begin{proof}
The identity $\zeta(\eta)=\zeta_{\mathrm{a}}(\eta)(1-2^{1-\eta})^{-1}$ is classical. Truncating the Laurent series $(1-2^{1-\eta})^{-1}$  near $\eta=1$  to the first order we have $\frac{1}{\log(2)(\eta-1)}+\frac{1}{2}+\frac{\log(2)}{12}(\eta-1)$. This is an upper bound for $(1-2^{1-\eta})^{-1}$ if $\eta>1$. In fact,  an application of l'H\^{o}pital's rule gives
\begin{align}
\lim_{\eta\to1}\frac{1}{\log(2)(\eta-1)}+\frac{1}{2}+\frac{\log(2)}{12}(\eta-1)-(1-2^{1-\eta})^{-1}=0
\end{align}
and we can check that for $\eta>1$
\begin{align}
\frac{\mathrm d}{\mathrm d\eta}\left(\frac{1}{\log(2)(\eta-1)}+\frac{1}{2}+\frac{\log(2)}{12}(\eta-1)-(1-2^{1-\eta})^{-1}\right)>0.
\end{align}
This yields the upper bound in \eqref{statement-lem-upper-bound-zeta}. A similar argument applied to $(1-2^{1-\eta})^{-1}-(\frac{1}{\log(2)(\eta-1)}+\frac{1}{2})$ yields the lower bound. 
\end{proof}
\begin{cor}\label{cor-upper-and-lower-bounds-for-zeta} 
For every $\eta_0>1$ there is a constant $c(\eta_0)>1$ such that if $1<\eta\leq \eta_0$ we have 
$\frac{1}{\eta-1}<\zeta(\eta)\leq \frac{c(\eta_0)}{\eta-1}$. 
\end{cor}
\begin{proof}
Note that functions $\eta\mapsto \zeta_{\mathrm{a}}(\eta)$ and $\eta\mapsto\frac{1}{\log(2)}+\ha(\eta-1)+\frac{\log(2)}{12}(\eta-1)^2$ are both positive and increasing for $\eta>1$. Therefore, for $\eta\in(1,\eta_0]$, their product is maximised at $\eta=\eta_0$. Similarly, since $\frac{1}{\log(2)}+\ha(\eta-1)$ is bounded below by its value at $\eta=1$ and $\zeta_{\mathrm{a}}(1)=\log(2)$ (this is the alternating harmonic series), we have $\zeta_{\mathrm{a}}(\eta)\left(\frac{1}{\log(2)(\eta-1)}+\frac{1}{2}\right)>1$. The statement then follows from Lemma \ref{lem-upper-bound-zeta} with 
\begin{align}
c(\eta_0)= \zeta_{\mathrm{a}}(\eta_0)\left(\frac{1}{\log(2)}+\ha(\eta_0-1)+\frac{\log(2)}{12}(\eta_0-1)^2\right).\label{formula-c(eta0)}
\end{align}
\end{proof}
\begin{remark}
Note that the lower bound in Corollary \ref{cor-upper-and-lower-bounds-for-zeta} is only meaningful for $1<\eta\leq2$ since for $\eta>2$ it is worse than the trivial bound $\zeta(\eta)>1$.
\end{remark}

\begin{cor}\label{cor-upper-and-lower-bounds-for-zeta5/4}
For $1<\eta\leq \frac{5}{4}$, we have $\frac{1}{\eta-1}<\zeta(\eta)\leq \frac{\sqrt{2}}{\eta-1}$.
\end{cor}
\begin{proof}
The statement follows immediately from Corollary \ref{cor-upper-and-lower-bounds-for-zeta} since, using \eqref{formula-c(eta0)}, it is easy to verify that $1.1487793\approx c(\tfrac{5}{4})<\sqrt{2}$.
\end{proof}

Note that Corollary \ref{cor-upper-and-lower-bounds-for-zeta5/4} can be used in the proofs of Theorems \ref{thm-rat-epsilon-explicit} and \ref{thm-irr-epsilon-explicit} since they use the ranges $1<\eta\leq \frac{4+\sqrt{10}}{6}$ and $1<\eta\leq\frac{11+\sqrt{73}}{16}$ respectively and $\frac{4+\sqrt{10}}{6},\frac{11+\sqrt{73}}{16}\leq \frac{5}{4}$.

%% file: tail-bound-MAIN-v9_BUMI.bbl
\begin{thebibliography}{10}

\bibitem{MR3330337}
J.~S. Athreya, J.~Chaika, and S.~Leli\`evre.
\newblock The gap distribution of slopes on the golden {L}.
\newblock In {\em Recent trends in ergodic theory and dynamical systems},
  volume 631 of {\em Contemp. Math.}, pages 47--62. Amer. Math. Soc.,
  Providence, RI, 2015.

\bibitem{Boca-Cobeli-Zaharescu-visible-lattice-points}
F.~P. Boca, C.~Cobeli, and A.~Zaharescu.
\newblock Distribution of lattice points visible from the origin.
\newblock {\em Comm. Math. Phys.}, 213(2):433--470, 2000.

\bibitem{Boca-Heersink-Spiegelhalter}
F.~P. Boca, B.~Heersink, and P.~Spiegelhalter.
\newblock Gap distribution of {F}arey fractions under some divisibility
  constraints.
\newblock {\em Integers}, 13:Paper No. A44, 15, 2013.

\bibitem{MR2274553}
F.~P. Boca and A.~Zaharescu.
\newblock The distribution of the free path lengths in the periodic
  two-dimensional {L}orentz gas in the small-scatterer limit.
\newblock {\em Comm. Math. Phys.}, 269(2):425--471, 2007.

\bibitem{Cellarosi-curlicue}
F.~Cellarosi.
\newblock Limiting curlicue measures for theta sums.
\newblock {\em Ann. Inst. Henri Poincar\'e Probab. Stat.}, 47(2):466--497,
  2011.

\bibitem{Cellarosi-Marklof}
F.~Cellarosi and J.~Marklof.
\newblock Quadratic weyl sums, automorphic functions and invariance principles.
\newblock {\em Proceedings of the London Mathematical Society},
  113(6):775--828, 2016.

\bibitem{Cellarosi-Osman-rational-tails}
F.~Cellarosi and T.~Osman.
\newblock Heavy tailed and compactly supported distributions of quadratic
  {W}eyl sums with rational parameters.
\newblock {\em arXiv:2210.09838}.

\bibitem{DemirciAkarsu}
E.~Demirci~Akarsu.
\newblock Short incomplete {G}auss sums and rational points on metaplectic
  horocycles.
\newblock {\em Int. J. Number Theory}, 10(6):1553--1576, 2014.

\bibitem{DAM2013}
E.~Demirci~Akarsu and J.~Marklof.
\newblock The value distribution of incomplete {G}auss sums.
\newblock {\em Mathematika}, 59(2):381--398, 2013.

\bibitem{MR3177379}
D.~Dolgopyat and B.~Fayad.
\newblock Deviations of ergodic sums for toral translations {I}. {C}onvex
  bodies.
\newblock {\em Geom. Funct. Anal.}, 24(1):85--115, 2014.

\bibitem{MR4179835}
D.~Dolgopyat and B.~Fayad.
\newblock Deviations of ergodic sums for toral translations {II}. {B}oxes.
\newblock {\em Publ. Math. Inst. Hautes \'{E}tudes Sci.}, 132:293--352, 2020.

\bibitem{dudley_2002}
R.~M. Dudley.
\newblock {\em Real Analysis and Probability}.
\newblock Cambridge Studies in Advanced Mathematics. Cambridge University
  Press, 2 edition, 2002.

\bibitem{Elkies-McMullen}
N.~D. Elkies and C.~T. McMullen.
\newblock Gaps in {${\sqrt n}\bmod 1$} and ergodic theory.
\newblock {\em Duke Math. J.}, 123(1):95--139, 2004.

\bibitem{Gut-Probability}
Allan Gut.
\newblock {\em Probability: a graduate course}.
\newblock Springer Texts in Statistics. Springer, New York, second edition,
  2013.

\bibitem{Kowalksi-Sawin}
E.~Kowalski and W.F. Sawin.
\newblock Kloosterman paths and the shape of exponential sums.
\newblock {\em Compos. Math.}, 152(7):1489--1516, 2016.

\bibitem{Lion-Vergne}
G.~Lion and M.~Vergne.
\newblock {\em The {W}eil representation, {M}aslov index and theta series},
  volume~6 of {\em Progress in Mathematics}.
\newblock Birkh\"auser, Boston, Mass., 1980.

\bibitem{Marklof-1999}
J.~Marklof.
\newblock Limit theorems for theta sums.
\newblock {\em Duke Math. J.}, 97(1):127--153, 1999.

\bibitem{Marklof2003b}
J.~Marklof.
\newblock Almost modular functions and the distribution of {$n^2x$} modulo one.
\newblock {\em Int. Math. Res. Not.}, (39):2131--2151, 2003.

\bibitem{Marklof2007b}
J.~Marklof.
\newblock Spectral theta series of operators with periodic bicharacteristic
  flow.
\newblock {\em Ann. Inst. Fourier (Grenoble)}, 57(7):2401--2427, 2007.
\newblock Festival Yves Colin de Verdi{\`e}re.

\bibitem{Marklof-Frobenius-numbers}
J.~Marklof.
\newblock The asymptotic distribution of {F}robenius numbers.
\newblock {\em Invent. Math.}, 181(1):179--207, 2010.

\bibitem{Marklof-Strombergsson-Annals-10}
J.~Marklof and A.~Str{\"o}mbergsson.
\newblock The distribution of free path lengths in the periodic {L}orentz gas
  and related lattice point problems.
\newblock {\em Ann. of Math. (2)}, 172(3):1949--2033, 2010.

\bibitem{Marklof-Strombergsson-Lorentz-gas-asymptotic-estimates}
J.~Marklof and A.~Str\"{o}mbergsson.
\newblock The periodic {L}orentz gas in the {B}oltzmann-{G}rad limit:
  asymptotic estimates.
\newblock {\em Geom. Funct. Anal.}, 21(3):560--647, 2011.

\bibitem{Marklof-Strombergsson-Lorentz-gas-power-law}
J.~Marklof and A.~Str\"{o}mbergsson.
\newblock Power-law distributions for the free path length in {L}orentz gases.
\newblock {\em J. Stat. Phys.}, 155(6):1072--1086, 2014.

\bibitem{Nandori-Szasz-Varju-tail-asymptotics}
P.~N\'{a}ndori, D.~Sz\'{a}sz, and T.~Varj\'{u}.
\newblock Tail asymptotics of free path lengths for the periodic {L}orentz
  process: on {D}ettmann's geometric conjectures.
\newblock {\em Comm. Math. Phys.}, 331(1):111--137, 2014.

\bibitem{Sarnak-equidistribution-horocycles}
P.~Sarnak.
\newblock Asymptotic behavior of periodic orbits of the horocycle flow and
  {E}isenstein series.
\newblock {\em Comm. Pure Appl. Math.}, 34(6):719--739, 1981.

\bibitem{MR3567252}
C.~Uyanik and G.~Work.
\newblock The distribution of gaps for saddle connections on the octagon.
\newblock {\em Int. Math. Res. Not. IMRN}, (18):5569--5602, 2016.

\end{thebibliography}
